\documentclass[intlimits]{article}

\usepackage{latexsym,amsfonts,amssymb,amsmath,amsthm,amscd}

\usepackage{mathtools}
\usepackage{accents}

\usepackage{enumitem}

\usepackage{mdwlist}

\providecommand{\keywords}[1]
{
  \small	
  \textbf{\textit{Keywords.}} #1
}

\providecommand{\MSC}[1]
{
  \small	
  \textbf{\textit{MSC.}} #1
}

\usepackage{bbm}
\usepackage{hyperref}

\usepackage{textgreek}

\DeclareRobustCommand\textprime{\leavevmode\raise.8ex\hbox{\the\scriptfont2 \char48 }}

\theoremstyle{plain}
\newtheorem{theorem}{Theorem}[section]
\newtheorem{lemma}{Lemma}[section]
\newtheorem{proposition}{Proposition}[section]

\newtheorem{definition}{Definition}[section]

\numberwithin{equation}{section}
\numberwithin{definition}{section}

\DeclareFontFamily{U}{mathx}{\hyphenchar\font45}
\DeclareFontShape{U}{mathx}{m}{n}{
      <5> <6> <7> <8> <9> <10>
      <10.95> <12> <14.4> <17.28> <20.74> <24.88>
      mathx10
      }{}
\DeclareSymbolFont{mathx}{U}{mathx}{m}{n}
\DeclareFontSubstitution{U}{mathx}{m}{n}
\DeclareMathAccent{\widecheck}{0}{mathx}{"71}
\DeclareMathAccent{\wideparen}{0}{mathx}{"75}

\newcounter{assumptioncounter}

\def\p{\partial}

\newcommand\JM{Mierczy\'nski}
\newcommand\RR{\ensuremath{\mathbb{R}}}
\newcommand\QQ{\ensuremath{\mathbb{Q}}}

\newcommand\NN{\ensuremath{\mathbb{N}}}

\newcommand{\abs}[1]{\ensuremath{\lvert#1\rvert}}

\newcommand{\norm}[1]{\ensuremath{\lVert#1\rVert}}

\newcommand{\dd}{\ensuremath{\mathrm{d}}}
\newcommand{\DD}{\ensuremath{\mathrm{D}}}

\DeclareMathOperator*{\esssup}{ess\,sup}

\begin{document}
\title {Systems of Parabolic Equations with Delays:  Continuous Dependence on Parameters}

\author{Marek Kryspin
\\
Faculty of Pure and Applied Mathematics \\
Wrocław University of Science and Technology \\
Wybrzeże Wyspiańskiego 27 \\
PL-50-370 Wrocław, Poland \\
\\
and
\\
Janusz Mierczyński
\\
Faculty of Pure and Applied Mathematics \\
Wrocław University of Science and Technology \\
Wybrzeże Wyspiańskiego 27 \\
PL-50-370 Wrocław, Poland \\
}

\date{} \maketitle

\begin{abstract}
    Results on continuous dependence on parameters, as well as on regularization, of solutions to linear systems of parabolic partial differential equations of second order with delay are given.  One of the main features is that the topology on the zero order and delay coefficients is the weak-* topology on an $L_{\infty}$ space, whereas the assumptions on the second- and first-order coefficients are slightly stronger.  This is important in the context of generation of linear skew-product dynamical systems by such equations.
\end{abstract}

\bigskip

\smallskip
\MSC{Primary 35R10, 35B30, 35B65, 46N20.  Secondary 35R05.}

\smallskip
This is an accepted version of the paper published in the Journal of Differential Equations, volume 409, 15 November 2024, pages 532--591, \newline https://doi.org/10.1016/j.jde.2024.07.039.

\smallskip
© 2024. This manuscript version is made available under the CC-BY-NC-ND 4.0 license https://creativecommons.org/licenses/by-nc-nd/4.0/.

\bigskip

\newpage
\section{Introduction}
\label{sect:into0}
In investigating linear nonautonomous evolution (differential) equations considerable role is played by skew\nobreakdash-\hspace{0pt}product dynamical systems.  To put it succinctly, assume that we have an evolution equation (it can be an ordinary differential equation, a partial differential equation, with or without time delay) which can be written as
\begin{equation}
\label{eq:pre-intro-1}
    \frac{\dd u}{\dd t} = \mathcal{A}(t) u,
\end{equation}
where $u$ is an unknown function, defined on some half\nobreakdash-\hspace{0pt}line having $\infty$ as its right endpoint and taking values in some function space endowed with the Banach space structure, and $\{\mathcal{A}(t)\}_{t \in \RR}$ is a given family of linear operators.

Let the above family be so regular that for a given initial condition there exists a unique solution of~\eqref{eq:pre-intro-1} satisfying that condition.  Also, assume that the solution operator has some compactness properties.  We want now to find a dynamical system generated by that solution operator.

The first approach, namely adding one additional dimension corresponding to uniform time flow has a big disadvantage, namely any conceivable compactness of the solution operator is irretrievably lost.  That is the reason why we need to try a different approach.

The idea is to find a metrizable topology on the set of all \emph{time\nobreakdash-\hspace{0pt}translates} of $\{\mathcal{A}(t)\}_{t \in \RR}$.  That topology needs to have the following two, slightly diverging properties:  first, it must not be too fine, in order for the closure of the set of all time\nobreakdash-\hspace{0pt}translates (the so\nobreakdash-\hspace{0pt}called \emph{hull}) to be compact, second, it must not be too coarse, in order for the solution operators to depend continuously, in an appropriate sense, on an element in the hull.

Let us start with the case of $\mathcal{A}(t)$ being an elliptic operator.  Then, if its coefficients depend Hölder\nobreakdash-\hspace{0pt}continuously on time, with the Hölder norm uniformly bounded in $t \in \RR$, one can make use of the theory of solutions based ultimately on Sobolevskiĭ--Tanabe approximations (\cite{Sob, Tan}).  For example, see A. Friedman \cite[Part~2]{Fried}, D. Henry \cite[Chpt.~7]{He}, D. Daners and P. Koch~Medina \cite[Sect.~I.2]{Dan-KM}, or, more generally, H. Amann~\cite[Chpt.~4]{Ama-book}, A. Lunardi~\cite[Chpt.~6]{Lunardi}.

We pass now to a general nonautonomous linear second order partial differential equation (PDE) of parabolic type with delay
\begin{align*}
    \frac{\p u}{\p t} & = \sum_{i, j = 1}^{N}  a_{ij}(t, x) \frac{\p^2 u}{\p x_{i} \p x_{j}}(t,x)  + \sum_{i=1}^{N} b_{i}(t,x) \frac{\p u}{\p x_{i}}(t,x)
    \\
    & + c_0(t,x)u(t,x) + c_1(t,x)u(t-1,x), \quad  0 < t \le T, \ x \in D,
\end{align*}
with second\nobreakdash-\hspace{0pt}order terms satisfying ellipticity conditions and $D \subset \RR^N$ being a bounded domain. When the coefficients are uniformly Hölder\nobreakdash-\hspace{0pt}continuous in $t$, one can apply one of the theories mentioned above to the terms without delay and represent a solution to an initial value problem as the solution of the corresponding integral equation (a \emph{mild solution}, in the terminology introduced by F. E. Browder \cite{Brow}).  For research along such lines, see for~example \cite{Fitz}, \cite{MarSm90}, \cite{NoNuObSa14}, \cite{OS19}.

In many of the papers quoted above, linear theory is usually a springboard to developing a theory of existence of solutions to nonlinear equations.  As a result, function spaces in which solutions are considered must allow some sort of differentiability of the Nemytskiĭ operator, which excludes $L_p$ spaces, for~instance.

\medskip
However, in many cases Hölder dependence of coefficients on time appears to be too strong an assumption.  First, it clearly does not cover the situation, not infrequent in applications, when the coefficients happen to depend on time in a discontinuous way.  Second, even if the original family is sufficiently regular the members of its hull need not be so regular.  In the present paper we put forth a viable theory of continuous dependence of solutions on parameters where the space of parameters is endowed with weak\nobreakdash-\hspace{0pt}* topology in $L_{\infty}$.

We consider systems of linear parabolic partial differential equations with delay
\begin{equation*}
    \begin{aligned}
        \frac{\p u^k}{\p t} = & \sum_{i=1}^{N} \frac{\p}{\p x_i} \biggl( \sum_{j=1}^{N} a^{k}_{ij}(t, x) \frac{\p u^k}{\p x_{j}}(t,x) + a^{k}_{i}(t,x) u^k(t,x) \biggr)
        \\
        & + \sum_{i=1}^N b^{k}_i(t,x)\frac{\p u^k}{\p x_i}(t,x) + \sum_{l=1}^{n} c_0^{kl}(t,x)u^l(t,x)
        \\
        & + \sum_{l=1}^{n} c_1^{kl}(t,x)u^l(t-1,x), \quad 1 \le k \le n, \ s < t \le T, \ x \in D,
    \end{aligned}
\end{equation*}
where $D \subset \RR^N$ is a bounded domain with sufficiently regular boundary $\partial D$.
Equations are complemented with appropriate boundary conditions (assume for~simplicity that they are Dirichlet) and initial conditions over the interval $[s - 1, s]$.

The restrictions of the coefficients of members of the hull to $[0, T] \times D$ serve as parameters in the above initial-boundary value problem.  The space of parameters, $Y$,
is a closed and bounded (in the norm topology) subset of a suitable $L_{\infty}$\nobreakdash-\hspace{0pt}space of (equivalence classes of) functions bounded a.e.\ on $(0, T) \times D$.  The topology on which we are working is the weak\nobreakdash-\hspace{0pt}* topology on $L_{\infty}$ as the dual space to $L_1$.

The assumptions on the second- and first\nobreakdash-\hspace{0pt}order coefficients are a~little stronger:  they are so regular that convergence in the topology of $Y$ entails convergence for a.e. $(t, x) \in (0, T) \times D$, whereas for (coupling) zeroth\nobreakdash-\hspace{0pt}order coefficients as well as time\nobreakdash-\hspace{0pt}delay terms convergence in the weak\nobreakdash-\hspace{0pt}* topology does suffice.

Regarding the solutions, we apply the (now) standard theory of weak solutions, as presented in \cite[Chpt.~XVIII]{DL}, or~\cite[Subsect.~7.1.1--7.1.2]{Evans}, for~example, to the (uncoupled) equations containing only second- and first\nobreakdash-\hspace{0pt}order coefficients.  Solutions of coupled systems of equations with zeroth order and delay terms are given as solutions of appropriate integral equations (mild solutions).

\medskip
One of the most important tools in the analysis of dynamical systems generated by various evolutionary equations are Lyapunov exponents, that is, exponential growth rates of solutions.

From the applications point of view, we are interested in exponential growth rates in the topology of $C([-1, 0], C(\bar{D}))$, or perhaps in $C([-1, 0], C^1(\bar{D}))$.  However, the numerical computation of Lyapunov exponents is much easier in the framework of a Hilbert space.  For~instance, when the dynamical system generated by the equation in~question enjoys some strong order\nobreakdash-\hspace{0pt}preservation property, the top Lyapunov exponent is realized on an invariant one\nobreakdash-\hspace{0pt}dimensional subbundle, and that exponent can be expressed as the integral of some computable function, for a survey see~\cite{Mi-NDS}, and for an example of application see~\cite{COS}.

For parabolic PDEs with delay such a natural Hilbert space is \break $L_2(D) \times L_2((-1, 0), L_2(D))$.  This is a reason why we show that a solution to a system, having started in $L_p(D) \times L_r((-1, 0), L_p(D))$ regularizes, after some waiting time, to a function belonging to $C([-1, 0], L_q(D))$.

\bigskip
The paper is organized as follows.  In Section~\ref{sect:intro} we introduce notation appearing throughout the paper as well as present general results from measure theory and functional analysis frequently needed in the sequel.   In Section~\ref{sect:assumptions} we put forward our assumptions.

Section~\ref{sect:higher-order}, of rather preparatory character, is devoted to properties of solutions to (scalar) linear nonautonomous second order PDEs of parabolic type, in variational (divergence) form, having only terms of second and first order.  The theory of (weak) solutions used by us is the theory initiated in $L_2$\nobreakdash-\hspace{0pt}type spaces, linked with the names of O. A. Ladyzhenskaya~\cite[Chpt.~III]{La}, and J.-L. Lions and E. Magenes~\cite[Chpt.~3]{Lions-M}, see also \cite[Chpt.~XVIII]{DL} and \cite[Chpt.~III]{LaSoUr}\cite[Subsect.~7.1.1--7.1.2]{Evans}.  Later, D. Daners showed in \cite{Dan0}, \cite{Dan2}, \cite{Dan6} that $L_2$\nobreakdash-\hspace{0pt}solutions take, after the initial time, values in more regular $L_p$\nobreakdash-\hspace{0pt}spaces, and gave corresponding estimates.  For later results in that direction, see \cite{Monogr}, \cite{MK-JM}.  It should be noted that, while some results from Section~\ref{sect:higher-order} may be spread in published sources (for example, \cite[p. 195]{He}), this is, to our best knowledge, a first attempt to systematize the material, and it can be of independent interest.

The main part of the paper, Section~\ref{sect:mild}, deals with properties of solutions to linear nonautonomous systems of second order parabolic PDEs with delay.  We have chosen to define a solution to such a system as a mild solution, that is, a solution of the integral equation arising from the Duhamel formula.  In Subsection~\ref{subsect:existence-mild} we give relevant definitions and prove the existence of solutions, defined as continuous functions into $L_p(D, \RR^n)$ under the assumption that the initial condition is in $L_p(D, \RR^n) \times L_r((-1, 0), L_p(D, \RR^n))$.  Subsection~\ref{subsect:estimates-mild} collects and proves results on various estimates used in the sequel.  Results on regularization of solutions are proved in Subsection~\ref{subsect:regularization}.  In Subsection~\ref{subsect:compact} results on compactness are given.  Finally, Subsection~\ref{subsubsect:continuous-dependence}, the main part, concludes the paper with results on continuous dependence.

To give a flavor of the main results, let us formulate an extremely simplified version.
\begin{theorem}
    Consider a system of $n$ linear nonautonomous reaction--diffusion equations with delay terms
    \begin{equation}
    \label{eq:example}
        \begin{aligned}
            \frac{\p u^k}{\p t} = {} & {\Delta}u^k(t,x) + \sum_{l=1}^{n} c_0^{kl}(t,x)u^l(t,x)
            \\
            & + \sum_{l=1}^{n} c_1^{kl}(t,x)u^l(t-1,x), \quad 1 \le k \le n, \ 0 < t \le T, \ x \in D,
    \end{aligned}
    \end{equation}
    where $T > 0$ is sufficiently large, $D \subset \RR^N$ is a bounded domain with boundary $\partial D$, $\Delta$ is the Laplace operator on $D$ with Dirichlet boundary conditions, and the zero\nobreakdash-\hspace{0pt}order coefficients $c_0^{kl}$ as well as delay coefficients $c_1^{kl}$ \textup{(}considered to be parameters\textup{)} belong to a closed and norm-bounded subset $\widetilde{Y}$ of $L_{\infty}((0, T) \times D, \RR^{2n^2})$.  Then
    \begin{itemize}
        \item
        for any $u_0 = (u_0^{(1)}, u_0^{(2)}) \in L_p(D, \RR^n) \times L_r((-1, 0), L_p(D, \RR^n))$, $1 \le p < \infty$ and $1 < r \le \infty$, there exists an $L_p(D, \RR^n)$\nobreakdash-\hspace{0pt}solution $u_p(\cdot)$, defined on $[-1, T]$, of~\eqref{eq:example} satisfying the initial condition
        \begin{equation*}
            \begin{cases}
                u_p(0) = u_0^{(1)}
                \\
                u_p(\tau) = u_0^{(2)}(\tau) & \text{for a.e. } \tau \in [-1, 0);
            \end{cases}
        \end{equation*}
        \item
        starting from some time $T_1 \in (0, T]$ the solution $u_p(\cdot)$ takes values in $L_q(D, \RR^n)$, where $1 \le p < q \le \infty$, where $T_1$ depends on $r$ and the $L_{\infty}$\nobreakdash-\hspace{0pt}norms of the coefficients only; further, the dependence is continuous, with the \break $\mathcal{L}(L_p(D, \RR^n), L_q(D, \RR^n))$\nobreakdash-\hspace{0pt}norms again depending on $r$ and the $L_{\infty}$\nobreakdash-\hspace{0pt}norms of the coefficients only;
        \item
        denote by $u(\cdot; c, u_0)$ the solution of~\eqref{eq:example} corresponding to the parameter $c = ((c_0^{kl})_{k, l = 1}^{n}, (c_1^{kl})_{k, l = 1}^{n}) \in \widetilde{Y}$ and the initial value $u_0 \in L_p(D, \RR^n) \times \newline L_r((-1, 0), L_p(D, \RR^n))$; then there is $T_2 \in [T_1, T)$ such that the mapping
        \begin{multline*}
            \Bigl[\, \widetilde{Y} \times \bigl(L_p(D, \RR^n) \times L_r((-1, 0), L_p(D, \RR^n))\bigr)   \ni (c, u_0)
            \\
            \mapsto u(\cdot; c, u_0){\restriction}_{[T_2, T]} \in C([T_2, T], L_q(D, \RR^n))\,\Bigr]
        \end{multline*}
        is continuous, where $1 \le p \le q < \infty$ and $1 < r \le \infty$.
    \end{itemize}
\end{theorem}
See Proposition~\ref{prop:delay-existence-hull}, Proposition~\ref{prop:regularization-q}, Proposition~\ref{prop:delay-estimates-2}, Proposition~\ref{prop:delay-estimates-infty} and Theorem~\ref{thm:continuity-wrt-parameters-regularizing}.

We need to mention that some of the results have been formulated earlier.  For example, continuity properties of reaction-diffusion equations without delay were considered in~\cite[Chpt.~6]{Monogr}, see also~\cite{MiShPart3}.

There is also vast literature on continuity properties for ordinary differential equations with delay, see, e.g., \cite{MiShPart2}, \cite{NOS3}, \cite{NOS2013}, \cite{COS}, \cite{MiNoOb18}, \cite{MiNoOb20}.

\section{Some notation and basic results}
\label{sect:intro}
We write $\RR^{+}$ for $[0, \infty)$, and $\QQ$ for the set of all rationals.

If $B \subset A$, we write $\mathbbm{1}_{\! B}$ for the \emph{indicator} of $B$: $\mathbbm{1}_{\! B}(a) = 1$ if $a \in B$ and $\mathbbm{1}_{\! B}(a) = 0$ if $a \in A \setminus B$.

For a metric space $(Y, d)$, $\mathfrak{B}(Y)$ denotes the
$\sigma$\nobreakdash-\hspace{0pt}algebra of all Borel subsets of $Y$.

For Banach spaces $X_1$, $X_2$ with norms $\lVert \cdot \rVert_{X_1}$, $\lVert \cdot \rVert_{X_2}$, we let $\mathcal{L}(X_1, X_2)$ stand for the Banach space of bounded linear mappings from $X_1$ into $X_2$, endowed with the standard norm.  Instead~of $\mathcal{L}(X, X)$ we write $\mathcal{L}(X)$.  $\mathcal{L}_{\mathrm{s}}(X_1, X_2)$ denotes the space of bounded linear mappings from $X_1$ into $X_2$ equipped with the strong operator topology. Instead~of $\mathcal{L}_{\mathrm{s}}(X, X)$ we write $\mathcal{L}_{\mathrm{s}}(X)$.

When $X^*$ is the dual of a Banach space $X$, we denote by $\langle \cdot, \cdot \rangle_{X, X^*}$ the duality pairing.  For $X$ being a Hilbert space, $( \cdot, \cdot )_X$ stands for the inner product in $X$.

For a sequence $(u_m)_{m = 1}^{\infty}$ in a Banach space $X$, we write $u_m \to u$ when $u_m$ are (norm) convergent to $u$, and $u_m \rightharpoonup u$ when $u_m$ are weakly convergent to $u$.

\medskip
Throughout the paper, $D \subset \RR^N$ stands for a bounded domain, with boundary $\partial D$.

For $(0, T)$, where $T > 0$, by $\mathfrak{L}((0, T))$ we understand the $\sigma$\nobreakdash-\hspace{0pt}algebra of all Lebesgue\nobreakdash-\hspace{0pt}measurable subsets of $(0, T)$.  The notations $\mathfrak{L}(D)$ and $\mathfrak{L}((0, T) \times D)$ are defined in a similar way.

In the sequel, we will refer to an equivalence class of (measurable, integrable, etc.) functions defined for Lebesgue\nobreakdash-\hspace{0pt}a.e.\ $t \in (0, T)$, etc, simply as a \textit{function}.

For $u$ belonging to a Banach space of functions defined on $D$ we will denote by $u[x]$ the value of $u$ at $x \in D$.  However, occasionally when speaking of differential equations we will use the standard notation $u(t, x)$ for the value of a \emph{solution} $u$ at time $t$ and point $x$.

$L_p(D, \RR^n) =: L_p(D)^n$ has the standard meaning, with the norm, for $1 \le p < \infty$, given by
\begin{equation*}
    \lVert u \rVert_{L_p(D)^n} := \biggl( \sum_{k = 1}^n \int_{D} \lvert u^{k}{[x]} \rvert^{p} \, \mathrm{d}x \biggr)^{\!\frac{1}{p}}, \quad u = (u^1, \dots, u^n),
\end{equation*}
and for $p = \infty$ given by
\begin{equation*}
    \lVert u \rVert_{L_{\infty}(D)^n} := \max\limits_{1 \le k \le n} \esssup\limits_{D} \lvert u^{k}\rvert, \quad u = (u^1, \dots, u^n).
\end{equation*}
Instead of $L_p(D)^1$ we write $L_p(D)$.  There holds
\begin{equation*}
    \lVert u \rVert_{L_p(D)^n} = \biggl( \sum_{k = 1}^n \lVert u^k \rVert_{L_p(D)}^p \biggr)^{\!\frac{1}{p}}, \quad u = (u^1, \dots, u^n),
\end{equation*}
for $1 \le p < \infty$, and
\begin{equation*}
    \lVert u \rVert_{L_{\infty}(D)^n} = \max\limits_{1 \le k \le n} \lVert u^k \rVert_{L_{\infty}(D)}, \quad u = (u^1, \dots, u^n).
\end{equation*}

For $1 \le p \le \infty$ let $p'$ stand for the H\"older conjugate of $p$.  The duality pairing between $L_p(D)^n$ and $L_{p'}(D)^n$ is given, for $1 < p < \infty$, or for $p =1$ and $p' = \infty$, by
\begin{multline*}
    \langle u, v \rangle_{L_p(D)^n, L_{p'}(D)^n} = \sum\limits_{k = 1}^n \int\limits_{D} u^k{[x]} \, v^k{[x]} \, \mathrm{d}x,
    \\
    u = (u^1, \dots, u^n) \in L_p(D)^n, \ v = (v^1, \dots, v^n) \in L_{p'}(D)^n.
\end{multline*}
It is clear that for $1\le p \le q <\infty$
\begin{equation*}
    \langle u, v \rangle_{L_p(D)^n, L_{p'}(D)^n} = \langle u, v \rangle_{L_q(D)^n, L_{q'}(D)^n}, \quad \text{for any } u \in L_q(D)^n, v\in L_{p'}(D)^n.
\end{equation*}
Therefore, the sub indexes in $ \langle \cdot,\cdot \rangle$ can occasionally be suppressed, however, not when clarity is needed.

The symbol $W_{2}^1(D)$ stands for the space of all functions in $L_2(D)$ whose weak partial derivatives of order $1$ belong to $L_2(D)$.  $W_{2}^1(D)$ becomes a Hilbert space when endowed with an inner product
\begin{equation*}
    ( u^k, v^k )_{W_{2}^1(D)} = \int\limits_{D} u^k[x] \, v^k[x] \, \dd x + \sum\limits_{j = 1}^N \int\limits_{D} \DD_j u^k[x] \, \DD_j v^k{[x]} \, \dd x.
\end{equation*}
$\norm{\cdot}_{W_{2}^1(D)}$ denotes the norm on $W_{2}^1(D)$ induced by $( \cdot, \cdot )_{W_{2}^1(D)}$.

$W_{2,0}^{1}(D)$ denotes the closure in $W_{2}^1(D)$ of the set $C^{\infty}_c(D)$ of all $C^{\infty}$ functions having compact support contained in $D$.

Let $u$ be a function defined for Lebesgue\nobreakdash-\hspace{0pt}a.e.\ $t \in (0, T)$ and taking values in $L_p(D)^n$, $1 \le p < \infty$.

\begin{itemize}
    \item
    $u$ is said to be \emph{measurable} if it is $(\mathfrak{L}((0, T)), \mathfrak{B}(L_p(D)^n))$\nobreakdash-\hspace{0pt}measurable, meaning that the preimage under $u$ of any open subset of $L_p(D)^n$ belongs to $\mathfrak{L}((0, T))$.
    \item
    $u$ is \emph{strongly measurable} (sometimes called \emph{Bochner measurable}) if there exists a sequence $(u_m)_{m = 1}^{\infty}$ of simple functions such that $\lim\limits_{m \to \infty} \lVert u_m(t) - u(t) \rVert_{L_p(D)^{n}} = 0$ for Lebesgue\nobreakdash-\hspace{0pt}a.e.\ $t \in (0, T)$.
    \item
    $u$ is \emph{weakly measurable} if for any $v \in L_{p'}(D)^n$ the function
    \begin{equation*}
        [\, t \mapsto \langle u(t), v \rangle_{L_p(D)^n, L_{p'}(D)^n} \,]
    \end{equation*}
    is $(\mathfrak{L}((0, T)), \mathfrak{B}(\RR))$\nobreakdash-\hspace{0pt}measurable.
\end{itemize}
\begin{theorem}
\label{thm:equiv-measurable}
    For $u \colon (0, T) \to L_p(D)^n$, $1 \le p < \infty$, measurability, strong measurability and weak measurability are equivalent.
\end{theorem}
The equivalence of strong and weak measurability is a consequence of Pettis's Measurability Theorem (see, e.g., \cite[Thm.~2.1.2]{DiUhl}).  For the fact that measurability implies strong measurability see, e.g., \cite[Thm.~1]{Vara}, whereas the proof of the reverse implication is a simple exercise.

For our purposes we will use the following definitions (see, e.g.\ \cite[Sect.~X.4]{AmEsch}).  A measurable $u \colon (0, T) \to L_p(D)^n$ belongs to $L_r((0, T), L_p(D)^n)$, $1 \le r \le \infty$, if $\norm{u(\cdot)}_{L_p(D)^n}$ belongs to $L_r((0, T))$, with
\begin{align*}
    \norm{u}_{L_r((0, T), L_p(D)^n)} & := \biggl( \int_{0}^{T} \norm{u(t)}^r_{L_p(D)^n} \, \dd t \biggr)^{\!\! 1/r}, \quad 1 \le r < \infty,
    \\
    \norm{u}_{L_{\infty}((0, T), L_p(D)^n)} & := \esssup_{t \in (0, T)} \norm{u(t)}_{L_p(D)^n}\ .
\end{align*}

The following result, a part of~\cite[Lemma~III.11.16]{DS}, compare~\cite[Thm.~4.4.2]{Vaeth}, will be used several times.
\begin{lemma}~
\label{lm:Dunford-Schwartz}
    \begin{enumerate}
        \item[\textup{(a)}]
        If $u \in L_1((0, T), L_1(D)^n)$ then the function
        \begin{equation*}
            \bigl[\, (0, T) \times D \ni (t, x) \mapsto u(t)[x] \in \RR^n \,\bigr]
        \end{equation*}
        belongs to $L_1((0, T) \times D, \RR^n)$.
        \item[\textup{(b)}]
        If $w$ is $(\mathfrak{L}((0, T) \times D), \mathfrak{B}(\RR^n))$\nobreakdash-\hspace{0pt}measurable, and for Lebesgue\nobreakdash-\hspace{0pt}a.e.\ $t \in (0, T)$ the \emph{$t$\nobreakdash-\hspace{0pt}section}
        $w(t, \cdot)$ belongs to $L_p(D)^n$, where $1 \le p < \infty$, then the function
        \begin{equation*}
            \Bigl[\, (0, T) \ni t \mapsto \bigl[\, D \ni x \mapsto w(t, x) \in \RR^n \,\bigr] \, \Bigr]
        \end{equation*}
        is $(\mathfrak{L}((0, T)), \mathfrak{B}(L_p(D)^n))$\nobreakdash-\hspace{0pt}measurable.
    \end{enumerate}
\end{lemma}

A function $u$ defined for Lebesgue\nobreakdash-\hspace{0pt}a.e.\ $t \in (0, T)$ and taking values in $L_{\infty}(D)^n$ is called \emph{weak\nobreakdash-\hspace{0pt}* measurable} if for each $v \in L_1(D)^n$ the function
\begin{equation*}
        \bigl[\, (0, T) \ni t \mapsto \langle v, u(t) \rangle_{L_1(D)^n, L_{\infty}(D)^n} \,\bigr]
\end{equation*}
is $(\mathfrak{L}((0, T)), \mathfrak{L}(\RR))$\nobreakdash-\hspace{0pt}measurable.

A function $u$ defined for Lebesgue\nobreakdash-\hspace{0pt}a.e.\ $t \in (0, T)$ and taking values in $L_{\infty}(D)^n$ is called \emph{weak\nobreakdash-\hspace{0pt}* integrable} if for each $v \in L_1(D)^n$ the function
\begin{equation*}
        \bigl[\, (0, T) \ni t \mapsto \langle v, u(t) \rangle_{L_1(D)^n, L_{\infty}(D)^n} \,\bigr]
\end{equation*}
belongs to $L_1((0, T))$.  In such a case there exists a unique $w \in L_{\infty}(D)^n$ (the \emph{Gel$\,\textprime$fand integral} of a weak\nobreakdash-\hspace{0pt}* integrable function $u$) such that
\begin{equation*}
   \langle v, w \rangle_{L_1(D)^n,L_{\infty}(D)^n}  = \int_0^T \langle v, u(t) \rangle_{L_1(D)^n,L_{\infty}(D)^n}  \,\dd t
\end{equation*}
for each $v \in L_1(D)^n$.  We will write the Gel$\,\textprime$fand integral of a weak\nobreakdash-\hspace{0pt}* integrable function $u$ as
\begin{equation*}
    \text{G -}\hspace{-1pt}\int_{0}^{T}  u(t) \, \mathrm{d}t.
\end{equation*}
For more on Gel$\textprime$fand integral, see \cite{Hash}.

\section{Assumptions}
\label{sect:assumptions}
Let $T > 0$ be fixed.  We set
\begin{equation*}
   \Delta := \{\, (s, t) \in \RR^2 : 0 \le s \le t \le T \,\}, \qquad \dot{\Delta} := \{\, (s, t) \in \RR^2 : 0 \le s < t \le T \,\}.
\end{equation*}
Observe that in Section~\ref{sect:into0} the symbol $\Delta$ has been employed to denote Laplacian.  As it will not be used in the sequel, this should not cause misunderstandings.

\medskip
Our main object of study are systems of $n$ linear second order partial differential equations with delay
\begin{equation}
\label{main-eq}
\begin{aligned}
\frac{\p u^k}{\p t} = & \sum_{i=1}^{N} \frac{\p}{\p x_i} \biggl( \sum_{j=1}^{N} a^{k}_{ij}(t, x) \frac{\p u^k}{\p x_{j}}(t,x)
+ a^{k}_{i}(t,x) u^k(t,x) \biggr)
\\
& + \sum_{i=1}^N b^{k}_i(t,x)\frac{\p u^k}{\p x_i}(t,x)
\\
& + \sum_{l=1}^{n} c_0^{kl}(t,x)u^l(t,x)
\\
& + \sum_{l=1}^{n} c_1^{kl}(t,x)u^l(t-1,x), \quad 1 \le k \le n, \ s < t \le T, \ x \in D,
\end{aligned}
\end{equation}
where $s \in [0, T)$ is an initial time and $D \subset \RR^N$ is a bounded domain with boundary $\p D$, complemented with boundary conditions
\begin{equation}
\label{main-bc}
\mathcal{B}^k u^k = 0, \qquad 1 \le k \le n, \ s < t \le T, \ x \in \p D,
\end{equation}
where $\mathcal{B}^k$ is either the Dirichlet boundary operator
\begin{equation*}
\mathcal{B}^k u^k = u^k, \quad s < t \le T, \ x \in \partial D
\end{equation*}
or the Robin boundary operator
\begin{equation*}
  \mathcal{B}^k u^k = \sum_{i=1}^N \biggl( \sum_{j=1}^N a_{ij}^{k}(t, x) \frac{\p u^k}{\p x_j} + a^{k}_{i}(t, x) u^{k} \biggr) \nu_i + d^{k}_0(t, x)u^{k}, \quad s < t \le T, \ x \in \partial D
\end{equation*}
(different kinds of boundary conditions are allowed for different $k$).
Above, $\nu = (\nu_1, \dots, \nu_{N})$ denotes the unit normal vector pointing out~of $\p D$.  When $d^{k}_0 \equiv 0$ in the Robin case, $\mathcal{B}^{k} u^k = 0$ is referred to as the Neumann boundary condition.

The initial condition is considered in the following way: for $s \in [0, T)$ and $u_0 = (u_0^{(1)}, u_0^{(2)}) \in L_p(D)^n \times L_r((- 1, 0), L_p(D)^n)$, where $1 \le p \le \infty$, $1 \le  r \le \infty$, find a solution of \eqref{main-eq}+\eqref{main-bc} satisfying
  \begin{equation}
  \label{main-ic}
  \begin{cases}
    u(s) = u_0^{(1)},
    \\
    u(\tau + s) = u_0^{(2)}(\tau) \quad \text{for Lebesgue-a.e. } \tau \in (- 1, 0).
  \end{cases}
  \end{equation}

We introduce some assumptions on $D$ and the coefficients of the problem~\eqref{main-eq}+\eqref{main-bc}.

\medskip
\begin{enumerate}[label=\textup{\textbf{(DA\arabic*)}},ref=\textbf{(DA\arabic*)},align=left]
    \item \label{AS-boundary}  (Boundary regularity)
    \itshape
    $D$ is a bounded domain.  If for some $k$ the boundary conditions are of the Robin type, $D$ has Lipschitz boundary.
\end{enumerate}

\medskip
When speaking of properties satisfied by points in $D$, we use the expression ``for Lebesgue\nobreakdash-\hspace{0pt}a.e.\ $x \in D$'' to indicate that the $N$\nobreakdash-\hspace{0pt}dimensional Lebesgue measure of the set of points not satisfying the property is zero.  Similarly, when speaking of properties satisfied by points in $\p D$, we use the expression ``for Hausdorff\nobreakdash-\hspace{0pt}a.e.\ $x \in \p D$'' to indicate that the $(N-1)$\nobreakdash-\hspace{0pt}dimensional Hausdorff measure of the set of points not satisfying the property is zero.  The expressions ``for Lebesgue\nobreakdash-\hspace{0pt}a.e.\ $(t,x) \in (0, T) \times D$,'' ``for Hausdorff\nobreakdash-\hspace{0pt}a.e.\ $(t,x) \in (0, T) \times \p D$'' are used in an analogous way.  We will skip ``Hausdorff'', ``Lebesgue'' if that cannot lead to misunderstandings.

\medskip
\begin{enumerate}[label=\textup{\textbf{(DA\arabic*)}},ref=\textbf{(DA\arabic*)},align=left,resume]
    \item \label{AS-ess-bounded} (Essential boundedness of coefficients)
    \itshape The functions
    \begin{itemize}[label=$\diamond$]
        \item
        $a^k_{ij} \colon (0, T) \times D \to \RR$ \textup{(}$k = 1, \dots, n$, $i, j = 1, \dots, N$\textup{)},
        \item
        $a^k_{i} \colon (0, T) \times D \to \RR$ \textup{(}$k = 1, \dots, n$, $i = 1, \dots, N$\textup{)},
        \item
        $b^k_{i} \colon (0, T) \times D \to \RR$ \textup{(}$k = 1, \dots, n$, $i = 1, \dots, N$\textup{)},
        \item
        $c_0^{kl} \colon (0, T) \times D \to \RR$ \textup{(}$k, l = 1, \dots, n$\textup{)},
        \item
        $c_1^{kl} \colon (0, T) \times D \to \RR$ \textup{(}$k, l = 1, \dots, n$\textup{)}
    \end{itemize}
    belong to $L_{\infty}((0, T) \times D)$.  For those $k$ for which the Robin boundary conditions hold the functions $d^k_0 \colon (0, T) \times \partial D \to \RR$ belong to $L_{\infty}((0, T) \times \partial D)$.
\end{enumerate}
In the sequel we will assume that $d^k_0 \equiv 0$ for those $k$ for which the Dirichlet boundary conditions hold.

From now on, when speaking of the weak-* topology we will understand $L_{\infty}((0, T) \times D)^{nN^2+2nN+2n^2} \times L_{\infty}((0, T) \times \p D)^{n}$ to be the dual of $L_1((0, T) \times D)^{nN^2+2nN+2n^2} \times L_1((0, T) \times \p D)^{n}$.

We put
\begin{equation*}
    a :=  \bigl( ((a_{ij}^{k})_{i,j=1}^{N})_{k=1}^{n}, ((a_{i}^{k})_{i=1}^{N})_{k=1}^{n}, ((b_{i}^{k})_{i=1}^{N})_{k=1}^{n}, (c_0^{kl})_{k,l=1}^{n}, (c_1^{kl})_{k,l=1}^{n}, (d_0^{k})_{k=1}^{n} \bigr).
\end{equation*}
In view of~\ref{AS-ess-bounded}, $a$ belongs to $L_{\infty}((0, T) \times D)^{nN^2+2nN+2n^2} \times L_{\infty}((0, T) \times \p D)^{n}$.

\smallskip
Denote by $Y$ the set of all possible parameters.  Our next assumption is the following.
\begin{enumerate}[label=\textup{\textbf{(DA\arabic*)}},ref=\textbf{(DA\arabic*)},align=left,resume]
    \item \label{AS-compact} (Compactness)
    \itshape
    $Y$ is a norm bounded subset of $L_{\infty}((0, T) \times D)^{nN^2+2nN+2n^2} \times L_{\infty}((0, T) \times \p D)^{n}$ that is closed in the weak-* topology.

    The function $d_0^k \ge 0$ if the Robin boundary condition holds. The function $d_0$ is interpreted as the zero function in the Dirichlet or Neumann cases.
\end{enumerate}
Unless specifically stated to the contrary, $Y$ is understood with the weak-* topology.  It follows from~\ref{AS-compact} via the Banach--Alaoglu theorem that $Y$ is a compact metrizable space.

\begin{enumerate}[label=\textup{\textbf{(DA\arabic*)}},ref=\textbf{(DA\arabic*)},align=left,resume]
    \item \label{AS-elliptic}  (Ellipticity)
    \itshape
    There exists $\alpha_0 > 0$ such that for each $a \in Y$ there holds
    \begin{equation*}
        \sum\limits_{i,j=1}^{N} a_{ij}^{k}(t, x) \xi_{i} \xi_{j} \ge \alpha_0 \sum\limits_{i=1}^{N} \xi_{i}^2, \qquad k = 1, \dots, n,\quad  \xi = (\xi_1, \dots, \xi_{N}) \in \RR^{N},
    \end{equation*}
    and
    \begin{equation*}
        a_{ij}^{k}(t, x) = a_{ji}^{k}(t, x), \qquad k = 1, \dots, n,\quad i, j = 1, 2, \dots, N,
    \end{equation*}
    for all $(t, x) \in (0, T) \times D$.
\end{enumerate}

\medskip
Throughout, we shall investigate (continuous) dependence of solutions of \eqref{main-eq}+\eqref{main-bc} on initial values and parameters.

\smallskip
For $a \in Y$ we define
\begin{equation*}
  a_{0} := \bigl( ((a_{ij}^{k})_{i,j=1}^{N})_{k=1}^{n}, ((a_{i}^{k})_{i=1}^{N})_{k=1}^{n}, ((b_{i}^{k})_{i=1}^{N})_{k=1}^{n}, (d_{0}^{k})_{k=1}^{n} \bigr),
\end{equation*}
considered an element of $L_{\infty}((0, T)  \times D)^{n (N^2 + 2N)} \times L_{\infty}((0, T)  \times \partial D)^n$.

\smallskip
Put $Y_0 := {\{\, a_0 : a \in Y \,\}}$.  The set $Y_0$ is a compact metrizable space.

For $a_0 \in Y_0$ and $1 \le k \le n$ we write
\begin{equation*}
  a_{0}^k := \bigl( (a_{ij}^{k})_{i,j=1}^{N}, (a_{i}^{k})_{i=1}^{N}, (b_{i}^{k})_{i=1}^{N}, d_{0}^{k} \bigr),
\end{equation*}
considered an element of $L_{\infty}((0, T)  \times D)^{N^2 + 2N} \times L_{\infty}((0, T)  \times \partial D)$.

For $1 \le k \le n$ we put $Y^k_0 := {\{\, a^k_0 : a_0 \in Y_0 \,\}}$.

\begin{enumerate}[label=\textup{\textbf{(DA\arabic*)}},ref=\textbf{(DA\arabic*)},align=left,resume]
    \item \label{AS-ae-converge} (Sequential compactness of $Y_0$ with respect to convergence a.e.)
    \itshape
    Any sequence $(a_{0,m})_{m = 1}^{\infty}$ of elements of $Y_0$, where
    \begin{equation*}
        a_{0,m} := \bigl( ((a_{ij,m}^{k})_{i,j=1}^{N})_{k=1}^{n}, ((a_{i,m}^{k})_{i=1}^{N})_{k=1}^{n}, ((b_{i,m}^{k})_{i=1}^{N})_{k=1}^{n}, (d_{0,m}^{k})_{k=1}^{n} \bigr),
    \end{equation*}
    convergent as $m \to \infty$ in the weak-* topology to $a_0 \in Y_0$ has the property that
        \begin{itemize}
        \item
        $\bigl( ((a_{ij,m}^{k})_{i,j=1}^{N})_{k=1}^{n}, ((a_{i,m}^{k})_{i=1}^{N})_{k=1}^{n}, ((b_{i,m}^{k})_{i=1}^{N})_{k=1}^{n}\bigr)$ converge pointwise a.e.\ on $(0, T) \times D$ to \allowbreak $\bigl( ((a_{ij}^{k})_{i,j=1}^{N})_{k=1}^{n}, ((a_{i}^{k})_{i=1}^{N})_{k=1}^{n}, ((b_{i}^{k})_{i=1}^{N})_{k=1}^{n}\bigr)$;
        \item
        $(d_{0,m}^{k})_{k=1}^{n}$ converge pointwise a.e.\ on $(0, T) \times \partial D$ to $(d_{0}^{k})_{k=1}^{n}$ .
        \end{itemize}
        \setcounter{assumptioncounter}{\value{enumi}}
\end{enumerate}

The assumption \ref{AS-ae-converge} is satisfied, for example, when the assumption \ref{AS-Hoelder} below holds.

\begin{enumerate}[label=\textup{\textbf{(E)}},align=left]
    \item \label{AS-Hoelder}  (Regularity and boundedness of higher order terms)
    \itshape
    \begin{itemize}
        \item
        There exists $\alpha \in (0,1)$ such that $\p D$ is an $(N-1)$\nobreakdash-\hspace{0pt}dimensional manifold of class $C^{3+\alpha}$;
        \item
        for each $a_0 \in Y_0$ the functions $a_{ij}^k$, $a_i^k$ \textup{(}$k = 1, \dots, n$, $i, j = 1, \dots, N$\textup{)} belong to $C^{2+\alpha, 2+\alpha}([0, T] \times \bar{D})$, with their $C^{2+\alpha, 2+\alpha}([0, T] \times \bar{D})$\nobreakdash-\hspace{0pt}norms bounded uniformly in $a_0 \in Y_0$;
        \item
        for each $a_0 \in Y_0$ the functions $b_{i}^{k}$ \textup{(}$k = 1, \dots, n$, $i= 1, \dots, N$\textup{)} belong to $C^{2+\alpha, 1+\alpha}([0, T] \times \bar{D})$, with their $C^{2+\alpha,1+\alpha}([0, T] \times \bar{D})$\nobreakdash-\hspace{0pt}norms bounded uniformly in $a_0 \in Y_0$;
        \item
        for each $a_0 \in Y_0$ the functions $d_{0}^{k}$ belong to $C^{2+\alpha,2+\alpha}([0, T] \times \p D)$, with their $C^{2+\alpha,2+\alpha}([0, T] \times \p D)$\nobreakdash-\hspace{0pt}norms bounded uniformly in $a_0 \in Y_0$.
    \end{itemize}
\end{enumerate}

\section{Equations with second and first order terms only}
\label{sect:higher-order}
We start with $n$ (uncoupled) PDEs, consisting of terms of order two and one only, parameterized by $a_0 \in Y_0$,
\begin{equation}
\label{eq-higher_orders}
    \begin{aligned}
        \frac{\p u^k}{\p t} = & \sum_{i=1}^{N} \frac{\p}{\p x_i} \biggl( \sum_{j=1}^{N} a^{k}_{ij}(t,x) \frac{\p u^k}{\p x_{j}}(t,x) + a^{k}_{i}(t,x) u^k(t,x) \biggr)
        \\
        & {} + \sum_{i=1}^N b^{k}_i(t,x)\frac{\p u^k}{\p x_i}(t,x), \quad 1 \le k \le n, \ s < t \le T, \ x \in D,
    \end{aligned}
\end{equation}
where $s \in [0, T)$ is an initial time, complemented with boundary conditions
\begin{equation}
\label{bc-higher_orders}
    \mathcal{B}^k u^k = 0, \quad 1 \le k \le n, \ s < t \le T, \ x \in \p D,
\end{equation}
where $\mathcal{B}^k$ is either the Dirichlet, or Neumann or else Robin boundary operator.

We will refer to the $k$\nobreakdash-\hspace{0pt}th, $1 \le k \le n$, equation in~\eqref{eq-higher_orders}, as \eqref{eq-higher_orders}$^k$, and similarly for~\eqref{bc-higher_orders}$^k$.

To stress the dependence of the equations and the boundary conditions on the parameter $a_0 \in Y_0$ we will write \eqref{eq-higher_orders}$_{a_0}$+\eqref{bc-higher_orders}$_{a_0}$.  Similarly, we use the notation \eqref{eq-higher_orders}$^k_{a^k_0}$+\eqref{bc-higher_orders}$^k_{a^k_0}$.

\subsection{Weak solutions:  Their existence and basic properties}
\label{subsect:weak}
Throughout the present subsection \ref{AS-boundary} through \ref{AS-elliptic} are assumed to hold.  At some places we assume \ref{AS-ae-converge}, too.

We recall a theory of existence of weak solutions to \eqref{eq-higher_orders}+\eqref{bc-higher_orders} as presented in~\cite{MK-JM}.

Let $1 \le k \le n$ be fixed.

We are looking for solutions of the problem \eqref{eq-higher_orders}$^k_{a^k_0}$+\eqref{bc-higher_orders}$^k_{a^k_0}$ for initial condition $u^k_0 \in L_2(D)$.  To define a solution we introduce the space $H^k$ as follows. Let
\begin{equation*}
V^k
\vcentcolon =
\begin{cases}
    W_{2,0}^1(D) &  \text{if  Dirichlet boundary condition holds}
    \\[1ex]
    W_{2}^1(D) & \text{if Neumann or Robin boundary condition holds}
\end{cases}
\end{equation*}
and
\begin{equation*}
\label{W-space-eq}
H^k = H(s,T;V^k) \vcentcolon= \{\,v \in L_2((s,T),V^k): \dot{v} \in L_2((s,T),(V^k)^{*})\,\}
\end{equation*}
equipped with the norm
\begin{equation*}
\lVert v \rVert_{H^k} \vcentcolon= \Bigl( \int_s^{T} \lVert v(\zeta) \rVert_{V^k}^2 \, \dd\zeta  + \int_s^{T} \lVert \dot{v}(\zeta) \rVert_{(V^k)^{*}}^2 \, \dd\zeta  \Bigr)^{\frac{1}{2}},
\end{equation*}
where $\dot{v} \vcentcolon= \dd v/\dd t$ is the time derivative in the sense of distributions taking values in $(V^k)^{*}$  (see \cite[Chpt.~XVIII]{DL} for definitions).

\medskip

For $a^k_0 \in Y^k_{0}$ define a bilinear form
\begin{equation*}
\begin{aligned}
B^k_{a^k_0}[t; u^k, v^k] & \vcentcolon= \int\limits_{D} \biggl( \, \sum\limits_{i, j = 1}^{N}a^k_{ij}(t, x) \frac{\partial u^k}{\partial x_{i}} \frac{\partial v^k}{\partial x_{j}} + \sum\limits_{i = 1}^{N} a^k_{i}(t, x) u \frac{\partial v^k}{\partial x_{i}}
\\
& - \sum\limits_{i = 1}^{N} b^k_{i}(t, x) \frac{\partial u^k}{\partial x_{i}} v^k \biggr) \, \dd x
\end{aligned}
\end{equation*}
in the Dirichlet or Neumann case, and
\begin{equation*}
\begin{aligned}
B^k_{a^k_0}[t; u^k, v^k] & \vcentcolon= \int\limits_{D} \biggl( \, \sum\limits_{i, j = 1}^{N} a^k_{ij}(t, x) \frac{\partial u^k}{\partial x_{i}} \frac{\partial v^k}{\partial x_{j}} + \sum\limits_{i = 1}^{N} a^k_{i}(t, x) u\frac{\partial v^k}{\partial x_{i}}
\\
& - \sum\limits_{i = 1}^{N} b^k_{i}(t, x) \frac{\partial u^k }{\partial x_{i}} v^k \biggr) \, \dd x + \int\limits_{\partial D} d^k_0(t, x) u^k v^k \, \dd H_{N - 1}
\end{aligned}
\end{equation*}
in the Robin case, where $H_{N - 1}$ stands for the $(N-1)$\nobreakdash-\hspace{0pt}dimensional Hausdorff measure.

\begin{definition}[Weak Solution]
\label{loc-weak-solution-def-null}
For $a^k_0 \in Y^k_{0}$, $0 \le s < T$ and $u^k_0 \in L_2(D)$ a function $u^k \in L_2([ s,T ],V^k)$ such that $\dot{u}^k \in L_2( [ s, T ], (V^k)^{*} )$ is a {\em weak solution\/} of \textup{(\ref{eq-higher_orders})$^k_{a^k_0}$+
(\ref{bc-higher_orders})$^k_{a^k_0}$} {\em on $[s, T]$ with initial condition $u^k(s) = u^k_0$} if
\begin{equation*}
  - \int\limits_{s}^{T} (u^k(\zeta), v^k)_{L_2(D)} \, \dot{\psi}(\zeta) \, \dd\zeta  + \int\limits_{s}^{T} B^k_{a^k_0}[\zeta; u^k(\zeta), v^k]  \psi(\zeta) \, \dd\zeta  = (u^k_0, v^k)_{L_2(D)} \, \psi(s)
\end{equation*}
for any $v^k \in V^k$ and any $\psi \in \mathcal{D}([s, T))$ where set $\mathcal{D}([s, T))$ is the space of all $C^{\infty}$ real functions having compact support in $[s,T)$ and $(\cdot, \cdot)_{L_2(D)}$ denotes the standard inner product in $L_2(D)$.
\end{definition}

\begin{proposition} [Existence of weak solution]
For any $a_0^k \in Y_0^k$ and any initial condition $u^k_0 \in L_2(D)$ there
exists a unique weak solution of~\eqref{eq-higher_orders}$^k_{a^k_0}$\textup{+}\eqref{bc-higher_orders}$^k_{a^k_0}$.
\end{proposition}
\begin{proof}
See \cite[Thm. 2.4]{Dan2} for a proof and \cite[Prop. 2.1.5]{Monogr} for a unified theory of weak solutions.
\end{proof}

For $a^k_0 \in Y_0^k$ and $0 \le s<T$ we write the unique weak solution of \textup{(\ref{eq-higher_orders})$^k_{a^k_0}$+
(\ref{bc-higher_orders})$^k_{a^k_0}$} with initial condition $u^k(s) = u^k_0$ as $u^k(t; a^k_0, u^k_0)$.

Below, the results from Proposition~\ref{prop:higher_orders-skew_product} to Proposition~\ref{prop:local-regularity-higher_order} are suitable reformulations of results in~\cite[Sect.~2]{MK-JM}, in~particular Prop.~2.2 through Prop.~2.8.
\begin{proposition}
  \label{prop:higher_orders-skew_product}
The mappings
  \begin{equation*}
    U^k_{a^k_0}(t,s) u^k_0 \vcentcolon = u^k(t; a^k_0, u^k_0), \quad 0 \le s\le t\le T, \ a^k_0 \in Y^k_0, \ u^k_0 \in L_2(D)
  \end{equation*}
  have the following properties.
  \begin{gather}
  \label{eq:cocycle2-1}
  U^k_{a^k_0}(s, s) = \mathrm{Id}_{L_2(D)}, \quad a^k_0 \in Y^k_0, \ s \in [0,T],
  \\
  \label{eq:cocycle2-2}
  U^k_{a^k_0}(t_2, t_1) \circ U^k_{a^k_0}(t_1, s) = U^k_{a^k_0}(t_2, s), \quad a^k_0 \in Y^k_0, \  0 \le s \le t_1 \le t_2\le T.
  \end{gather}
\end{proposition}

\begin{proposition}
~
  \label{prop:higher_orders-skew_product-p}
   \begin{enumerate}[label=\textup{(\roman*)}, ref=\textup{\ref{prop:higher_orders-skew_product-p}}\textup{(\roman*)}]
     \item\label{prop:higher_orders-skew_product-p_i}
      Let $1 \le p \le \infty$ and $0 \le s < T$.  For any $a^k_0 \in Y^k_0$ and $t \in [s, T]$ there exists $U^k_{a^k_{0},p}(t,s) \in \mathcal{L}(L_p(D))$ such that
      \begin{equation*}
          U^k_{a^k_{0},p}(t,s) u^k_0 = U^k_{a^k_{0}}(t,s) u^k_0, \quad u^k_0 \in L_2(D) \cap L_p(D),
      \end{equation*}
      where $U^k_{a^k_{0},p}(s,s)$ is interpreted as $\mathrm{Id}_{L_p(D)}$.
      \item\label{prop:higher_orders-skew_product-p_ii}
      Let $1 < p < \infty$ and $a^k_{0} \in Y^k_0$.  Then the mapping
      \begin{equation*}
          \bigl[\, [s, T] \ni t \mapsto U^k_{a^k_{0},p}(t,s) \in \mathcal{L}_{\mathrm{s}}(L_p(D)) \,\bigr]
      \end{equation*}
      is continuous.
 \end{enumerate}
\end{proposition}

For $p=1$ we have an analog of Proposition~\ref{prop:higher_orders-skew_product-p_ii}.
\begin{proposition}
\label{prop:higher_orders-p-q-continuity-1}
 Let $1 \le p < \infty$, $0 \le s < T$ and $a^k_0 \in Y_0^{k}$. Then the mapping
    \begin{equation*}
        \big[\, (s, T] \ni t \mapsto U^k_{a^k_0,p}(t, s) \in \mathcal{L}_{\mathrm{s}}(L_p(D)) \, \big]
    \end{equation*}
    is continuous.
\end{proposition}

\begin{proposition}
\label{Extension_identyty_cocycle}
For any $a^k_0 \in Y^k_{0}$, $0 \le s\le t_1 \le t_2 \le T$ and any $1 \le p \le \infty$
\begin{equation}
\label{eq:cocycle2-2_p}
    U^k_{a^k_0,p}(t_2,t_1) \circ U^k_{a^k_0,p} (t_1,s)= U^k_{a^k_0,p}(t_2,s)
\end{equation}
\end{proposition}

\begin{proposition}
  \label{prop:higher_orders-skew_product-p-q-estimates}~
  \begin{enumerate}[label=\textup{(\roman*\textup)},ref=\ref{prop:higher_orders-skew_product-p-q-estimates}\textup{(\roman*)}]
    \item\label{prop:higher_orders-skew_product-p-q-estimates_1} For any $a^k_{0} \in Y^k_0$, any $(s,t) \in \dot{\Delta}$ and any $1 \le p \le q \le  \infty$ there holds $U^k_{a^k_{0},p}(t,s) \in \mathcal{L}(L_p(D), L_q(D))$.
    \item\label{prop:higher_orders-skew_product-p-q-estimates_2} There are constants $M \ge 1$ and $\gamma \in \RR$ such that
    \begin{equation}
    \label{eq:L_p-L_q estymation}
           \bigl\|U^k_{a^k_{0},p}(t, s)\bigr\|_{\mathcal{L}(L_{p}(D), L_{q}(D))} \leq M (t-s)^{-\frac{N}{2}\left(\frac{1}{p}-\frac{1}{q}\right)} e^{\gamma (t-s)}
    \end{equation}
    for $1 \le p \le q \le \infty$, $a^k_{0} \in Y^k_0$ and $(s,t)\in\dot{\Delta}$.
\end{enumerate}
\end{proposition}

In view of the above, from now on we will suppress the subscript $p$ in the notation for $U^k_{a^k_0}(\cdot, \cdot)$.

 In the sequel we will frequently assume that $\gamma \ge 0$ in Proposition~\ref{prop:higher_orders-skew_product-p-q-estimates} and its derivates.  Also, as there finitely many $k$, we assume that $M$ and $\gamma$ are good for all $k$.

\begin{proposition}
\label{prop:local-regularity-higher_order}
Let $1 \le p \le \infty$ and $0 \le s < T$. Then for any $t_0 \in (s, T]$ there exists $\alpha \in (0,1)$ such that for any $a^k_0 \in Y^k_0$, any $u^k_0 \in L_p(D)$, and any compact subset $D_0 \subset D$ the function
$\bigl[\,[t_0,T] \times D_0 \ni (t,x) \mapsto (U^k_{a^k_0}(t) u_0)[x]\,\bigr]$
belongs to $C^{\alpha/2,\alpha}([t_0,T] \times D_0)$. Moreover, for fixed $t_0$, and $D_0$, the $C^{\alpha/2,\alpha}([t_0,T] \times D_0)$\nobreakdash-norm of the above restriction is bounded above by a constant depending on $\norm{u^k_0}_{L_p(D)}$ only.
\end{proposition}

\begin{proposition}
\label{prop:higher_orders-continuity}
    Let $0 \le s < T$, $1 \le p < \infty$ and $a^k_0 \in Y_0^k$.  Then the mapping
    \begin{equation*}
        \bigl[\, (s, T] \ni t \mapsto U^k_{a^k_0}(t,s) \in \mathcal{L}_{\mathrm{s}}(L_p(D))\,\bigr]
    \end{equation*}
    is continuous.
\end{proposition}
For a proof, see~\cite[Prop.~2.2.3]{Monogr}.

\begin{proposition}
\label{prop:higher_orders-continuity-1}
    Let $0 \le s < T$, $1 < p < \infty$ and $a^k_0 \in Y_0^k$.  Then the mapping
    \begin{equation*}
        \bigl[\, [s, T] \ni t \mapsto U^k_{a^k_0}(t,s) \in \mathcal{L}_{\mathrm{s}}(L_p(D))\,\bigr]
    \end{equation*}
    is continuous.
\end{proposition}
For a proof, see~\cite[Thm.~5.1]{Dan2}.

\begin{proposition}
  \label{prop:higher_order-compactness}
  For any $(s,t_0) \in \dot{\Delta}$, $1 \le p < \infty$, $1 \le q < \infty$ and a bounded $E \subset L_p(D)$ the set
  \begin{equation*}
    \bigl\{\, \bigl[\, [t_0, T] \ni t \mapsto U^k_{a^k_0}(t,s) u^k_0 \in L_p(D) \bigr] \vcentcolon a^k_0 \in Y^k_0, \ u^k_0 \in E \,\bigr\}
  \end{equation*}
  is precompact in $C([t_0, T], L_q(D))$.
\end{proposition}
\begin{proof}
  Fix $(s,t_0)\in \dot{\Delta}$, $1 \le p < \infty$, $1 \le q < \infty$ and a bounded $E \subset L_p(D)$. Let $(a^k_{0,m})_{m=1}^{\infty} \subset Y^k_0$ and $(u^k_{0,m})_{m=1}^{\infty} \subset E$.  Put, for $m = 1, 2, \dots$,
  \begin{equation*}
      u^k_m(t) \vcentcolon = U^k_{a^k_{0,m}}(t,s) u^k_{0,m}, \quad t \in [t_0, T].
  \end{equation*}
  It follows from Proposition~\ref{prop:local-regularity-higher_order} via the Arzel\`a--Ascoli theorem by diagonal process that, after possibly taking a subsequence, $(u^k_m)_{m = 1}^{\infty}$ converges as $m \to \infty$ to some function $\tilde{u}^k$ defined on $[t_0, T]$ and taking values in the set of continuous real functions on $D$ in such a way that for any compact $D_0 \subset D$ the functions $[\, t \mapsto u^k_m(t)\!\!\restriction_{D_0} \,]$ converge to $[\, t \mapsto \tilde{u}^k (t)\!\!\restriction_{D_0} \,]$ in $C([t_0, T], C(D_0))$.

  We claim that $u^k_m$ converge to $\tilde{u}^k$ in the $C([t_0, T], L_q(D))$\nobreakdash-\hspace{0pt}norm. By Proposition~\ref{prop:higher_orders-skew_product-p-q-estimates}, there is $M > 0$ such that $\norm{u^k_m(t)}_{L_{\infty}(D)} \le M$ and $\norm{\tilde{u}^k(t)}_{L_{\infty}(D)} \le M$ for all $m = 1, 2, \ldots$ and all $t \in [t_0, T]$.  For $\epsilon > 0$ take a compact $D_0 \subset D$ such that $\lambda(D \setminus D_0) < (\epsilon/(4 M))^q$, where $\lambda$ denotes the $N$\nobreakdash-\hspace{0pt}dimensional Lebesgue measure. We have
  \begin{equation*}
      \norm{(u^k_m(t) - \tilde{u}^k(t)) \, \mathbbm{1}_{\! D \setminus D_0}}_{L_q(D)} \le \frac{\epsilon}{2}
  \end{equation*}
  for all $m = 1, 2, \ldots$ and all $t \in [t_0, T]$.  Further, since $[\, t \mapsto u^k_m(t)\!\!\restriction_{D_0} \,]$ converge to $ [\, t \mapsto \tilde{u}^k(t)\!\!\restriction_{D_0} \,]$ in the $C([t_0, T], C(D_0))$\nobreakdash-\hspace{0pt}norm, there is $m_0$ such that
  \begin{equation*}
      \norm{(u^k_m - \tilde{u}^k) \mathbbm{1}_{\! D_0}}_{C([t_0, T], L_q(D))} \le \frac{\epsilon}{2}
  \end{equation*}
  for all $m \ge m_0$ (here $\mathbbm{1}_{\! D_0}$ stands for the function constantly equal to $\mathbbm{1}_{\! D_0}$).  Consequently,
  \begin{equation*}
      \norm{u^k_m - \tilde{u}^k}_{C([t_0, T], L_q(D))} \le \epsilon
  \end{equation*}
  for all $m \ge m_0$.
 \end{proof}

\begin{proposition}
\label{cor:p-q}
    Let $0 \le s < T$, $1 \le p < \infty$, $1 \le q < \infty$, $a^k_0 \in Y^k_0$.   Then the mapping
    \begin{equation*}
        \bigl[\, (s,T] \times L_p(D) \ni (t, u^k_0) \mapsto U^k_{a^k_0}(t, s) u^k_0 \in L_q(D) \,\bigr]
    \end{equation*}
    is continuous.
\end{proposition}
\begin{proof}
   Let $(u^k_m)_{m = 1}^{\infty}$ converge in $L_p(D)$ to $u^k_0$ and let $(t_m)_{m = 1}^{\infty}$ converge to $t > s$.  Take $\epsilon > 0$.

   We claim that there is $m_1$ such that $\norm{U^k_{a^k_0}(t_m, s) u^k_0 - U^k_{a^k_0}(t, s) u^k_0}_{L_q(D)} < \epsilon/2$ for $m \ge m_1$.  Indeed, if $1 \le q \le p < \infty$ this is a consequence of Proposition~\ref{prop:higher_orders-continuity} and the continuity of the embedding of $L_p(D)$ into $L_q(D)$.  So assume $1 \le p < q < \infty$.  It follows from Proposition~\ref{prop:higher_order-compactness} that, by passing to a subsequence if~necessary, $U^k_{a^k_0}(t_m, s) u^k_0$ converge, as $m \to \infty$, in the $L_q(D)$\nobreakdash-\hspace{0pt}norm, to some $\hat{u} \in L_q(D)$.  As $L_q(D)$ embeds continuously in $L_p(D)$ and, again by Proposition~\ref{prop:higher_orders-continuity}, $U^k_{a^k_0}(t_m, s) u^k_0$ converge, as $m \to \infty$, in the $L_p(D)$\nobreakdash-\hspace{0pt}norm to $U^k_{a^k_0}(t, s) u^k_0$, we have $\hat{u} = U^k_{a^k_0}(t, s) u^k_0$ and the claim holds.

   Further, Proposition~\ref{prop:higher_orders-skew_product-p-q-estimates_2} implies the existence of $m_2$ such that $\norm{U^k_{a^k_0}(t_m, s) u^k_m - U^k_{a^k_0}(t_m, s) u^k_0}_{L_q(D)} < \epsilon/2$ for $m \ge m_2$.  Consequently,
   \begin{multline*}
       \norm{U^k_{a^k_0}(t_m, s) u^k_m - U^k_{a^k_0}(t, s) u^k_0}_{L_q(D)} \\
       \le \norm{U^k_{a^k_0}(t_m, s) u^k_m - U^k_{a^k_0}(t_m, s) u^k_0}_{L_q(D)}
       \\
       + \norm{U^k_{a^k_0}(t_m, s) u^k_0 - U^k_{a^k_0}(t, s) u^k_0}_{L_q(D)} < \epsilon
   \end{multline*}
    for $m \ge  \max\{m_1, m_2\}$.
\end{proof}

\subsection{The Adjoint Operator}
For a fixed $0 < s \le T$, $1 \le k \le N$, together with \textup{(\ref{eq-higher_orders})$_{a^k_0}^{k}$+
(\ref{bc-higher_orders})$_{a^k_0}^{k}$} we consider the adjoint equations, that is the backward parabolic equations
\begin{align}
\label{eq-adjoint-hull}
-\frac{\partial u^k}{\partial t} = &\sum_{i=1}^{N}
\frac{\p }{\p x_i}\Bigl (\sum_{j=1}^{N} a^k_{ji}(t,x) \frac{\p u^k}{\p
x_j} - b^k_i(t,x) u^k \Bigr ) \nonumber \\
& - \sum_{i=1}^{N} a^k_i(t,x) \frac{\p u^k}{\p x_i}, \quad  0 \le t < s, \ x \in D,
\end{align}
complemented with the boundary conditions:
\begin{equation}
\label{bc-adjoint-hull}
\mathcal{B}^{k,*} u^k = 0, \quad  0 \le t < s, \ x \in \p D,
\end{equation}
where $\mathcal{B}^{k,*}$ is as in \eqref{bc-higher_orders} with
$a_0$ replaced by
\begin{equation*}
    a_0^{*} \vcentcolon= (((a^k_{ji})_{i,j=1}^N)_{k=1}^n, (-(b^k_i)_{i=1}^N)_{k=1}^n, (-(a^k_i)_{i=1}^N)_{k=1}^n, (d^k_0)_{k=1}^n).
\end{equation*}
Similarly, we write
\begin{equation*}
    a_0^{k,*} \vcentcolon= ((a^k_{ji})_{i,j=1}^N, -(b^k_i)_{i=1}^N, (a^k_i)_{i=1}^N, d^k_0).
\end{equation*}
To stress the dependence on the parameter we write \eqref{eq-adjoint-hull}$_{a^{k}_0}^{k}$+\eqref{bc-adjoint-hull}$_{a^{k}_0}^{k}$.

Since all analogs of the assumptions \ref{AS-boundary} and \ref{AS-ess-bounded} are satisfied for~\eqref{eq-adjoint-hull}+\eqref{bc-adjoint-hull}, we can define, for $u^k_0 \in L_2(D)$, a (weak) solution of~\eqref{eq-adjoint-hull}$_{a^{k}_0}^{k}$ \allowbreak +\eqref{bc-adjoint-hull}$_{a^{k}_0}^{k}$, defined on $[0,s]$, with the \textit{final condition} $u^k(s) = u^k_0$.  The following analog of Proposition~\ref{prop:higher_orders-skew_product} holds.
\begin{proposition}
    \label{prop:higher_orders-skew_product-adjoint}
  For $a^k_0 \in Y_0^{k}$, $0< s \le T$ and $u^k_0 \in L_2(D)$ there is precisely one weak solution
  \begin{equation*}
   \bigl[\, [0,s] \ni t \mapsto  U^{k,*}_{a^k_0}(t,s) u^k_0 \in L_2(D) \,\bigr]
  \end{equation*}
  of \textup{\eqref{eq-adjoint-hull}$_{a^{k}_0}^{k}$+%
  \eqref{bc-adjoint-hull}}$_{a^{k}_0}^{k}$ satisfying the final condition $u^{*}(s; a^k_0, u^k_0) = u^k_0$. This  mapping has the following properties
\begin{align}
  \label{eq:cocycle2-1-adjoint-semiprocess}
  U^{k,*}_{a^k_{0}}(t, t) = \mathrm{Id}_{L_2(D)}, & \quad a^k_{0} \in Y^k_0, \ t \in [0,s],
  \\
  \label{eq:cocycle2-2-adjoint-semiprocess}
  U^{k,*}_{a^k_{0}}(t_1, t_2) \circ U^{k,*}_{a^k_{0}}(t_2, s) = U^{k,*}_{a^k_{0}}(t_1, s), & \quad a^k_{0} \in Y^k_0, \ 0 \le t_1 \le t_2 \le s.
  \end{align}
\end{proposition}

From now on, in order to connect to  \textup{(\ref{eq-higher_orders})$_{a_0}$+
(\ref{bc-higher_orders})$_{a_0}$} $s$ and $t$ will  be interchanged: $s$ is less than or equal to $t$.

Below we formulate an analog of Proposition~\ref{prop:higher_orders-skew_product-p}.
\begin{proposition}
~
\label{prop:higher_orders-skew_product-p-adjoint}
\begin{enumerate}[label=\textup{(\roman*)},ref=\ref{prop:higher_orders-skew_product-p-adjoint}\textup{(\roman*)}]
      \item\label{prop:higher_orders-skew_product-p-adjoint_(i)}
      Let $1 \le p < \infty$ and $(s,t)\in\dot{\Delta}$. Then $U^{k,*}_{a^k_0}(s,t)$ extends to a linear operator in $\mathcal{L}(L_{p}(D)).$
      \item\label{prop:higher_orders-skew_product-p-adjoint_(ii)}
      Let $1 < p < \infty$, $0 < s \le T$ and $a^k_{0} \in Y^k_0$.  Then the mapping
      \begin{equation*}
           \bigl[\, [0,s] \ni t \mapsto U^{k,*}_{a^k_{0}}(s,t) \in \mathcal{L}_{\mathrm{s}}(L_{p}(D)) \,\bigr]
      \end{equation*}
      is continuous.
  \end{enumerate}
\end{proposition}

The following analog of Proposition~\ref{prop:higher_orders-skew_product-p-q-estimates} holds.
\begin{proposition}
\label{prop:higher_orders-skew_product-p-q-adjoint}~
    \begin{enumerate}[label=\textup{(\roman*\textup)},ref=\ref{prop:higher_orders-skew_product-p-q-adjoint}\textup{(\roman*)}]
    \item\label{prop:higher_orders-skew_product-p-q-adjoint1}
    For any $a^k_{0} \in Y^k_0$, any $0 \le t < s \le T$ and any $1 \le p \le q \le \infty$ there holds $U^{k,*}_{a^k_{0}}(s,t) \in \mathcal{L}_{\mathrm{s}}(L_p(D), L_q(D))$.
    \item\label{prop:higher_orders-skew_product-p-q-adjoint_2}
    There are constants $M \ge 1$ and $\gamma \in \RR$, the same as in Proposition~\textup{\ref{prop:higher_orders-skew_product-p-q-estimates_2}}, such that
    \begin{equation*}
      \lVert U^{k,*}_{a^k_{0}}(s,t) \rVert_{\mathcal{L}(L_p(D), L_q(D))} \le M (t-s)^{-\frac{N}{2} (\frac{1}{p} - \frac{1}{q})} e^{\gamma (t-s)}
    \end{equation*}
    for $1 \le p \le q \le \infty$, $a^k_{0} \in Y^k_{0}$ and $(s,t) \in \dot{\Delta}$.
    \end{enumerate}
\end{proposition}

\begin{proposition}
\label{prop:dual-L2}
    For  $a^k_0 \in Y^k_0$ there holds
    \begin{multline}
    \label{eq:dual-L2-semiprocess}
        ( U^k_{a^k_0}(t,s) u^k_0, v^k_0 )_{L_2(D)} = ( u^k_0, U^{k,*}_{a^k_0}(s,t) v^k_0 )_{L_2(D)}
        \\
        \text{for any } 0\le s \le t \le T, \ u^k_0, \ v^k_0 \in L_{2}(D).
    \end{multline}
\end{proposition}
Proposition~\ref{prop:dual-L2} states that the linear operator $U^{k,*}_{a^k_0}(s, t) \in \mathcal{L}(L_2(D))$ is the dual (in the functional-analytic sense) of $U^k_{a^k_0}(t,s) \in \mathcal{L}(L_2(D))$.  For a proof, see~\cite[Prop.~2.3.3]{Monogr}.

\begin{proposition}
\label{prop:dual-Lp-q}
    For $1 < p \le q < \infty$ and $a^k_0 \in Y^k_0$ there holds
\begin{equation}
    \label{eq:dual-Lp-q-semiprocess}
        \langle U^k_{a^k_0}(t,s) u^k_0, v^k_0 \rangle_{L_q(D), L_{q'}(D)} = \langle u^k_0, U^{k,*}_{a^k_0}(s,t) v^k_0 \rangle_{L_p(D), L_{p'}(D)}
\end{equation}
for any $(s,t) \in \dot{\Delta}$, $u^k_0 \in L_p(D)$ and $v^k_0 \in L_{q'}(D)$.
\end{proposition}
\begin{proof}
     Fix $(s,t) \in \dot{\Delta}$, $u^k_0 \in L_p(D)$ and $v^k_0 \in L_{q'}(D)$.  From Propositions~\ref{prop:higher_orders-skew_product-p-q-estimates_1} and \ref{prop:higher_orders-skew_product-p-q-adjoint} it follows that $U^k_{a^k_0}(\zeta,s) u^k_0, U^{k,*}_{a^k_0}(\zeta,t) v^k_0 \in L_2(D)$ for all $\zeta \in (s, t)$, consequently $\langle U^k_{a^k_0}(\zeta,s) u^k_0, U^k_{a^k_0}(\zeta,t) v^k_0 \rangle_{L_2(D)}$ is well defined for such $\zeta$.
    An application of \eqref{eq:cocycle2-2}, Proposition~\ref{prop:dual-L2} and \eqref{eq:cocycle2-2-adjoint-semiprocess} gives that for any $s < \zeta_1 \le \zeta_2 < t$ there holds
    \begin{multline*}
        ( U^k_{a^k_0}(\zeta_2,s) u^k_0, U^{k,*}_{a^k_0}(\zeta_2,t) v^k_0 )_{L_2(D)}
        \\
        = \langle U^k_{a^k_0}(\zeta_2,\zeta_1) U_{a^k_0}(\zeta_1,s) u^k_0, U^{k,*}_{a^k_0}(\zeta_2,t) v^k_0 \rangle_{L_2(D)}
        \\
        = \langle  U^k_{a^k_0}(\zeta_1,s) u^k_0, U^{k,*}_{a^k_0}(\zeta_1, \zeta_2) U^{k,*}_{a^k_0}(\zeta_2,t) v^k_0 \rangle_{L_2(D)}
        \\
        = ( U^k_{a^k_0}(\zeta_1,s) u^k_0, U^{k,*}_{a^k_0}(\zeta_1,t) v^k_0 )_{L_2(D)}.
    \end{multline*}
    Therefore the assignment
    \begin{align*}
        (s, t) \ni \zeta \mapsto ( & U^k_{a^k_0}(\zeta,s) u^k_0, U^{k,*}_{a^k_0}(\zeta,t) v^k_0 )_{L_2(D)}
        \\
        = {} & \langle U_{a^k_0}(\zeta,s) u^k_0, U^{k,*}_{a^k_0}(\zeta,t) v^k_0 \rangle_{L_q(D), L_{q'}(D)}
        \\
        = {} & \langle U^k_{a^k_0}(\zeta,s) u^k_0, U^{k,*}_{a^k_0}(\zeta,t) v^k_0 \rangle_{L_p(D), L_{p'}(D)}
    \end{align*}
    is constant (denote its value by $A$).  If we let $\zeta \nearrow t$, then $U^k_{a^k_0}(\zeta,s) u^k_0$ converges, by Proposition~\ref{prop:higher_orders-continuity}, in the $L_q(D)$\nobreakdash-\hspace{0pt}norm to $U^k_{a^k_0}(t,s) u^k_0$ and $U^{k,*}_{a^k_0}(\zeta,t) v^k_0$ converges, by Proposition~\ref{prop:higher_orders-skew_product-p-adjoint_(ii)}, in the $L_{q'}(D)$\nobreakdash-\hspace{0pt}norm to $v^k_0$, consequently $\langle U^k_{a^k_0}(t,s) u^k_0, v^k_0 \rangle_{L_q(D), L_{q'}(D)} = A$.  If we let $\zeta \searrow s$, then $U^k_{a^k_0}(\zeta,s) u^k_0$ converges, by Proposition~\ref{prop:higher_orders-skew_product-p_ii}, in the $L_p(D)$\nobreakdash-\hspace{0pt}norm to $u^k_0$ and $U^{k,*}_{a^k_0}(\zeta,t) v_0$ converges, by Propositions~\ref{prop:higher_orders-skew_product-p-q-adjoint} and \ref{prop:higher_orders-skew_product-p-adjoint_(ii)}, in the $L_{p'}(D)$\nobreakdash-\hspace{0pt}norm to $U^{k,*}_{a^k_0}(s,t) v^k_0$, consequently $\langle u^k_0, U^{k,*}_{a^k_0}(s,t) v^k_0 \rangle_{L_p(D), L_{p'}(D)} = A$.  This concludes the proof.
\end{proof}

As $\Delta$ is the closure of $\dot{\Delta}$, the following is a consequence of Propositions~\ref{prop:dual-L2} and~\ref{prop:dual-Lp-q} together with Proposition~\ref{prop:higher_orders-continuity-1} and its adjoint equation version.
\begin{proposition}
\label{prop:dual-Lp}
    For $1 < p < \infty$ and $a^k_0 \in Y^k_0$ there holds
\begin{equation}
    \label{eq:dual-Lpsemiprocess}
        \langle U^k_{a^k_0}(t,s) u^k_0, v^k_0 \rangle_{L_p(D), L_{p'}(D)} = \langle u^k_0, U^{k,*}_{a^k_0}(s,t) v^k_0 \rangle_{L_p(D), L_{p'}(D)}
\end{equation}
for any $(s,t) \in \Delta$, $u^k_0 \in L_p(D)$ and $v^k_0 \in L_{p'}(D)$.
\end{proposition}

It follows from Proposition~\ref{prop:dual-Lp-q} that if $0 \le s < t \le T$ then the linear operator $\allowbreak U^{k,*}_{a^k_0}(s, t) \in \mathcal{L}(L_{q'}(D), L_{p'}(D))$ is the dual (in the functional-analytic sense) of $U^k_{a^k_0}(t,s) \in \mathcal{L}(L_p(D), L_q(D))$.  If $1 < p = q < \infty$, by Proposition~\ref{prop:dual-Lp}, that property holds for $0 \le s \le t \le T$.

\smallskip
We pass now to the case $p = 1$.  For this we will need the following result.
\begin{lemma}
\label{lm:weak-*-aux}
    For $a^k_0 \in Y^k_0$ there holds
    \begin{multline*}
     \langle U^k_{a^k_0}(t,s) u^k_0, v^k_0 \rangle_{L_1(D), L_{\infty}(D)} = ( U^k_{a^k_0}(t,s) u^k_0, v^k_0 )_{L_2(D)}
    \\
    = \langle u^k_0, U^{k,*}_{a^k_0}(s,t) v^k_0 \rangle_{L_1(D), L_{\infty}(D)}
    \end{multline*}
    for any  $(s, t) \in \dot{\Delta}$, $u^k_0 \in L_1(D)$ and $v^k_0 \in L_{\infty}(D)$.
\end{lemma}
\begin{proof}
    Fix $(s, t) \in \dot{\Delta}$, $u^k_0 \in L_1(D)$ and $v^k_0 \in L_{\infty}(D)$.  That all the terms above are well defined follows from Proposition~\ref{prop:higher_orders-skew_product-p-q-estimates_1}, and the first equality is obvious.  Let $(u_{0,m}^k)_{m = 1}^{\infty} \subset L_2(D)$ be a sequence approximating $u_0^k$ in the $L_1(D)$\nobreakdash-\hspace{0pt}norm.  For each $m \in \NN$ the equality $( U^k_{a^k_0}(t,s) u^k_{0,m}, v^k_0 )_{L_2(D)} = (  u^k_{0,m}, U^{k,*}_{a^k_0}(s,t) v^k_0 )_{L_2(D)}$ follows from Proposition~\ref{prop:dual-L2} and the equality $(  u^k_{0,m}, U^{k,*}_{a^k_0}(s,t) v^k_0 )_{L_2(D)} = \break \langle u^k_{0,m}, U^{k,*}_{a^k_0}(s,t) v^k_0 \rangle_{L_1(D), L_{\infty}(D)}$ is straightforward.  As $m \to \infty$, $( U^k_{a^k_0}(t,s) u^k_{0,m}, v^k_0 )_{L_2(D)}$ converge to $( U^k_{a^k_0}(t,s) u^k_0, v^k_0 )_{L_2(D)}$, which is equal to $(  u^k_{0}, U^{k,*}_{a^k_0}(s,t) v^k_0 )_{L_2(D)}$ by Proposition~\ref{prop:dual-L2}.  Finally, $\langle u^k_{0,m}, U^{k,*}_{a^k_0}(s,t) v^k_0 \rangle_{L_1(D), L_{\infty}(D)}$ converge to \newline $\langle u^k_{0}, U^{k,*}_{a^k_0}(s,t) v^k_0 \rangle_{L_1(D), L_{\infty}(D)}$.
\end{proof}
\begin{proposition}
\label{prop:dual-all}
    For $1 \le p \le q < \infty$ and $a^k_0 \in Y^k_0$ there holds
    \begin{equation}
    \label{eq:dual-L1-q}
        \langle U^k_{a^k_0}(t,s) u^k_0, v^k_0 \rangle_{L_q(D), L_{q'}(D)} = \langle u^k_0, U^{k,*}_{a^k_0}(s,t) v^k_0 \rangle_{L_p(D), L_{p'}(D)}
    \end{equation}
    for any $(s,t) \in \dot{\Delta}$, $u^k_0 \in L_p(D)$ and $v^k_0 \in L_{q'}(D)$.
\end{proposition}
\begin{proof}
    For $1 = p = q$ this is a consequence of the equality of the first and last terms in Lemma~\ref{lm:weak-*-aux}, for $1 = p < q < \infty$ this is a consequence of the first equality in Lemma~\ref{lm:weak-*-aux} and Proposition~\ref{prop:dual-Lp-q}, and for $1 < p < \infty$ this is just Proposition~\ref{prop:dual-Lp-q}.
\end{proof}

\subsection{Continuous Dependence of Weak Solutions}
\label{subsect:dependence-weak-solutions}

\begin{lemma}
\label{lm:continuity-in-strong-top}
    Let $1 < p \le q  < \infty$ and $a^k_0 \in Y^k_0$.  Then the mapping
    \begin{equation*}
        \bigl[\, \dot{\Delta} \ni (s, t) \mapsto U^k_{a^k_{0}}(t, s) \in \mathcal{L}_{\mathrm{s}}(L_p(D), L_q(D)) \,\bigr]
    \end{equation*}
    is continuous.
\end{lemma}
\begin{proof}
   Assume that $s_m \to s < t$ as $m \to \infty$.  Fix $u^k_0 \in L_p(D)$ and $v^k_0 \in L_{q'}(D)$.  We have, by Proposition~\ref{prop:dual-all} and the adjoint equation analog of~Proposition~\ref{prop:higher_orders-continuity},
    \begin{multline*}
        \langle U^k_{a^k_{0}}(t, s_m)u^k_0, v^k_0 \rangle_{L_q(D), L_{q'}(D)} = \langle u^k_0, U_{a^k_{0}}^{k,*}(s_m, t) v^k_0 \rangle_{L_p(D), L_{p'}(D)}
        \\
        \to \langle u^k_0,U_{a^k_{0}}^{k,*}(s, t)v_0^k \rangle_{L_p(D), L_{p'}(D)} = \langle U^k_{a^k_{0}}(t, s) u^k_0, v^k_0 \rangle_{L_q(D), L_{q'}(D)},
    \end{multline*}
    so $U^k_{a^k_{0}}(t, s_m) u^k_0 \rightharpoonup U^k_{a^k_{0}}(t, s) u^k_0$ in $L_q(D)$.

    As $\{\, U^k_{a^k_{0}}(t, s_m) u^k_0 : m \in \NN \,\}$ is, by Proposition~\ref{prop:higher_order-compactness}, precompact in $L_q(D)$, the convergence is in the norm.

    Assume that $s_m \to s$ and $t_m \to t$ with $s < t$, and fix $u^k_0 \in L_p(D)$.  We can assume that $s_m < (s + t)/2 < t$ for all $m$.  By the previous paragraph, $U^k_{a^k_{0}}((s + t)/2, s_m) u^k_0 \to U^k_{a^k_{0}}((s + t)/2, s) u^k_0$ in $L_q(D)$.  Proposition~\ref{cor:p-q} implies that
    \begin{gather*}
        U^k_{a^k_{0}}(t_m, s_m) u^k_0 = U^k_{a^k_{0}} (t_m, \tfrac{1}{2} (s + t)) (U^k_{a^k_{0}}(\tfrac{1}{2}(s + t), s_m) u^k_0)
        \\
        \to U^k_{a^k_{0}}(t, \tfrac{1}{2}(s + t)) (U^k_{a^k_{0}}(\tfrac{1}{2}(s + t), s) u^k_0) = U^k_{a^k_{0}}(t, s) u^k_0,
    \end{gather*}
    where the convergence is in $L_q(D)$.
\end{proof}

\begin{lemma}
\label{lm:weak-*-continuity}
    Let $a^k_0 \in Y^k_0$.  For any $v_0^k \in L_{\infty}(D)$ the mapping
    \begin{equation*}
        \bigl[\, \dot{\Delta} \times L_1(D) \ni (s, t, u_0^k) \mapsto \langle U^k_{a^k_{0}}(t, s) u_0^k, v_0^k \rangle_{L_1(D), L_{\infty}(D)} \in \RR \,\bigr]
    \end{equation*}
    is continuous.
\end{lemma}
\begin{proof}
    Fix $v_0^k \in L_{\infty}(D)$. Take sequences $(s_m)_{m = 1}^{\infty}$, $(t_m)_{m = 1}^{\infty}$, $(s_m, t_m) \in \dot{\Delta}$, convergent to $s < t$, and $(u_{0,m}^k)_{m = 1}^{\infty} \subset L_1(D)$ convergent to $u_0^k$.  Let $(\tilde{u}_l)_{l = 1}^{\infty} \subset L_2(D)$ be a sequence approximating $u_0^k$ in the $L_1(D)$\nobreakdash-\hspace{0pt}norm.   By Lemma~\ref{lm:continuity-in-strong-top}, $( U^k_{a^k_{0}}(t_m, s_m) \tilde{u}_l, v_0^k )_{L_2(D)}$ converge to $( U^k_{a^k_{0}}(t, s) \tilde{u}_l, v_0^k )_{L_2(D)}$ as $m \to \infty$, for all $l$, which is, by Lemma~\ref{lm:weak-*-aux}, the same as saying that
    $\langle \tilde{u}_l, U^{k,*}_{a^k_{0}}(s_m,t_m) v^k_0 \rangle_{L_1(D), L_{\infty}(D)}$ converge to \newline $\langle \tilde{u}_l, U^{k,*}_{a^k_0}(s,t) v^k_0 \rangle_{L_1(D), L_{\infty}(D)}$ as $k \to \infty$, for all $l$.  It follows from Proposition~\ref{prop:higher_orders-skew_product-p-q-adjoint_2} that there is $M_1 > 0$ such that $\norm{U^{k,*}_{a^k_{0}}(s_m,t_m) v^k_0}_{L_{\infty}(D)} \le M_1$ for all $m$ and \newline $\norm{U^{k,*}_{a^k_0}(s,t) v^k_0}_{L_{\infty}(D)} \le M_1$.  Take $\epsilon > 0$.  There is $l_1$ such that $\norm{u_0^k - \tilde{u}_{l_1}}_{L_1(D)} < \epsilon/(4 M_1)$.  Let $m_1$ be such that, first, $\norm{u_{0,m}^k - \tilde{u}_{l_1}}_{L_1(D)} < \epsilon/(4 M_1)$ and, second,  $\abs{( U^k_{a^k_{0}}(t_m, s_m) \tilde{u}_{l_1}, v_0^k )_{L_2(D)} - ( U^k_{a^k_{0}}(t, s) \tilde{u}_{l_1}, v_0^k )_{L_2(D)}} < \epsilon/2$, for $m \ge m_1$.  We estimate
    \begin{multline*}
        \abs{\langle U^k_{a^k_{0}}(t_m, s_m) u_{0,m}^k, v_0^k \rangle_{L_1(D), L_{\infty}(D)} - \langle U^k_{a^k_{0}}(t, s) u_0^k, v_0^k \rangle_{L_1(D), L_{\infty}(D)}}
        \\
        \le \abs{\langle u_{0,m}^k - \tilde{u}_{l_1}, U^{k,*}_{a^k_{0}}(s_m,t_m) v^k_0 \rangle_{L_1(D), L_{\infty}(D)}}
        \\
        + \abs{( U^k_{a^k_{0}}(t_m, s_m) \tilde{u}_{l_1}, v_0^k )_{L_2(D)} - ( U^k_{a^k_{0}}(t, s) \tilde{u}_{l_1}, v_0^k )_{L_2(D)}}
        \\
        + \abs{\langle u_0^k - \tilde{u}_{l_1}, U^{k,*}_{a^k_0}(s,t) v^k_0 \rangle_{L_1(D), L_{\infty}(D)}}
        \\
        < \frac{\epsilon}{4 M_1} M_1 + \frac{\epsilon}{2} + \frac{\epsilon}{4 M_1} M_1 = \epsilon
    \end{multline*}
    for all $m \ge m_1$.
\end{proof}

\begin{lemma}
\label{lm:continuity-in-strong-top-L1}
    Let $1 = p \le q < \infty$ and $a^k_0 \in Y^k_0$.  Then the mapping
    \begin{equation*}
        \bigl[\, \dot{\Delta} \ni (s, t) \mapsto U^k_{a^k_{0}}(t, s) \in \mathcal{L}_{\mathrm{s}}(L_1(D), L_q(D)) \,\bigr]
    \end{equation*}
    is continuous.
\end{lemma}
\begin{proof}
   Assume that $s_m \to s < t$ as $m \to \infty$.  Fix $u^k_0 \in L_1(D)$ and $v^k_0 \in L_{q'}(D)$.  We have, by Proposition~\ref{prop:dual-all} and Lemma~\ref{lm:weak-*-continuity},
    \begin{multline*}
        \langle U^k_{a^k_{0}}(t, s_m)u^k_0, v^k_0 \rangle_{L_q(D), L_{q'}(D)} = \langle u^k_0, U_{a^k_{0}}^{k,*}(s_m, t) v^k_0 \rangle_{L_1(D), L_{\infty}(D)}
        \\
        \to \langle u^k_0,U_{a^k_{0}}^{k,*}(s, t)v_0^k \rangle_{L_1(D), L_{\infty}(D)} = \langle U^k_{a^k_{0}}(t, s) u^k_0, v^k_0 \rangle_{L_q(D), L_{q'}(D)},
    \end{multline*}
    so $U^k_{a^k_{0}}(t, s_m) u^k_0 \rightharpoonup U^k_{a^k_{0}}(t, s) u^k_0$ in $L_q(D)$.  As $\{\, U^k_{a^k_{0}}(t, s_m) u^k_0 : m \in \NN \,\}$ is, by Proposition~\ref{prop:higher_order-compactness}, precompact in $L_q(D)$, the convergence is in the norm.

    Assume that $s_m \to s$ and $t_m \to t$ with $s < t$, and fix $u^k_0 \in L_1(D)$.  We can assume that $s_m < (s + t)/2 < t$ for all $m$.  By the previous paragraph, $U^k_{a^k_{0}}((s + t)/2, s_m) u^k_0 \to U^k_{a^k_{0}}((s + t)/2, s) u^k_0$ in $L_q(D)$.  Proposition~\ref{cor:p-q} implies that
    \begin{gather*}
        U^k_{a^k_{0}}(t_m, s_m) u^k_0 = U^k_{a^k_{0}} (t_m, \tfrac{1}{2} (s + t)) (U^k_{a^k_{0}}(\tfrac{1}{2}(s + t), s_m) u^k_0)
        \\
        \to U^k_{a^k_{0}}(t, \tfrac{1}{2}(s + t)) (U^k_{a^k_{0}}(\tfrac{1}{2}(s + t), s) u^k_0) = U^k_{a^k_{0}}(t, s) u^k_0,
    \end{gather*}
    where the convergence is in $L_q(D)$.
\end{proof}
Putting together Lemmas~\ref{lm:continuity-in-strong-top} and~\ref{lm:continuity-in-strong-top-L1} gives the following.
\begin{proposition}
    Let $1 \le p \le q  < \infty$ and $a^k_0 \in Y^k_0$.  Then the mapping
    \begin{equation*}
        \bigl[\, \dot{\Delta} \ni (s, t) \mapsto U^k_{a^k_{0}}(t, s) \in \mathcal{L}_{\mathrm{s}}(L_p(D), L_q(D)) \,\bigr]
    \end{equation*}
    is continuous.
\end{proposition}

\begin{lemma}
\label{lm:continuity-in-strong-top-p}
    Let $1 < p < \infty$ and $a^k_0 \in Y^k_0$.  Then the mapping
    \begin{equation*}
        \big[\, \Delta \ni (s, t) \mapsto U^k_{a^k_{0}}(t, s) \in \mathcal{L}_{\mathrm{s}}(L_p(D)) \,\big]
    \end{equation*}
    is continuous.
\end{lemma}
\begin{proof}
   Assume that $s_m < t$, $s_m \to t$ as $m \to \infty$.  Fix $u^k_0 \in L_p(D)$ and $v^k_0 \in L_{p'}(D)$.  We have, by Proposition~\ref{prop:dual-Lp} and the adjoint equation analog of~Proposition~\ref{prop:higher_orders-continuity-1},
    \begin{multline*}
        \langle U^k_{a^k_{0}}(t, s_m)u^k_0, v^k_0 \rangle_{L_p(D), L_{p'}(D)} = \langle u^k_0, U_{a^k_{0}}^{k,*}(s_m, t) v^k_0 \rangle_{L_p(D), L_{p'}(D)}
        \\
        \to \langle u^k_0, v^k_0 \rangle_{L_p(D), L_{p'}(D)},
    \end{multline*}
    so $U^k_{a^k_{0}}(t, s_m) u^k_0 \rightharpoonup u^k_0$ in $L_p(D)$.  We claim that $\lim\limits_{m \to \infty} \norm{U^k_{a^k_{0}}(t, s_m) u^k_0}_{L_p(D)} = \norm{u^k_0}_{L_p(D)}$.  To prove the claim it suffices, in view of weak convergence, to prove that \newline $\liminf\limits_{m \to \infty}\norm{U^k_{a^k_{0}}(t, s_m) u^k_0}_{L_p(D)} \ge \norm{u^k_0}_{L_p(D)}$.  Suppose to the contrary that there are $\eta > 0$ and a subsequence $(m_l)_{l = 1}^{\infty}$ such that $\norm{U^k_{a^k_{0}}(t, s_{m_l}) u^k_0}_{L_p(D)} \le (1 - \eta) \norm{u^k_0}_{L_p(D)}$.  The Hahn--Banach theorem guarantees the existence of $w \in L_{p'}(D)$ such that $\norm{w}_{L_{p'}(D)} = 1$ and $\langle u^k_0, w \rangle_{L_p(D), L_{p'}(D)} = \norm{u^k_0}_{L_p(D)}$.  Proposition~\ref{prop:higher_orders-skew_product-p-adjoint_(ii)} (with final time $t$) gives that for $l$ sufficiently large one has $\langle u^k_0, U_{a^k_{0}}^{k,*}(s_{m_l}, t) w \rangle_{L_p(D), L_{p'}(D)} > (1 - \tfrac{\eta}{2}) \norm{u^k_0}_{L_p(D)}$, which contradicts \newline $\abs{\langle U^k_{a^k_{0}}(t, s_{m_l})u^k_0, w \rangle_{L_p(D), L_{p'}(D)}} \le \norm{U^k_{a^k_{0}}(t, s_{m_l}) u^k_0}_{L_p(D)} \, \norm{w}_{L_{p'}(D)} \le (1 - \eta) \norm{u^k_0}_{L_p(D)}$ for all $l$.  As $L_p(D)$, being uniformly convex, has the Radon--Riesz property, the norm convergence $U^k_{a^k_{0}}(t, s_m) u^k_0 \to u^k_0$ follows.

    Assume that $s_m < t_m$, $s_m \to t$, $t_m \to t$, and fix $u^k_0 \in L_p(D)$.  The case $s = t < s_m < t_m$ is covered by Proposition~\ref{prop:higher_orders-continuity-1}, so assume that $s_m < t < t_m$ for all $m$.  By the previous paragraph, $U^k_{a^k_{0}}(t, s_m) u^k_0 \to u^k_0$ in $L_p(D)$.  Proposition~\ref{prop:higher_orders-continuity-1} implies that
    \begin{equation*}
        U^k_{a^k_{0}}(t_m, s_m) u^k_0 = U^k_{a^k_{0}} (t_m, t) (U^k_{a^k_{0}}(t, s_m) u^k_0)
        \to u^k_0 \quad \text{in } L_p(D).
    \end{equation*}
    The remaining cases are covered by Lemma~\ref{lm:continuity-in-strong-top}.
\end{proof}

\begin{proposition}
\label{prop:continuity-in-norm-top}
    Let $1 \le p \le q < \infty$ and $a^k_0 \in Y^k_0$. Then the mapping
    \begin{equation*}
       \big[\, \dot{\Delta} \ni (s, t) \mapsto U^k_{a^k_0}(t, s) \in \mathcal{L}(L_p(D), L_q(D)) \,\big]
    \end{equation*}
    is continuous.
\end{proposition}
\begin{proof}
     Let $s_m \to s$ and $t_m \to t$ with $s < t$.  Suppose to the contrary that there are $\epsilon > 0$ and $(u^k_{0,m})_{m = 1}^{\infty} \subset L_p(D)$, $\norm{u^k_{0,m}}_{L_p(D)} = 1$, such that
    \begin{equation*}
        \norm{U^k_{a^k_0}(t_m, s_m) u^k_{0,m} - U_{a^k_0}(t, s) u^k_{0,m}}_{L_p(D)} \ge \epsilon, \quad m = 1, 2, 3, \dots.
    \end{equation*}
    It follows from Proposition~\ref{prop:higher_order-compactness} that, after possibly taking a subsequence and relabelling, we can assume that $U^k_{a^k_0}(t_m, s_m) u^k_{0,m}$ converge to $\tilde{u}$ and $U^k_{a^k_0}(t, s) u^k_{0,m}$ converge to $\hat{u}$, both in $L_q(D)$.  For any $v^k_0 \in L_{q'}(D)$ we have, by Proposition~\ref{prop:dual-all},
    \begin{multline*}
        \langle (U^k_{a^k_0}(t_m, s_m) - U^k_{a^k_0}(t, s)) u^k_{0,m}, v^k_0 \rangle_{L_q(D), L_{q'}(D)}
        \\
        = \langle u^k_{0,m}, (U_{a^k_0}^{k,*}(s_m, t_m) - U_{a^k_0}^{k,*}(s, t)) v^k_0 \rangle_{L_p(D), L_{p'}(D)}.
    \end{multline*}
    Since $\norm{u^k_{0,m}}_{L_p(D)} = 1$, we conclude from the adjoint equation analog of Lemma~\ref{lm:continuity-in-strong-top} (for $1 < p < \infty$) or from Lemma~\ref{lm:weak-*-continuity} (for $p = 1$) that the above expression converges to zero as $m \to \infty$.  Consequently $\tilde{u} = \hat{u}$, a contradiction.
\end{proof}

\begin{proposition}
  \label{prop:higher_orders-skew_product-continuity}
    Assume, in addition, \textup{\ref{AS-ae-converge}}. For $1 < p < \infty$ the mapping
    \begin{equation*}
      \bigl[\, Y^k_0 \times \dot{\Delta} \times  L_p(D) \ni (a^k_{0}, s, t, u^k_0) \mapsto U^k_{a^k_{0}}(t, s) u^k_0 \in L_p(D) \, \bigl]
    \end{equation*}
    is continuous.
\end{proposition}
\begin{proof}
    It follows from~\cite[Props.~2.2.12 and~2.2.13]{Monogr} that, for $2 \le p < \infty$, the mapping
    \begin{equation*}
        \big[\, Y^k_0 \times \dot{\Delta} \times  L_2(D) \ni (a^k_{0}, s, t, u^k_0) \mapsto U^k_{a^k_{0}}(t, s) u^k_0 \in L_p(D) \,\big]
    \end{equation*}
  is continuous.

  In view of the above paragraph and \eqref{eq:cocycle2-2}, to conclude the proof it suffices to show that for any $1 < p < 2$ the mapping
  \begin{equation*}
     \big[\, Y^k_0 \times \dot{\Delta} \times  L_p(D) \ni (a^k_0, s, t, u^k_0) \mapsto U^k_{a^k_0}(t, s) u^k_0 \in L_2(D) \, \big]
    \end{equation*}
    is continuous.  Observe that if we have $a^k_{0,m} \to a^k_0 \in Y^k_0$, $s_m \to s$, $t_m \to t$ with $s_m < t_{m}$ and $s < t$, and $u^k_{0,m} \to u^k_0 \in L_p(D)$, then from Proposition~\ref{prop:higher_order-compactness} it follows that, after possibly choosing a subsequence, there is $w \in L_2(D)$  such that $U^k_{a^k_{0,m}}(t_m, s_m) u^k_{0,m} \to w$ in $L_2(D)$.  Consequently, $( U^k_{a^k_{0,m}}(t_m, s_m) u^k_{0,m}, v_{0}^k )_{L_2(D)} \to ( w, v_0^k )_{L_2(D)}$ as $m \to \infty$, for any $v_0^k \in L_2(D)$.  On the other hand, one has, by Proposition~\ref{prop:dual-Lp-q},
    \begin{equation*}
        ( U^k_{a_{0,m}^{k}}(t_m, s_m) u^k_{0,m}, v_0^k )_{L_2(D)} = \langle  u^k_{0,m}, U^{k,*}_{a^k_{0,m}}(s_m, t_m) v_0^k \rangle_{L_p(D), L_{p'}(D)}.
    \end{equation*}
    As $2 < p' < \infty$, an application of the result in the first paragraph to the adjoint equation yields that $U^{k,*}_{a^k_{0,m}}(s_m,t_m) v^k$ converges, as $m \to \infty$, to $U^{k,*}_{a^k_0}(s,t) v_0^k$ in $L_{p'}(D)$.  As $u^k_{0,m}$ converges to $u^k_0$ in $L_p(D)$, we have that \newline $\langle  u^k_{0,m}, U^{*}_{a^k_{0,m}}(s_m,t_m) v_0^k \rangle_{L_p(D), L_{p'}(D)}$ converges to $\langle  u^k_0, U^{k,*}_{a^k_0}(s,t) v_0^k \rangle_{L_p(D), L_{p'}(D)}$, which is, by Proposition~\ref{prop:dual-Lp-q}, equal to $\langle  U^k_{a^k_0}(t,s) u^k_0,  v_0^k \rangle_{L_2(D)}$.  As $v_0^k \in L_2(D)$ is arbitrary, we have $w = U^k_{a^k_0}(t,s) u^k_0$.
\end{proof}

\begin{proposition}
  \label{prop:continuity-p}
    Assume, in addition, \textup{\ref{AS-ae-converge}}. For $1 < p < \infty$ the mapping
    \begin{equation*}
      \bigl[\, Y^k_0 \times \Delta \times  L_p(D) \ni (a^k_{0}, s, t, u^k_0) \mapsto U^k_{a^k_{0}}(t, s) u^k_0 \in L_p(D) \, \bigl]
    \end{equation*}
    is continuous.
\end{proposition}
\begin{proof}
   In view of Proposition~\ref{prop:higher_orders-skew_product-continuity} it suffices to assume $a^k_{0,m} \to a^k_0 \in Y^k_0$, $s_m \to s$, $t_m \to t$ with $s_m \le t_{m}$ and $s = t$, and $u^k_{0,m} \to u^k_0 \in L_p(D)$.  We estimate
    \begin{multline*}
        \norm{U^k_{a^k_{0,m}}(t_m, s_m) u^k_{0,m} - u^k_{0,m}}_{L_p(D)}
        \\
        \le \norm{U^k_{a^k_{0,m}}(t_m, s_m) (u^k_{0,m} - u^k_{0})}_{L_p(D)} + \norm{U^k_{a^k_{0,m}}(t_m, s_m) u^k_{0} - u^k_0}_{L_p(D)}.
    \end{multline*}
     The first summand on the right\nobreakdash-\hspace{0pt}hand side is estimated, by Proposition~\ref{prop:higher_orders-skew_product-p-q-estimates_2}, by $M e^{\gamma (t_m - s_m)} \norm{u^k_{0,m} - u^k_{0}}_{L_p(D)}$, hence it converges to zero as $m \to \infty$.   The second summand converges to zero by Lemma~\ref{lm:continuity-in-strong-top-p}.
\end{proof}

\begin{proposition}
  \label{prop:higher_orders-skew_product-norm_continuity}
    Assume, in addition, \textup{\ref{AS-ae-converge}}. For $1 < p \le q < \infty$ the mapping
    \begin{equation*}
      \bigl[\, Y^k_{0} \times \dot{\Delta} \ni (a^k_{0}, s, t) \mapsto U^k_{a^k_{0}}(t, s) \in \mathcal{L}(L_p(D), L_q(D)) \, \bigr]
    \end{equation*}
    is continuous.
\end{proposition}
\begin{proof}
    Let $(a^k_{0, m})_{m = 1}^{\infty} \subset Y^k_0$ be a sequence converging to $a^k_0$ as $m \to \infty$, and let $(s_m)_{m = 1}^{\infty}$ and $(t_m)_{m = 1}^{\infty}$ be sequences converging, respectively, to $s$ and $t$ with $0 \le s < t \le T$, as $m \to \infty$.  Suppose to the contrary that $\norm{U^k_{a^k_{0,m}}(t_m, s_m) - U^k_{a^k_{0}}(t, s)}_{\mathcal{L}(L_p(D), L_q(D))}$ does not converge to $0$, that is, there exist $\epsilon > 0$ and a sequence $({u^k_{0,m}})_{m = 1}^{\infty} \subset L_p(D)$, $\norm{u^k_{0,m}}_{L_p(D)} = 1$ for all $m$, such that
    \begin{equation*}
        \norm{U^k_{a^k_{0,m}}(t_m, s_m) u^k_{0,m} - U^k_{a^k_{0}}(t, s) u^k_{0,m}}_{L_q(D)} \ge \epsilon
    \end{equation*}
    for all $m$.

    It follows from Proposition~\ref{prop:higher_order-compactness} that, after possibly extracting a subsequence, we can assume that $U^k_{a^k_{0,m}}((s_m+t_m)/2, s_m) u^k_{0,m}$ and $U^k_{a^k_{0}}((s+t)/2, s) u^k_{0,m}$ converge, as $m \to \infty$, in the $L_q(D)$\nobreakdash-\hspace{0pt}norm. We claim that both converge to the same $\tilde{u}$.  Indeed, it suffices to check that the difference $(U^k_{a^k_{0,m}}((s_m+t_m)/2, s_m) - U^k_{a^k_{0}}((s+t)/2, s)) u^k_{0,m}$ converges to zero in $L_q(D)$, which is, in~light of the equalities (Proposition~{\ref{prop:dual-all}})
    \begin{multline*}
        \langle (U^k_{a^k_{0,m}}(s_m+t_m)/2, s) - U^k_{a^k_{0}}((s+t)/2, s)) u^k_{0,m}, v_{0}^k \rangle_{L_q(D), L_{q'}(D)}
        \\
        = \langle u^k_{0,m}, (U^{k,*}_{a^k_{0,m}}(s_m, (s_m+t_m)/2) - U^{k,*}_{a^k_{0}}(s, (s+t)/2)) v_{0}^k \rangle_{L_p(D), L_{p'}(D)}, \quad v_{0}^k \in L_{p'}(D),
    \end{multline*}
    a consequence of the analog for the adjoint equation of Proposition~\ref{prop:higher_orders-skew_product-continuity}.

    Proposition~\ref{prop:higher_orders-skew_product-continuity} implies that
    \begin{multline*}
        \norm{U^k_{a^k_{0,m}}((t_m, s_m) u^k_{0,m} - U^k_{a^k_{0}}(t, (s+t)/2) \tilde{u}}_{L_q(D)}
        \\
        = \norm{U^k_{a^k_{0,m}}(t_m, (s_m+t_m)/2) (U^k_{a^k_{0,m}}((s_m+t_m)/2, s_m) u^k_{0,m})- U^k_{a^k_{0}}(t, (s+t)/2) \tilde{u}}_{L_q(D)} \to 0,
    \end{multline*}
    and
    \begin{multline*}
        \norm{U^k_{a^k_{0}}(t, s) u^k_{0,m} - U^k_{a^k_{0}}(t, (s+t)/2)) \tilde{u}}_{L_q(D)}
        \\
        = \norm{U^k_{a^k_{0}}(t, (s+t)/2) (U^k_{a^k_{0}}((s+t)/2, s) u^k_{0,m}) - U^k_{a^k_{0}}(t, (s+t)/2) \tilde{u}}_{L_q(D)}
        \\
        \le \norm{U^k_{a^k_{0}}(t, (s+t)/2)}_{\mathcal{L}(L_q(D))} \norm{U^k_{a^k_{0}}((s+t)/2, s) u^k_{0,m} - \tilde{u}}_{L_q(D)} \to 0
      \end{multline*}
    therefore $\norm{U^k_{a^k_{0,m}}(t_m, s_m) u^k_{0,m} - U^k_{a^k_{0}}(t, s) u^k_{0,m}}_{L_q(D)}$ converges to zero, a contradiction.
\end{proof}

We proceed now to consider the case $p = 1$.  The next result is an extension of Lemma~\ref{lm:weak-*-continuity}.
\begin{lemma}
\label{lm:weak-*-continuity-parameters}
    Assume, in addition, \textup{\ref{AS-ae-converge}}.  For any $v_0^k \in L_{\infty}(D)$ the mapping
    \begin{equation*}
        \bigl[\, Y_0^k \times \dot{\Delta} \times L_1(D) \ni (a_0^k, s, t, u_0^k) \mapsto \langle U^k_{a^k_{0}}(t, s) u_0^k, v_0^k \rangle_{L_1(D), L_{\infty}(D)} \in \RR \,\bigr]
    \end{equation*}
    is continuous.
\end{lemma}
\begin{proof}
    Fix $v_0^k \in L_{\infty}(D)$. Take sequences $(a_{0, m}^k)_{m = 1}^{\infty} \subset Y_0^k$ convergent to $a_{0}^k$, $(s_m)_{m = 1}^{\infty}$, $(t_m)_{m = 1}^{\infty}$, $(s_m, t_m) \in \dot{\Delta}$, convergent to $s < t$, and $(u_{0,m}^k)_{m = 1}^{\infty} \subset L_1(D)$ convergent to $u_0^k$.  Let $(\tilde{u}_l)_{l = 1}^{\infty} \subset L_2(D)$ be a sequence approximating $u_0^k$ in the $L_1(D)$\nobreakdash-\hspace{0pt}norm.   By Proposition~\ref{prop:higher_orders-skew_product-norm_continuity}, $( U^k_{a^k_{0,m}}(t_m, s_m) \tilde{u}_l, v_0^k )_{L_2(D)}$ converge to $( U^k_{a^k_{0}}(t, s) \tilde{u}_l, v_0^k )_{L_2(D)}$ as $m \to \infty$, for all $l$, which is, by Lemma~\ref{lm:weak-*-aux}, the same as saying that
    $\langle \tilde{u}_l, U^{k,*}_{a^k_{0,m}}(s_m,t_m) v^k_0 \rangle_{L_1(D), L_{\infty}(D)}$ converge to \newline $\langle \tilde{u}_l, U^{k,*}_{a^k_0}(s,t) v^k_0 \rangle_{L_1(D), L_{\infty}(D)}$ as $k \to \infty$, for all $l$.  It follows from Proposition~\ref{prop:higher_orders-skew_product-p-q-adjoint_2} that there is $M_1 > 0$ such that $\norm{U^{k,*}_{a^k_{0,m}}(s_m,t_m) v^k_0}_{L_{\infty}(D)} \le M_1$ for all $m$ and \newline $\norm{U^{k,*}_{a^k_0}(s,t) v^k_0}_{L_{\infty}(D)} \le M_1$.  Take $\epsilon > 0$.  There is $l_1$ such that $\norm{u_0^k - \tilde{u}_{l_1}}_{L_1(D)} < \epsilon/(4 M_1)$.  Let $m_1$ be such that, first, $\norm{u_{0,m}^k - \tilde{u}_{l_1}}_{L_1(D)} < \epsilon/(4 M_1)$ and, second,  $\abs{( U^k_{a^k_{0,m}}(t_m, s_m) \tilde{u}_{l_1}, v_0^k )_{L_2(D)} - ( U^k_{a^k_{0}}(t, s) \tilde{u}_{l_1}, v_0^k )_{L_2(D)}} < \epsilon/2$, for $m \ge m_1$.  We estimate
    \begin{multline*}
        \abs{\langle U^k_{a^k_{0,m}}(t_m, s_m) u_{0,m}^k, v_0^k \rangle_{L_1(D), L_{\infty}(D)} - \langle U^k_{a^k_{0}}(t, s) u_0^k, v_0^k \rangle_{L_1(D), L_{\infty}(D)}}
        \\
        \le \abs{\langle u_{0,m}^k - \tilde{u}_{l_1}, U^{k,*}_{a^k_{0,m}}(s_m,t_m) v^k_0 \rangle_{L_1(D), L_{\infty}(D)}}
        \\
        + \abs{( U^k_{a^k_{0,m}}(t_m, s_m) \tilde{u}_{l_1}, v_0^k )_{L_2(D)} - ( U^k_{a^k_{0}}(t, s) \tilde{u}_{l_1}, v_0^k )_{L_2(D)}}
        \\
        + \abs{\langle u_0^k - \tilde{u}_{l_1}, U^{k,*}_{a^k_0}(s,t) v^k_0 \rangle_{L_1(D), L_{\infty}(D)}}
        \\
        < \frac{\epsilon}{4 M_1} M_1 + \frac{\epsilon}{2} + \frac{\epsilon}{4 M_1} M_1 = \epsilon
    \end{multline*}
    for all $m \ge m_1$.
\end{proof}

\begin{lemma}
  \label{lm:L1-L2-norm_continuity}
    Assume, in addition, \textup{\ref{AS-ae-converge}}. The mapping
    \begin{equation*}
      \bigl[\, Y^k_{0} \times \dot{\Delta} \ni (a^k_{0}, s, t) \mapsto U^k_{a^k_{0}}(t, s) \in \mathcal{L}(L_1(D), L_2(D)) \, \bigr]
    \end{equation*}
    is continuous.
\end{lemma}
\begin{proof}
    Let $(a^k_{0, m})_{m = 1}^{\infty} \subset Y^k_0$ be a sequence converging to $a^k_0$ as $m \to \infty$, and let $(s_m)_{m = 1}^{\infty}$ and $(t_m)_{m = 1}^{\infty}$ be sequences converging, respectively, to $s$ and $t$ with $0 \le s < t \le T$, as $m \to \infty$.  Suppose to the contrary that $\norm{U^k_{a^k_{0,m}}(t_m, s_m) - U^k_{a^k_{0}}(t, s)}_{\mathcal{L}(L_1(D), L_2(D))}$ does not converge to $0$, that is, there exist $\epsilon > 0$ and a sequence $({u^k_{0,m}})_{m = 1}^{\infty} \subset L_1(D)$, $\norm{u^k_{0,m}}_{L_1(D)} = 1$ for all $m$, such that
    \begin{equation*}
        \norm{U^k_{a^k_{0,m}}(t_m, s_m) u^k_{0,m} - U^k_{a^k_{0}}(t, s) u^k_{0,m}}_{L_2(D)} \ge \epsilon
    \end{equation*}
    for all $m$.

    It follows from Proposition~\ref{prop:higher_order-compactness} that, after possibly extracting a subsequence, we can assume that $U^k_{a^k_{0,m}}((s_m+t_m)/2, s_m) u^k_{0,m}$ and $U^k_{a^k_{0}}((s+t)/2, s) u^k_{0,m}$ converge, as $m \to \infty$, in the $L_q(D)$\nobreakdash-\hspace{0pt}norm. We claim that both converge to the same $\tilde{u}$.  Indeed, it suffices to check that the difference $(U^k_{a^k_{0,m}}((s_m+t_m)/2, s_m) - U^k_{a^k_{0}}((s+t)/2, s)) u^k_{0,m}$ converges to zero in $L_2(D)$, which is, in~light of the equalities (Proposition~{\ref{prop:dual-all}})
    \begin{multline*}
        \bigl( (U^k_{a^k_{0,m}}(s_m+t_m)/2, s) - U^k_{a^k_{0}}((s+t)/2, s)) u^k_{0,m}, v_{0}^k \bigr)_{L_2(D)}
        \\
        = \langle u^k_{0,m}, (U^{k,*}_{a^k_{0,m}}(s_m, (s_m+t_m)/2) - U^{k,*}_{a^k_{0}}(s, (s+t)/2)) v_{0}^k \rangle_{L_1(D), L_{\infty}(D)}, \quad v_{0}^k \in L_{\infty}(D),
    \end{multline*}
    a consequence of Lemma~\ref{lm:weak-*-continuity-parameters}.

    Proposition~\ref{prop:higher_orders-skew_product-continuity} implies that
    \begin{multline*}
        \norm{U^k_{a^k_{0,m}}((t_m, s_m) u^k_{0,m} - U^k_{a^k_{0}}(t, (s+t)/2) \tilde{u}}_{L_2(D)}
        \\
        = \norm{U^k_{a^k_{0,m}}(t_m, (s_m+t_m)/2) (U^k_{a^k_{0,m}}((s_m+t_m)/2, s_m) u^k_{0,m})- U^k_{a^k_{0}}(t, (s+t)/2) \tilde{u}}_{L_2(D)} \to 0,
    \end{multline*}
    and
    \begin{multline*}
        \norm{U^k_{a^k_{0}}(t, s) u^k_{0,m} - U^k_{a^k_{0}}(t, (s+t)/2)) \tilde{u}}_{L_2(D)}
        \\
        = \norm{U^k_{a^k_{0}}(t, (s+t)/2) (U^k_{a^k_{0}}((s+t)/2, s) u^k_{0,m}) - U^k_{a^k_{0}}(t, (s+t)/2) \tilde{u}}_{L_2(D)}
        \\
        \le \norm{U^k_{a^k_{0}}(t, (s+t)/2)}_{\mathcal{L}(L_q(D))} \norm{U^k_{a^k_{0}}((s+t)/2, s) u^k_{0,m} - \tilde{u}}_{L_2(D)} \to 0
      \end{multline*}
    therefore $\norm{U^k_{a^k_{0,m}}(t_m, s_m) u^k_{0,m} - U^k_{a^k_{0}}(t, s) u^k_{0,m}}_{L_2(D)}$, converges to zero, a contradiction.
\end{proof}
A consequence of Proposition~\ref{prop:higher_orders-skew_product-norm_continuity} and Lemma~\ref{lm:L1-L2-norm_continuity} together with Proposition~\ref{Extension_identyty_cocycle} is the following.
\begin{proposition}
  \label{prop:norm_continuity-all}
    Assume, in addition, \textup{\ref{AS-ae-converge}}. For $1 \le p \le q < \infty$ the mapping
    \begin{equation*}
      \bigl[\, Y^k_{0} \times \dot{\Delta} \ni (a^k_{0}, s, t) \mapsto U^k_{a^k_{0}}(t, s) \in \mathcal{L}(L_p(D), L_q(D)) \, \bigr]
    \end{equation*}
    is continuous.
\end{proposition}

\section{Mild Solutions}
\label{sect:mild}

After all those preparations, we proceed now to the main part of the paper.  The solutions of~\eqref{main-eq}+\eqref{main-bc} are defined as solutions of a perturbation, via the Duhamel formula, of an uncoupled system having only second- and first-order terms.
In Subsection~\ref{subsect:existence-mild} the existence of such solutions is proved.  Useful estimates are proved in Subsection~\ref{subsect:estimates-mild}.  Subsection~{\ref{subsect:regularization}} presents results on regularization of solutions, including some further estimates. In Subsection~\ref{subsect:compact} results on compactness are given.  The last part, Subsection~\ref{subsubsect:continuous-dependence}, concludes the paper with results on continuous dependence on initial values and parameters.

\smallskip
We assume \ref{AS-boundary}--\ref{AS-elliptic}.
\subsection{Existence of mild solutions}
\label{subsect:existence-mild}
We will look for solutions of~\eqref{main-eq}+\eqref{main-bc} as solutions of some integral equation.  In order to do so we introduce some notation.

For $a_0 \in Y_0$ and $(s, t) \in \dot{\Delta}$ we write
\begin{equation*}
  U^0_{a_0}(t, s) := (U^1_{a^1_0}(t, s), \dots, U^n_{a^n_0}(t, s)).
\end{equation*}
%In other words, for $u_0 = (u^1_0, \dots, u^n_0)$ the symbol $\widetilde{U}^0_{\tilde{a}_0}(t, s) u_0$ denotes
%\begin{equation*}
%  (\tilde{u}^1(t; s, \tilde{a}_0, u^1_0), \dots, \tilde{u}^n(t; s, \tilde{a}_0, u^n_0)).
%\end{equation*}
Depending on the context, $U^0_{a_0}(t, s)$ is considered an element of $\mathcal{L}(L_p(D)^n, L_q(D)^n)$, etc.

For further reference we formulate the following corollary of Proposition~\ref{prop:higher_orders-skew_product-p-q-estimates}.
\begin{proposition}~
\label{prop:system-estimates}
   \begin{enumerate}
    \item[\textup{(i)}]
    For any $a_0 \in Y_0$, any $(s, t) \in \dot{\Delta}$ and any $1 \le p \le q \le \infty$ there holds $U^{0}_{a_0}(t, s) \in \mathcal{L}(L_p(D)^n, L_q(D)^n)$.
    \item[\textup{(ii)}]
    There are constants $M \ge 1$ and $\gamma \in \RR$, the same as in Proposition~\ref{prop:higher_orders-skew_product-p-q-estimates}(ii), such that
    \begin{equation*}
      \lVert U^{0}_{a_0}(t, s) \rVert_{\mathcal{L}(L_p(D)^n, L_q(D)^n)} \le M (t - s)^{-\tfrac{N}{2} \bigl(\tfrac{1}{p} - \tfrac{1}{q}\bigr)} e^{\gamma (t - s)}
    \end{equation*}
    for $1 \le p \le q \le \infty$, $a_0 \in Y_0$ and $(s, t) \in \dot{\Delta}$.
  \end{enumerate}
\end{proposition}
\begin{proof}
    Part (i) is straightforward.  To prove (ii), assume first $1 < p \le q < \infty$.  Observe that, for $u_0 = (u_0^1, \ldots, u_0^n)$, there holds
    \begin{align*}
        \lVert U^{0}_{a_0}(t, s) u_0 \rVert_{L_q(D)^n} & = \biggl( \sum\limits_{k=1}^{n} \lVert U^{k}_{a_0^k}(t, s) u^k_0 \rVert^q_{L_q(D)} \biggr)^{\!\frac{1}{q}}
        \\[1ex]
        & \le M (t - s)^{-\tfrac{N}{2} \bigl(\tfrac{1}{p} - \tfrac{1}{q}\bigr)} e^{\gamma (t - s)} \,
        {\biggl( \sum\limits_{k=1}^{n} \lVert u^k_0 \rVert^q_{L_p(D)} \biggr)^{\!\frac{1}{q}}}
        \\[1ex]
        & \le M (t - s)^{-\tfrac{N}{2} \bigl(\tfrac{1}{p} - \tfrac{1}{q}\bigr)} e^{\gamma (t - s)} \,
        \biggl( \sum\limits_{k=1}^{n} \lVert u^k_0 \rVert^p_{L_p(D)} \biggr)^{\!\frac{1}{p}}
        \\
        & = M (t - s)^{-\tfrac{N}{2} \bigl(\tfrac{1}{p} - \tfrac{1}{q}\bigr)} e^{\gamma (t - s)} \lVert u_0 \rVert_{L_p(D)^n}.
    \end{align*}
    The remaining cases are considered in a similar way.
\end{proof}

Below is a corollary of Propositions~\ref{prop:continuity-in-norm-top} and~\ref{prop:norm_continuity-all}.
\begin{proposition}
\label{prop:continuiuty-to-C}
    Assume moreover~\textup{\ref{AS-ae-converge}}.  Let $1 \le p \le q < \infty$ and $0 \le s < T$.  Then
    \begin{enumerate}[label=\textup{(\roman*)},ref=\ref{prop:continuiuty-to-C}\textup{(\roman*)}]
        \item\label{prop:continuiuty-to-C(i)}
        for any $t_0 \in (s, T)$ the mapping
        \begin{equation*}
            \bigl[\, Y_0 \times L_p(D)^n \ni (a_0, u_0) \mapsto U_{a_0}(\cdot, s) u_0 \in C([t_0, T], L_q(D)^n) \, \bigr]
        \end{equation*}
        is continuous;
        \item\label{prop:continuiuty-to-C(ii)}
        for $1 < p = q$ the mapping
        \begin{equation*}
            \bigl[\, Y_0 \times L_p(D)^n \ni (a_0, u_0) \mapsto U_{a_0}(\cdot, s) u_0 \in C([s, T], L_p(D)^n) \, \bigr]
        \end{equation*}
        is continuous.
    \end{enumerate}
\end{proposition}

Let $a \in Y$ be fixed.  For Lebesgue\nobreakdash-\hspace{0pt}a.e.\ $t \in [0, T]$ and $1 \le p \le \infty$ we define multiplication operators from $L_p(D)^n$ into $L_p(D)^n$ by the formulas
\begin{align*}
    (\mathcal{C}^0_{a}(t) u)[x] & = \Bigl(\sum\limits_{l = 1}^{n} c_0^{1l}(t, x) u^l{[x]}, \dots, \sum\limits_{l = 1}^{n} c_0^{nl}(t, x) u^l{[x]}\Bigr),
    \\
    (\mathcal{C}^1_{a}(t) u)[x] & = \Bigl(\sum\limits_{l = 1}^{n} c_1^{1l}(t, x) u^l{[x]}, \dots, \sum\limits_{l = 1}^{n} c_1^{nl}(t, x) u^l{[x]}\Bigr)
\end{align*}
for Lebesgue\nobreakdash-\hspace{0pt}a.e.\ $x \in D$.

In other words, if we define, for Lebesgue\nobreakdash-\hspace{0pt}a.e.\ $(t, x) \in (0, T) \times D$, matrices $c_i(t, x) := [c_i^{kl}(t, x)]_{k,l=1}^n$, $i = 0, 1$, we can write
\begin{equation*}
    \mathcal{C}_a^{i}(t) u = c_i(t, \cdot) u(\cdot) \quad \text{for Lebesgue-a.e. } t \in (0, T).
\end{equation*}

\begin{definition}
\label{def:mild-solution}
    By a \emph{mild $L_p$\nobreakdash-\hspace{0pt}solution}, $1 < p < \infty$, of \eqref{main-eq}$_{a}$\textup{+}\eqref{main-bc}$_{a}$ satisfying the initial conditions \eqref{main-ic}, where $0 \le s < T$, $u_0 = (u_0^{(1)}, u_0^{(2)}) \in L_p(D)^n \times L_r((- 1, 0), L_p(D)^n)$, $1 \le r \le \infty$, we understand a continuous  function \newline $u_p = u_p(\cdot; s, a, u_0) \colon [s, T] \to L_p(D)^n$ satisfying the integral equation
    \begin{equation}
    \label{eq:integral}
        \begin{aligned}
        u_p(t) = {} & U_{a_0}^0(t, s) u_0^{(1)} + \int\limits_{s}^{t} U_{a_0}^0(t,\zeta) \bigl( \mathcal{C}^0_a(\zeta) u_p(\zeta) \bigr) \, \mathrm{d}\zeta
        \\
        & + \int\limits_{s}^{t} U_{a_0}^0(t, \zeta) \bigl( \mathcal{C}^1_a(\zeta) u_p(\zeta - 1) \bigr) \, \mathrm{d}\zeta, \quad s \le t \le T,
        \end{aligned}
    \end{equation}
where in the third term on the right-hand side $u_p(\zeta - 1)$ is, for Lebesgue\nobreakdash-\hspace{0pt}a.e.\ $\zeta \in (s, s + 1)$, understood as $u_0^{(2)}(\zeta -s - 1)$.
\end{definition}
\begin{definition}
\label{def:mild-solution-L1}
    By a \emph{mild $L_1$\nobreakdash-\hspace{0pt}solution} of \eqref{main-eq}$_{a}$\textup{+}\eqref{main-bc}$_{a}$ satisfying the initial conditions \eqref{main-ic}, where $0 \le s < T$, $u_0 = (u_0^{(1)}, u_0^{(2)}) \in L_1(D)^n \times L_r((- 1, 0), L_1(D)^n)$, $1 \le r \le \infty$, we understand a bounded continuous function $u_1 = u_1(\cdot; s, a, u_0) \colon (s, T] \to L_1(D)^n$ satisfying the integral equation
    \begin{equation}
    \label{eq:integral-L1}
        \begin{aligned}
        u_{1}(t) = {} & U_{a_0}^0(t, s) u_0^{(1)} + \int\limits_{s}^{t} U_{a_0}^0(t,\zeta) \bigl( \mathcal{C}^0_a(\zeta) u_{1}(\zeta) \bigr) \, \mathrm{d}\zeta
        \\
        & + \int\limits_{s}^{t} U_{a_0}^0(t, \zeta) \bigl( \mathcal{C}^1_a(\zeta) u_{1}(\zeta - 1) \bigr) \, \mathrm{d}\zeta, \quad s < t \le T,
        \end{aligned}
    \end{equation}
where in the third term on the right-hand side $u_{1}(\zeta - 1)$ is, for Lebesgue\nobreakdash-\hspace{0pt}a.e.\ $\zeta \in (s, s + 1)$, understood as $u_0^{(2)}(\zeta -s - 1)$.
\end{definition}

From now on, generic theorems will be stated under the assumption that the initial time $s$ is taken to be $0$.  Passage to a case of a general $s \in [0, T)$ is straightforward.

\medskip
We introduce some notation.  Let $1 \le \xi \le \infty$, $1 \le \eta \le \infty$.  For an $n \times n$ matrix $g = [g^{kl}]_{k,l=1}^n$ with entries in $L_{\infty}(D)$ we denote
\begin{multline}
\label{eq:definition-matrix-norms}
   \norm{g}_{\xi,\eta}
   \\
   := \begin{cases}
           \displaystyle \Biggl( \sum_{k = 1}^n  \Bigl( \sum_{l = 1}^n \norm{g^{kl}}_{L_{\infty}(D)}^{\xi} \Bigr)^{\frac{\eta}{\xi}} \Biggr)^{\frac{1}{\eta}}, & \text{ for } 1 \le \xi < \infty, \ 1 \le \eta < \infty,
            \\[3ex]
            \displaystyle \Biggl( \sum_{k = 1}^n  \Bigl( \max_{1\le l\le n} \norm{g^{kl}}_{L_{\infty}(D)} \Bigr)^{\eta} \Biggr)^{\frac{1}{\eta}}, & \text{ for } \xi = \infty, \ 1 \le \eta < \infty,
            \\[3ex]
            \displaystyle \max_{1 \le k\le n} \Bigl( \sum_{l = 1}^n \norm{g^{kl}}_{L_{\infty}(D)}^{\xi} \Bigr)^{\frac{1}{\xi}},  & \text{ for } 1 \le \xi < \infty, \ \eta = \infty,
            \\[4ex]
            \displaystyle  \max_{1\le k\le n}  \Bigl( \max_{1\le l\le n} \norm{g^{kl}}_{L_{\infty}(D)} \Bigr), & \text{ for } \xi = \eta = \infty.
            \\[3ex]
          \end{cases}
\end{multline}
From the inequalities on the norms in finite dimensional real vector spaces it follows that
\begin{equation}
\label{eq:matrix-norms-inequalities}
    \norm{g}_{\xi,\eta} \le \norm{g}_{1,1}, \quad 1 \le \xi \le \infty, \ 1 \le \eta \le \infty.
\end{equation}

The next result is, at~least for $n = 1$, a standard exercise (see, e.g., \cite[Example~III.2.2]{Con}).  We give it for the sake of completeness.
\begin{lemma}
\label{lm:global-mult-estim}
    Let $a \in Y$.  Then for each $i = 0, 1$ and $1 \le p \le \infty$ there holds
    \begin{equation*}
        \lVert \mathcal{C}^i_{a}(t) \rVert_{\mathcal{L}(L_p(D)^n)} \le  \lVert c_i(t,\cdot)\rVert_{p',p}
    \end{equation*}
    for Lebesgue\nobreakdash-\hspace{0pt}a.e.\ $t \in (0, T)$.
\end{lemma}
\begin{proof}
    Assume first $1 < p < \infty$.  For $u \in L_p(D)^n$  and $t \in (0, T)$ we have
    \begin{equation*}
      \lVert \mathcal{C}^i_{a}(t) u \rVert_{L_p(D)^n} = \biggl( \int_{D} \sum_{k = 1}^n \bigl\lvert \sum_{l = 1}^n c^{kl}_i(t, x) u^l{[x]} \bigr\rvert^{p} \, \mathrm{d}x \biggr)^{\frac{1}{p}}.
  \end{equation*}
  By the H\"older inequality,
\begin{align*}
          \bigl\lvert \sum_{l = 1}^n c^{kl}_i(t, x) u^l{[x]} \bigr\rvert &\le \sum_{l = 1}^n \lvert c^{kl}_i(t, x) u^l{[x]} \rvert \\
          & \le \Bigl( \sum_{l = 1}^n \abs{c^{kl}_i(t, x)}^{p'} \Bigr)^{\frac{1}{p'}} \Bigl( \sum_{l = 1}^n \lvert u^{l}{[x]} \rvert^{p} \Bigr)^{\frac{1}{p}}.
\end{align*}
  Further,
  \begin{equation*}
      \Bigl( \sum_{l = 1}^n \lvert c^{kl}_i(t, x) \rvert^{p'} \Bigr)^{\frac{p}{p'}} \le \Bigl( \sum_{l = 1}^n \norm{ c^{kl}_i(t, \cdot ) }_{L_{\infty}{(D)}}^{p'} \Bigr)^{\frac{p}{p'}}
  \end{equation*}
  so
\begin{align*}
    \sum_{k = 1}^n \Bigl( \sum_{l = 1}^n \lvert c^{kl}_i(t, x) \rvert^{p'} \Bigr)^{\frac{p}{p'}} & \le  \sum_{k = 1}^n \Bigl( \sum_{l = 1}^n \norm{ c^{kl}_i(t, \cdot ) }_{L_{\infty}{(D)}}^{p'} \Bigr)^{\frac{p}{p'}}\\
    & = \Biggl(\sum_{k = 1}^n \Bigl( \sum_{l = 1}^n \norm{ c^{kl}_i(t, \cdot ) }_{L_{\infty}{(D)}}^{p'} \Bigr)^{\frac{p}{p'}} \Biggr)^{\frac{1}{p}\cdot p}
    \\
    & = \lVert [c^{kl}_i(t,\cdot)]_{k,l=1}^{n} \rVert_{p',p}^p.
\end{align*}
 We have thus
\begin{align*}
          \lVert \mathcal{C}^i_{a}(t) u \rVert_{L_p(D)^n} &\le \biggl(  \lVert [c^{kl}_i(t,\cdot)]_{k,l=1}^{n} \rVert^p_{p',p}   \int_{D} \Bigl( \sum_{l = 1}^n \lvert u^{l}[x] \rvert^{p} \Bigr)  \, \mathrm{d}x \biggr)^{\frac{1}{p}} \\
          & =  \lVert [c^{kl}_i(t,\cdot)]_{k,l=1}^{n} \rVert_{p',p}  \biggl( \int_{D} \sum_{l = 1}^n \lvert u^{l}[x] \rvert^p \, \mathrm{d}x \biggr)^{\frac{1}{p}}  \\
          & =  \lVert [c^{kl}_i(t,\cdot)]_{k,l=1}^{n} \rVert_{p',p} \, \norm{u}_{L_p(D)^n}.
\end{align*}
For $p = 1$ we estimate
\begin{align*}
    \lVert \mathcal{C}^i_{a}(t) u \rVert_{L_1(D)^n} & =  \int_{D} \sum_{k = 1}^n \bigl\lvert \sum_{l = 1}^n c^{kl}_i(t, x) u^l{[x]} \bigr\rvert \, \mathrm{d}x
    \\
    & \le \sum_{k = 1}^n  \biggl( \int_{D} \sum_{l = 1}^n \lvert c^{kl}_i(t, x) u^l{[x]} \rvert \, \mathrm{d}x \biggr)
    \\
    & \le \sum_{k = 1}^n \biggl( \int_{D}  \sum_{l = 1}^n \lVert c^{kl}_i(t, \cdot) \lVert_{L_{\infty}(D)} \lvert u^l{[x]} \rvert \, \mathrm{d}x \biggr)
    \\
    & \le \sum_{k = 1}^n \biggl( \max\limits_{1 \le l \le n} \lVert c^{kl}_i(t, \cdot) \lVert_{L_{\infty}(D)} \sum_{l = 1}^n \int_{D} \lvert u^l{[x]} \rvert \, \mathrm{d}x \biggr)
    \\
    & = \biggl( \sum_{k = 1}^n \max\limits_{1 \le l \le n} \lVert c^{kl}_i(t, \cdot) \lVert_{L_{\infty}(D)} \biggr) \biggl( \sum_{l = 1}^n \int_{D} \lvert u^l{[x]} \rvert \, \mathrm{d}x \biggr)
    \\
    & = \norm{[c^{kl}_i(t,\cdot)]_{k,l=1}^n}_{\infty,1} \, \lVert u \rVert_{L_1(D)^n}.
\end{align*}
For $p = \infty$ we estimate
\begin{align*}
    \lVert \mathcal{C}^i_{a}(t) u \rVert_{L_{\infty}(D)^n} & = \max\limits_{1 \le k \le n} \biggl( \esssup\limits_{x \in D} \Bigl\lvert \sum_{l = 1}^n c^{kl}_i(t, x) u^l{[x]} \Bigr\rvert \biggr)
    \\
    & \le \max\limits_{1 \le k \le n} \biggl( \sum_{l = 1}^n \esssup\limits_{x \in D}{ \lvert c^{kl}_i(t, x) u^l{[x]} \rvert} \biggr)
    \\
    & \le \max\limits_{1 \le k \le n}  \biggl( \sum_{l = 1}^n \bigl(\esssup\limits_{x \in D} \abs{c^{kl}_i(t, x)} \bigr) \lVert u^l \rVert_{L_{\infty}(D)} \biggr)
    \\
    & \le \max\limits_{1 \le k \le n} \biggl( \Bigl( \sum\limits_{l = 1}^n \lVert c^{kl}_i(t, \cdot) \lVert_{L_{\infty}(D)} \Bigr) \max_{1 \le l \le n} \lVert u^l \rVert_{L_{\infty}(D)} \biggr)
    \\
    & = \max\limits_{1 \le k \le n} \Bigl( \sum_{l = 1}^n \norm{c^{kl}_i(t, \cdot)}_{L_{\infty}(D)} \Bigr) \max_{1 \le l \le n} \lVert u^l \rVert_{L_{\infty}(D)}
    \\
    & = \norm{[c^{kl}_i(t,\cdot)]_{k,l=1}^n}_{1,\infty} \, \lVert u \rVert_{L_{\infty}(D)^n}.
\end{align*}
\end{proof}

Let $K$ stand for the supremum of the $L_{\infty}((0, T) \times D)$\nobreakdash-\hspace{0pt}norms of $c_i^{kl}$ over all $a \in Y$.  In~view of~\ref{AS-compact}, $K < \infty$.

As a consequence of Lemma~\ref{lm:global-mult-estim} and~\eqref{eq:matrix-norms-inequalities} the following holds.
\begin{lemma}
\label{lm:mult-estim}
    For each $i = 0, 1$, each $1 \le p \le \infty$ and each $a \in Y$ we have
    \begin{equation*}
        \lVert \mathcal{C}^i_{a}(t) \rVert_{\mathcal{L}(L_p(D)^n)} \le n^2 K
    \end{equation*}
    for Lebesgue\nobreakdash-\hspace{0pt}a.e.\ $t \in (0, T)$.
\end{lemma}

\begin{lemma}
\label{lm:mult-estim-1}
    Assume $1 \le p < \infty$ and $a \in Y$.  Then for each $u \in L_1((0, T), L_p(D)^n)$ the function     assigning to Lebesgue\nobreakdash-\hspace{0pt}a.e.\ $t \in (0, T)$ the value $\mathcal{C}^{i}_{a}(t) u(t) \in L_p(D)^n$ is $(\mathfrak{L}((0, T)), \mathfrak{B}(L_p(D)^n))$\nobreakdash-\hspace{0pt}measurable.  Moreover, the assignment
    \begin{equation*}
        \bigl[\, u(\cdot) \mapsto \mathcal{C}^{i}_{a}(\cdot) u(\cdot) \,\bigr]
    \end{equation*}
    has $\mathcal{L}(L_r((0, T), L_p(D)^n))$\nobreakdash-\hspace{0pt}norm not larger than $n^2 K$, for any $1 \le r \le \infty$.
\end{lemma}
\begin{proof}
    Since $L_1((0, T), L_p(D)^n)$ embeds continuously in $L_1((0, T), L_1(D)^n)$, for any $u \in L_1((0, T), L_p(D)^n)$ the function
    \begin{equation*}
        \bigl[\, (0, T) \times D \ni (t, x) \mapsto u(t)[x] \in \RR^n \,\bigr]
    \end{equation*}
    belongs, by Lemma~\ref{lm:Dunford-Schwartz}(a), to $L_1((0, T) \times D, \RR^n)$, so is $(\mathfrak{L}((0, T) \times D), \mathfrak{B}(\RR^n))$\nobreakdash-\hspace{0pt}measurable.  Therefore the function
    \begin{equation*}
        \bigl[\, (0, T) \times D \ni (t, x) \mapsto c_i(t, x) \cdot u(t)[x] \in \RR^n \,\bigr],
    \end{equation*}
    being the product of measurable functions, is $(\mathfrak{L}((0, T) \times D), \mathfrak{B}(\RR^n))$\nobreakdash-\hspace{0pt}measurable, too.  By Lemma~\ref{lm:mult-estim}, for Lebesgue\nobreakdash-\hspace{0pt}a.e.\ $t \in (0, T)$, the section $c_i(t, \cdot) u(t) = \mathcal{C}^i_{a}(t) u(t)$ belongs to $L_p(D)^n$, so, applying Lemma~\ref{lm:Dunford-Schwartz}(b) we obtain the first part of the lemma.

    Since $L_r((0, T), L_p(D)^n)$, $1 \le r < \infty$, equals the set of $(\mathfrak{L}((0, T), \mathfrak{B}(L_p(D)^n))$\nobreakdash-\hspace{0pt}measurable functions $u$ for which $\int\limits_{0}^{T} \norm{u(t)}_{L_p(D)^n}^r \, \mathrm{d}t < \infty$, for such $u$ there holds, by Lemma~\ref{lm:mult-estim},
    \begin{align*}
        & \biggl( \int_{0}^{T} \norm{\mathcal{C}^i_{a}(t) u(t)}_{L_p(D)^n}^r \, \mathrm{d}t  \biggr)^{\!\! 1/r}
        \\
        \le {} & \esssup\{\, \norm{{\mathcal{C}}^i_{a}(t)}_{\mathcal{L}(L_p(D)^n)} : t \in (0, T) \,\} \cdot \biggl( \int_{0}^{T} \norm{u(t)}_{L_p(D)^n}^r \, \mathrm{d}t  \biggr)^{\!\! 1/r}
        \\
        \le {} & n^2 K \lVert u \rVert_{L_r((0, T), L_p(D)^n)}.
    \end{align*}
    Similarly, for $u \in L_{\infty}((0, T), L_p(D)^n)$,
    \begin{align*}
        & \esssup\{\, \norm{{\mathcal{C}}^i_{a}(t) u(t)}_{L_p(D)^n} : t \in (0, T) \,\}
        \\
        \le {} & \esssup\{\, \norm{{\mathcal{C}}^i_{a}(t)}_{\mathcal{L}(L_p(D)^n)} : t \in (0, T) \,\} \cdot \esssup\{\, \norm{u(t)}_{L_p(D)^n} : t \in (0, T) \,\}
        \\
        \le {} & n^2 K \norm{u}_{L_{\infty}((0, T), L_p(D)^n)}.
    \end{align*}
\end{proof}

\begin{lemma}
\label{lm:I_i}
    Assume $1 \le p \le q < \infty$, $1 \le r \le \infty$, $a_0 \in Y_0$ and $u \in L_r((0, T), L_p(D)^n)$, with
    \begin{equation}
    \label{eq:ass-p-q-0}
        \frac{N}{2}\biggl(\frac{1}{p}- \frac{1}{q}\biggr) < \frac{1}{r'}
    \end{equation}
    or
    \begin{equation*}
        p = q.
    \end{equation*}
    Then
    \begin{enumerate}
    [label=\textup{(\roman*)},ref=\ref{lm:I_i}\textup{(\roman*)}]
        \item
        \label{lm:I_i1}
        for any $0 < t \le T$ the function
        \begin{equation}
        \label{eq:I_i}
            [\, \zeta \mapsto U^0_{a_0}(t, \zeta) u(\zeta)\,]
        \end{equation}
        belongs to $L_{w}((0, t), L_{q}(D)^n)$,  where
        \begin{equation*}
            w \ge 1, \ \frac{1}{w} > \frac{1}{r} + \frac{N}{2}\biggl(\frac{1}{p}- \frac{1}{q}\biggr),
        \end{equation*}
        if~\eqref{eq:ass-p-q-0} is satisfied, and $w = r$ if $p = q$; moreover, the linear operator assigning~\eqref{eq:I_i} to $u$ belongs to $\mathcal{L}(L_r((0,t), L_p(D)^n), L_w((0,t), L_p(D)^n))$, with the norm bounded uniformly in $a_0 \in Y_0$;
        \item
        \label{lm:I_i2}
        the mapping
        \begin{equation}
        \label{eq:I_i-1}
            \biggl[\, [0, T] \ni t \mapsto \int_{0}^{t} U^0_{a_0}(t, \zeta) u(\zeta) \, \mathrm{d}\zeta \,\biggr]
        \end{equation}
        belongs to $C([0, T], L_q(D)^n)$.
    \end{enumerate}
\end{lemma}
\begin{proof}
    Fix $0 < t \le T$.  We show first that \eqref{eq:I_i} defines a $(\mathfrak{L}((0, t)), \mathfrak{B}(L_{q}(D)^n))$\nobreakdash-\hspace{0pt}measurable function.

    Consider first the case $p > 1$.  The property is equivalent, by Theorem~\ref{thm:equiv-measurable}, to showing that for each $v \in L_{q'}(D)^n$ the function
        \begin{equation*}
            \bigl[\, \zeta \mapsto \langle U^0_{a_0}(t, \zeta) u(\zeta), v \rangle_{L_{q}(D)^n, L_{q'}(D)^n} \,\bigr]
        \end{equation*}
    is $(\mathfrak{L}((0, t)), \mathfrak{B}(\RR))$\nobreakdash-\hspace{0pt}measurable.  By Proposition~\ref{prop:dual-all} applied twice, we need to show that the function
        \begin{equation*}
            \bigl[\, \zeta \mapsto \langle u(\zeta), U^{0,*}_{a_0}(\zeta, t) v \rangle_{L_{p}(D)^n, L_{p'}(D)^n} \,\bigr]
        \end{equation*}
    is $(\mathfrak{L}((0, t)), \mathfrak{B}(\RR))$\nobreakdash-\hspace{0pt}measurable, which follows from the fact that $u$ is \newline $(\mathfrak{L}((0, t)), \mathfrak{B}(L_p(D)^n))$\nobreakdash-\hspace{0pt}measurable (by assumption) and the function
        \begin{equation*}
            \bigl[\, [0, t) \ni \zeta \mapsto U^{0,*}_{a_0}(\zeta, t) v \in L_{p'}(D)^n \,\bigr]
        \end{equation*}
    is continuous (by the adjoint equation analog of Lemma~\ref{lm:continuity-in-strong-top}).

    We proceed to the case $p = 1$.  Let, first, $q = 1$.  The mapping~\eqref{eq:I_i} can be represented as $Q \circ P$, where $P \colon (0, t) \to \mathcal{L}_{\mathrm{s}}(L_1(D)^n) \times L_1(D)^n$ is given as $P(\zeta)  = (U^0_{a_0}(t, \zeta), u(\zeta))$ and $Q$ is the evaluation map of an operator in $\mathcal{L}_{\mathrm{s}}(L_1(D)^n)$ on an element of $L_1(D)^n$, restricted to the set
    \begin{equation*}
        \{\, A \in \mathcal{L}_{\mathrm{s}}(L_1(D)^n) : \norm{A}_{\mathcal{L}(L_1(D)^n} \le M e^{\gamma T} \,\} \times L_1(D)^n,
    \end{equation*}
    where $M$ and $\gamma$ are as in Proposition~\ref{prop:system-estimates}(ii).  The mapping $Q$ is continuous, cf. \cite[Lemma~A.6(1)]{GT-Q}.  The first coordinate of $P$ is continuous, by Proposition~\ref{prop:continuity-in-norm-top}, hence $\mathfrak{L}((0, t)), \mathfrak{B}(\mathcal{L}_{\mathrm{s}}(L_1(D)^n))$\nobreakdash-\hspace{0pt}measurable (for more on Borel sets in $\mathcal{L}_{\mathrm{s}}(L_1(D)^n)$ see, e.g., \cite[Appendix~A]{GT-Q}), and the second coordinate is $(\mathfrak{L}((0, t)), \mathfrak{B}(L_1(D)^n))$\nobreakdash-\hspace{0pt}measurable by~assumption.  Consequently, $P$ is $(\mathfrak{L}((0, t)), \mathfrak{B}(\mathcal{L}_{\mathrm{s}}(L_1(D)^n)) \otimes \mathfrak{B}(L_1(D)^n))$-measurable, which is equivalent to the $(\mathfrak{L}((0, t)), \mathfrak{B}(\mathcal{L}_{\mathrm{s}}(L_1(D)^n) \times L_1(D)^n))$-measurability (see \cite[Lemma~6.4.2(i)]{Bogachev}).  If $q > 1$ we observe that, for Lebesgue\nobreakdash-\hspace{0pt}a.e. \ $\zeta \in (0, t)$, \newline $\langle U^0_{a_0}(t, \zeta) u(\zeta), v \rangle_{L_{q}(D)^n, L_{q'}(D)^n} = \langle U^0_{a_0}(t, \zeta) u(\zeta), v \rangle_{L_{1}(D)^n, L_{\infty}(D)^n}$ and proceed as in the case $p > 1$.

    Assume~\eqref{eq:ass-p-q-0}.  It follows from Proposition~\ref{prop:higher_orders-skew_product-p-q-estimates}(ii) that the function $[\, (0, t) \ni \zeta \mapsto \norm{U^0_{a_0}(t, \zeta)}_{\mathcal{L}(L_p(D)^n, L_q(D)^n)} \,]$ belongs to $L_z((0, t))$ with $z = 1/(\frac{N}{2}(\frac{1}{p}- \frac{1}{q}))$.  Then the membership of~\eqref{eq:I_i} in $L_{w}((0, t), L_{q}(D)^n)$ as well as the bound on its norm follow from the generalized H\"older inequality.  The case $p = q$ is the consequence of the fact that, by Proposition~\ref{prop:higher_orders-skew_product-p-q-estimates}(ii), the function $[\, (0, t) \ni \zeta \mapsto \norm{U^0_{a_0}(t, \zeta)}_{\mathcal{L}(L_p(D)^n)} \,]$ belongs to $L_{\infty}((0, t))$. The proof of part (i) is thus completed.

    \smallskip
    We proceed to the proof of part (ii).  By part (i), the function~\eqref{eq:I_i-1} is well defined.  Let $0 \le t_1 \le t_2 \le T$.  We write
    \begin{multline*}
        \int_{0}^{t_2} U^{0}_{a_0}(t_2, {\zeta}) u({\zeta}) \, \mathrm{d}{\zeta} - \int_{0}^{t_1} U^{0}_{a_0}(t_1, {\zeta}) u({\zeta}) \, \mathrm{d}{\zeta}
        \\
        = \int_{0}^{t_1} \bigl( U^{0}_{a_0}(t_2, {\zeta}) - U^{0}_{a_0}(t_1, {\zeta})\bigr) u({\zeta}) \, \mathrm{d}{\zeta} + \int_{t_1}^{t_2} U^{0}_{a_0}(t_2, {\zeta}) u({\zeta}) \, \mathrm{d}{\zeta}.
    \end{multline*}
    Let $\epsilon > 0$.  As \eqref{eq:I_i} belongs to $L_{w}((0, t), L_{q}(D)^n)$, it is a consequence of~\cite[Thm.~II.2.4(i)]{DiUhl} that the $L_q(D)^n$\nobreakdash-\hspace{0pt}norm of the second term on the right-hand side can be made $< \epsilon/3$ by taking $t_1, t_2$ sufficiently close to each other.  Regarding the first term, we write
    \begin{multline*}
        \int_{0}^{t_1} \bigl( U^{0}_{a_0}(t_2, {\zeta}) - U^{0}_{a_0}(t_1, {\zeta})\bigr) u({\zeta}) \, \mathrm{d}{\zeta}
        \\
        = \int_{0}^{t_1 - \eta} \bigl( U^{0}_{a_0}(t_2, {\zeta}) - U^{0}_{a_0}(t_1, {\zeta})\bigr) u({\zeta}) \, \mathrm{d}{\zeta} + \int_{t_1 - \eta}^{t_1} \bigl( U^{0}_{a_0}(t_2, {\zeta}) - U^{0}_{a_0}(t_1, {\zeta})\bigr) u({\zeta}) \, \mathrm{d}{\zeta}.
    \end{multline*}
    Again by \cite[Thm.~II.2.4(i)]{DiUhl}, for $\eta > 0$ sufficiently small there holds
    \begin{equation*}
            \biggl\lVert \int_{t_1 - \eta}^{t_1} \bigl( U^{0}_{a_0}(t_2, {\zeta}) - U^{0}_{a_0}(t_1, {\zeta})\bigr) u({\zeta}) \, \mathrm{d}{\zeta} \biggr\rVert_{L_q(D)^n} < \frac{\epsilon}{3}.
    \end{equation*}
    It follows from Proposition~\ref{prop:continuity-in-norm-top} that the assignment
    \begin{equation*}
        \bigl[\, \{\, \dot{\Delta} : \eta \le \zeta + \eta \le t \le T \,\} \ni (\zeta, t) \mapsto U^0_{a_0}(t, \zeta) \in \mathcal{L}(L_p(D)^n, L_q(D)^n) \,\bigr]
   \end{equation*}
    is uniformly continuous, consequently there exists $\delta > 0$ such that if $\eta \le \zeta + \eta \le t_1 \le t_2$, $t_2 - t_1 < \delta$, then
    \begin{equation*}
        \lVert U^0_{a_0}(t_2, \zeta) - U^0_{a_0}(t_1, \zeta) \rVert_{\mathcal{L}(L_p(D)^n, L_q(D)^n)}  < \frac{\epsilon}{3 \norm{u}_{L_1((0, T), L_p(D)^n)}}.
    \end{equation*}
   Therefore
    \begin{equation*}
        \biggl\lVert \int_{0}^{t_1 - \eta} \bigl( U^0_{a_0}(t_2, {\zeta}) - U^0_{a_0}(t_1, {\zeta})\bigr) u({\zeta}) \, \mathrm{d}{\zeta} \biggr\rVert_{L_q(D)^n}
         < \frac{\epsilon}{3}.
    \end{equation*}
    This concludes the proof of part (ii).
\end{proof}

We will now present a counterpart of the above lemma for $q = \infty$. Since we lack separability of $L_{\infty}(D)^n$, in the text below the integrals of $L_{\infty}(D)^n$-valued functions will be understood in the Gel$\textprime$fand sense. See~\cite{Hash} for more details.
\begin{lemma}
\label{lm:L-_infty}
     Assume $1 \le p < \infty$, $1 < r \le \infty$, $a_0 \in Y_0$ and $u \in L_r((0, T), L_p(D)^n)$, with
    \begin{equation}
    \label{eq:ass-p-q-infty}
       \frac{N}{2p} < \frac{1}{r'}
    \end{equation}
    Then for any $t\in[0,T]$ the mapping
    \begin{equation*}
        [\, (0,t) \ni \zeta \mapsto U^0_{a_0}(t, \zeta) u(\zeta) \in L_{\infty}(D)^n \,]
    \end{equation*}
    is weak-* integrable. In other words
        \begin{equation}
        \label{eq:integral-q-infty}
           \textup{G -}\hspace{-1pt} \int_{0}^{t} U^0_{a_0}(t, \zeta) u(\zeta) \, \mathrm{d}\zeta \in L_{\infty}(D)^n.
        \end{equation}
\end{lemma}
\begin{proof}
    We start by demonstrating the weak-* measurability of the mapping
    \begin{equation*}
         \bigl[\, [0, t] \ni \zeta \mapsto U^0_{a_0}(t, \zeta) u(\zeta)  \in L_{\infty}(D)^n \,\bigr]
    \end{equation*}
or equivalently $L_1(D)^n$-measurability in the sense of~\cite[see p. 304]{Hash}.  Fix  $v\in L_1(D)^n$.  It suffices to show that the mapping
        \begin{equation}
        \label{eq:L_infty-L_1-par}
            \bigl[\, \zeta \mapsto \langle v, U^0_{a_0}(t, \zeta) u(\zeta) \rangle_{L_{1}(D)^n, L_{\infty}(D)^n} \,\bigr]
        \end{equation}
    is $(\mathfrak{L}((0, t)), \mathfrak{B}(\RR))$\nobreakdash-\hspace{0pt}measurable. It can be done via approximation. Let $(v_m)_{m = 1}^{\infty} \subset L_{2}(D)^n$ be a sequence  approximating $v$ in the $L_1(D)^n$-norm. Since for any $m\in\NN$ we have
    \begin{equation*}
        ( v_m, U^0_{a_0}(t, \zeta) u(\zeta) )_{L_{2}(D)^n} = \langle v_m, U^0_{a_0}(t, \zeta) u(\zeta) \rangle_{L_{1}(D)^n, L_{\infty}(D)^n}
    \end{equation*}
    and the mapping
            \begin{equation}
            \label{eq:L_p-L_p'-aprox}
            \bigl[\, \zeta \mapsto ( v_m, U^0_{a_0}(t, \zeta) u(\zeta) )_{L_{2}(D)^n} \,\bigr]
        \end{equation}
        is $(\mathfrak{L}((0, t)), \mathfrak{B}(\RR))$\nobreakdash-\hspace{0pt}measurable by Theorem~\ref{thm:equiv-measurable}, it remains to observe that the mapping~\eqref{eq:L_infty-L_1-par} is the point-wise limit of mappings~\eqref{eq:L_p-L_p'-aprox}. Moreover, the weak-* integrability follows from the following estimate, using Proposition~\ref{prop:system-estimates}(ii) and the Hölder inequality
        \begin{multline*}
                      \int_{0}^{t} \abs{\langle v, U^0_{a_0}(t, \zeta) u(\zeta) \rangle_{L_{1}(D)^n, L_{\infty}(D)^n}} \, \mathrm{d}\zeta   \\
             \le  \int_{0}^{t} \norm{ U^0_{a_0}(t, \zeta) }_{\mathcal{L}(L_p(D)^n, L_{\infty}(D)^n)} \norm{u(\zeta)}_{L_{p}(D)^n} \norm{v}_{L_{1}(D)^n} \, \mathrm{d}\zeta \\
            \le M e^{\gamma T} \biggl(  \int_{0}^t (t - \zeta)^{- \frac{N r'}{2p}} \, \dd\zeta \biggr)^{\!\! 1/r'} \biggl(  \int_{0}^t \norm{u(\zeta)}_{L_p(D)^n}^{r} \, \dd\zeta \biggr)^{\!\! 1/r}\norm{v}_{L_{1}(D)^n}
            \\
             =\frac{M e^{\gamma T }}{\bigl(1 -  \frac{N r'}{2p} \bigr)^{1/r'}} t^{\frac{1}{r'} -  \frac{N }{2p}} \biggl(  \int_{0}^t \norm{u(\zeta)}_{L_p(D)^n}^{r} \, \dd\zeta \biggr)^{\!\! 1/r} \norm{v}_{L_{1}(D)^n} .
        \end{multline*}
       The proof is concluded by application of~\cite[Lemma 1.1]{Hash}.
\end{proof}

\begin{proposition}
  \label{prop:delay-existence-hull}
  For any $1 \le p < \infty$, $1 \le r \le \infty$, $a \in Y$ and $u_{0} = (u_0^{(1)}, u_0^{(2)}) \in L_p(D)^n \times L_r((- 1, 0), L_p(D)^n)$ there exists a unique mild $L_p$-solution $u_p(\cdot; 0, a, u_0) \colon [0, T] \to L_p(D)^n$ of \eqref{main-eq}$_{a}$\textup{+}%
  \eqref{main-bc}$_{a}$ satisfying the initial condition
  \begin{equation}
  \label{eq-hull-ic}
  \begin{cases}
    u(0) = u_0^{(1)},
    \\
    u(\tau) = u_0^{(2)}(\tau) \text{ for a.e. } \tau \in (- 1, 0).
  \end{cases}
  \end{equation}
\end{proposition}

The idea of the proof runs as~follows.
The solution is obtained as the sum of the function $[\, (0, T] \ni t \mapsto U_{a_0}^0(t, 0) u_0^{(1)} \,]$ and a fixed point of the contraction mapping $\mathfrak{G}_{a, u_0}$ of a closed subset of $C([0, T], L_p(D)^n)$ into itself defined as
\begin{equation}
\label{eq:Gothic-G}
  \begin{aligned}
  \mathfrak{G}_{a, u_0}(u)[t] := {} & \int\limits_{0}^{t} U_{a_0}^0(t, \zeta) \bigl( \mathcal{C}^0_{a}(\zeta) (u(\zeta) + U_{a_0}^0(\zeta, 0) u_0^{(1)}) \bigr) \, \mathrm{d}\zeta
  \\
      & + \int\limits_{0}^{t} U_{a_0}^0(t, \zeta) \bigl( \mathcal{C}^1_{a}(\zeta) (u_0^{(2)}(\zeta - 1) + \mathbbm{1}_{\! (1, T)}(\zeta)  U_{a_0}^0(\zeta, 0) u_0^{(1)})\bigr) \, \mathrm{d}\zeta, \quad 0 \le t \le T.
  \end{aligned}
\end{equation}

Until revoking, $u$ stands for a generic function in $C([0, T], L_p(D)^n)$.

We introduce a metric $d_{\mu}$ on $C([0, T], L_p(D)^n)$ by the formula
\begin{equation*}
    d_{\mu}(u, v) := \sup{\{\, e^{-{\mu}t} \norm{u(t) - v(t)}_{L_p(D)^n} : t \in [0, T] \,\}}, \quad u, v \in C([0, T], L_p(D)^n),
\end{equation*}
where $\mu > 0$ is a positive number (to be chosen later).  The metric $d_{\mu}$ is complete, and equivalent with the metric induced by the usual norm on $C([0, T], L_p(D)^n)$.

\begin{proposition}
\label{prop:contraction}
There exists $\mu > 0$ such that for any $a \in Y$ and any $u_0 \in L_p(D)^n \times L_r((- 1, 0), L_p(D)^n)$ the operator $\mathfrak{G}_{a, u_0}$ on the closed set $\{\, u \in C([0, T],L_p(D)^n) : u(0) = 0 \,\}$ into itself defined by~\eqref{eq:Gothic-G} is a contraction map in the $d_{\mu}$ metric, with the contraction coefficient independent of $1 \le p < \infty$ and $1 \le r \le \infty$.
\end{proposition}
\begin{proof}
Let $a \in Y$ be fixed.

We start by showing that $\mathfrak{G}_{a, u_0}$ takes $\{\, u \in C([0, T],L_p(D)^n) : u(0) = 0 \,\}$ into itself.  Indeed, as a consequence of Lemmas~\ref{lm:mult-estim-1} and~\ref{lm:I_i}, the functions
\begin{gather*}
    \biggl[\, [0, T] \ni t \mapsto \int\limits_{0}^{t} U_{a_0}^0(t, \zeta) \bigl(\mathcal{C}^0_{a}(\zeta) (u(\zeta) + U_{a_0}^0(\zeta, 0) u_0^{(1)}) \bigr) \, \mathrm{d}\zeta \,\biggr], \\
    \biggl[\, [0, T] \ni  t \mapsto \int\limits_{0}^{t} U_{a_0}^0(t, \zeta) \bigl(\mathcal{C}^1_{a}(\zeta) (u_0^{(2)}(\zeta - 1) + \mathbbm{1}_{\! (1, T)}(\zeta)  U_{a_0}^0(\zeta, 0) u_0^{(1)}) \bigr) \, \mathrm{d}\zeta \,\biggr]
\end{gather*}
belong to $C([0, T], L_p(D)^n)$.  The fact that $\mathfrak{G}_{a, u_0}(u)[0] = 0$ is straightforward.

For any $u, v \in C([0, T], L_p(D)^n)$ such that $u(0) = v(0) = 0$ and $0 \le t \le T$, by the results from Proposition~\ref{prop:system-estimates} and Lemma~\ref{lm:mult-estim} under the assumption $p = q$, the following estimate holds:
\begin{equation*}
  \begin{aligned}
  \lVert \mathfrak{G}_{a, u_0}(u)[t] & - \mathfrak{G}_{a, u_0}(v)[t] \rVert_{L_p(D)^n}
  \\
  \le
  {} & \int\limits_{0}^{t} \lVert U_{a_0}^0(t, \zeta) \bigl( \mathcal{C}^0_{a}(\zeta)( u(\zeta) - v(\zeta)) \bigr) \rVert_{L_p(D)^n} \, \mathrm{d}\zeta
  \\
  + {} & \int\limits_{1}^{t} \lVert  U_{a_0}^0(t, \zeta) \bigl(\mathcal{C}_{a}^{1}(\zeta)( u(\zeta - 1) - v(\zeta - 1)) \bigr) \rVert_{L_p(D)^n} \, \mathrm{d}\zeta
  \\
  \le {} & n^2 K M \int\limits_{0}^{t} e^{\gamma (t-\zeta)} \lVert u(\zeta) - v(\zeta) \rVert_{L_p(D)^n} \, \mathrm{d}\zeta  \\
  {} + {} & n^2 K M \int\limits_{0}^{t - 1} e^{\gamma (t-\zeta - 1)} \lVert u(\zeta) - v(\zeta) \rVert_{L_p(D)^n} \, \mathrm{d}\zeta  \\
  \le {} & n^2 K M e^{\gamma T} \int\limits_{0}^{t} \lVert u(\zeta) - v(\zeta) \rVert_{L_p(D)^n} \, \mathrm{d}\zeta  \\
  {} + {} & n^2 K M e^{\gamma T} \int\limits_{0}^{t  - 1} \lVert u(\zeta) - v(\zeta) \rVert_{L_p(D)^n} \, \mathrm{d}\zeta  \\
 \le {} & 2 n^2 K M e^{\gamma T}  \int\limits_{0}^{t} e^{\mu \zeta} d_{\mu}(u, v) \, \mathrm{d}\zeta \\
 = & \biggl(  \frac{2 n^2 K M e^{\gamma T} (e^{\mu t} - 1)}{\mu} \biggr) d_{\mu}(u, v) .
 \end{aligned}
\end{equation*}
As a consequence,
\begin{equation*}
    \begin{aligned}
    d_{\mu}(\mathfrak{G}_{a, u_0}(u), \mathfrak{G}_{a, u_0}(v)) & = \sup{\{\, e^{-{\mu}t} \norm{\mathfrak{G}_{a, u_0}(u)[t] - \mathfrak{G}_{a, u_0}(v)[t]}_{L_p(D)^n} : t \in [0, T] \,\}}
    \\[1ex]
    & \le \sup{\biggl\{\, e^{-{\mu}t} \biggl(  \frac{2 n^2 K M e^{\gamma T} (e^{\mu t} - 1)}{\mu} \biggr) : t \in [0, T] \,\biggr\}} \, d_{\mu}(u, v)
    \\[1ex]
    & \le \frac{2 n^2 K M e^{\gamma T}}{\mu} d_{\mu}(u, v).
    \end{aligned}
    \end{equation*}
It suffices to take $\mu > 2 n^2 K M e^{\gamma T}$.

Observe that in the formula for the contraction coefficient only $K$, $M$ and $\gamma$ appear, which are, by Proposition~\ref{prop:system-estimates} and Lemma~\ref{lm:mult-estim}, independent of $1 \le p < \infty$.  The independence of $1 \le r \le \infty$ is straightforward.
\end{proof}
The Contraction Mapping Principle guarantees the existence and uniqueness of the fixed point of $\mathfrak{G}_{a, u_0}$, which, on adding the term $U_{a_0}^0(t, 0) u_0^{(1)}$, becomes the unique mild $L_p$\nobreakdash-\hspace{0pt}solution of \eqref{main-eq}$_{a}$\textup{+}%
\eqref{main-bc}$_{a}$ on $[0, T]$ satisfying the initial condition \eqref{eq-hull-ic}.  This concludes the proof of Proposition~\ref{prop:delay-existence-hull}.

For Lebesgue\nobreakdash-\hspace{0pt}a.e.\ $t \in (-1, 0)$, $u_p(t; 0, a, u_0)$ is interpreted as $u_0^{(2)}(t)$.

\begin{lemma}
\label{lm:solution-translate}
    For any $1 \le p < \infty$, $1 \le r < \infty$, $a \in Y$, $u_0 \in L_p(D)^n \times L_r((- 1, 0), L_p(D)^n)$ and $\theta \in [0, T)$ there holds
    \begin{equation*}
        \label{eq:integralhull-p-p}
       \begin{aligned}
            u_p(t; 0, a, u_0) = {} & U_{a_0}^0(t, \theta) u_p(\theta; 0, a, u_0) + \int\limits_{\theta}^{t} U_{a_0}^0(t, \zeta) \bigl( \mathcal{C}^0_{a}(\zeta) u_{p}(\zeta; 0, a, u_0) \bigr) \, \mathrm{d}\zeta
            \\
            & + \int\limits_{\theta}^{t} U_{a_0}^0(t, \zeta) \bigl( \mathcal{C}^1_{a}(\zeta) u_{p}(\zeta - 1; 0, a, u_0) \bigr) \, \mathrm{d}\zeta, \quad t \in [\theta, T].
        \end{aligned}
    \end{equation*}
\end{lemma}
\begin{proof}
For $t \in [\theta, T]$, in view of\eqref{eq:cocycle2-2_p} and the Hille theorem (\cite[Thm.~II.2.6]{DiUhl}) there holds
\begin{equation*}
    \begin{aligned}
    u_p(t; 0, a, u_0) = {} & U_{a_0}^0(t, 0) u_0^{(1)} + \int\limits_{0}^{t} U_{a_0}^0(t, \zeta) \bigl( \mathcal{C}^0_{a}(\zeta) u_p(\zeta; 0, a, u_0) \bigr) \, \mathrm{d}\zeta
    \\
    & + \int\limits_{0}^{t} U_{a_0}^0(t, \zeta) \bigl( \mathcal{C}^1_{a}(\zeta) u_p(\zeta - 1; 0, a, u_0) \bigr) \, \mathrm{d}\zeta
    \\
    = {} & U_{a_0}^0(t, \theta) \Biggl( U_{a_0}^0(\theta, 0) u_0^{(1)} +
    \int\limits_{0}^{\theta} U_{a_0}^0(\theta, \zeta) \bigl( \mathcal{C}^0_{a}(\zeta) u_p(\zeta; 0, a, u_0) \bigr) \, \mathrm{d}\zeta
    \\
    & + \int\limits_{0}^{\theta} U_{a_0}^0(\theta, \zeta) \bigl( \mathcal{C}^1_{a}(\zeta) u_p(\zeta - 1; 0, a, u_0) \bigr) \, \mathrm{d}\zeta \Biggr)
    \\
    & +  \int\limits_{\theta}^{t} U_{a_0}^0(t, \zeta) \bigl( \mathcal{C}^0_{a}(\zeta) u_p(\zeta; 0, a, u_0) \bigr) \, \mathrm{d}\zeta
    +  \int\limits_{\theta}^{t} U_{a_0}^0(t, \zeta) \bigl( \mathcal{C}^1_{a}(\zeta) u_p(\zeta - 1; 0, a, u_0) \bigr) \, \mathrm{d}\zeta
    \\
    = {} & U_{a_0}^0(t, \theta) u_p(\theta; 0, a, u_0) + \int\limits_{\theta}^{t} U_{a_0}^0(t, \zeta) \bigl( \mathcal{C}^0_{a}(\zeta) u_{p}(\zeta; 0, a, u_0) \bigr) \, \mathrm{d}\zeta
            \\
            & + \int\limits_{\theta}^{t} U_{a_0}^0(t, \zeta) \bigl( \mathcal{C}^1_{a}(\zeta) u_{p}(\zeta - 1; 0, a, u_0) \bigr) \, \mathrm{d}\zeta.
    \end{aligned}
\end{equation*}
\end{proof}
For $\theta \in [0, T)$ we define
\begin{multline*}
    R(u_p(\cdot; 0, a, u_0))[\theta] := (u_p(\theta; 0, a, u_0), u_p(\cdot + \theta; 0, a, u_0){\restriction}_{(- 1, 0)})
    \\
    \in L_p(D)^n \times L_r((-1, 0), L_p(D)^n).
\end{multline*}
The formula in the statement of Lemma~\ref{lm:solution-translate} can now be written as:
\begin{equation}
\label{eq:translation}
    u_p(t; 0, a, u_0) = u_p(t; \theta, a, R(u_p(\cdot; 0, a, u_0))[\theta]), \quad 0 \le \theta \le t \le T.
\end{equation}

\subsection{Estimates of solutions}
\label{subsect:estimates-mild}

\begin{proposition}
\label{prop:estimate-mild}
For any $a \in Y$, $1 \le p < \infty$, $1 \le r \le \infty$ and $u_0 \in L_p(D)^n \times L_r((-1,0),L_p(D)^n)$  there holds
\begin{multline*}
            \lVert u_p(t;a, u_0) \rVert_{L_p(D)^n}
            \\
             \le \Bigl( Me^{\gamma }(1+n^2K)\exp{\bigl( n^2KMe^{\gamma } \bigr)} \Bigr)^{\lceil T \rceil} \lVert u_0 \rVert_{L_p(D)^n \times L_r((-1,0), L_p(D)^n)}
\end{multline*}
for all $t \in [0, T]$, where $M \ge 1$ and  $\gamma \ge 0$ are the constants from Proposition~\ref{prop:system-estimates}~\textup{(ii)}.
\end{proposition}
\begin{proof}
Let $a \in Y$, $1 \le p < \infty$, $1 \le r \le \infty$ and $u_0 \in L_p(D)^n \times L_r((-1,0),L_p(D)^n)$ be fixed. To keep notation simple we write $u_{p}(t)$ for $u_p(t; 0, a, u_0)$,
and
\begin{equation}
\label{eq:I}
    \begin{aligned}
       I_0(t) & :=  U_{a_0}^0(t, 0) u_0^{(1)},
       \\[1ex]
       I_1(t) & := \int\limits_{0}^{t} U_{a_0}^0(t, \zeta) \bigl(\mathcal{C}^0_{a}(\zeta) u_{p}(\zeta) \bigr) \, \mathrm{d}\zeta,
       \\[1ex]
       I_2(t) & := \int\limits_{0}^{t} U_{a_0}^0(t, \zeta) \bigl(\mathcal{C}^1_{a}(\zeta) u_{p}(\zeta - 1) \bigr) \, \mathrm{d}\zeta.
    \end{aligned}
\end{equation}

First, note that there is an inequality
\begin{align*}
       \norm{u_{p}(t)}_{L_p(D)^n} & \le \norm{I_0(t)}_{L_p(D)^n} + \norm{I_1(t)}_{L_p(D)^n} + \norm{I_2(t)}_{L_p(D)^n}.
\end{align*}
for $t \in [0,T]$. By Proposition~\ref{prop:system-estimates}(ii), there are constants $M  \ge 1$ and $\gamma \ge 0$ such that
\begin{align*}
    \norm{I_0(t)}_{L_p(D)^n}
    &\le M  e^{\gamma T} \norm{u_0^{(1)}}_{L_p(D)^n}\\
    &\le M e^{\gamma T} \norm{u_{0}}_{L_p(D)^n \times L_r((-1,0), L_p(D)^n)}.
\end{align*}
Furthermore the inequalities
\begin{align*}
    \norm{I_1(t)}_{L_p(D)^n} & \le  \int\limits_{0}^{t} \norm{ U_{a_0}^0(t, \zeta) \bigl(\mathcal{C}^0_{a}(\zeta) u_{p}(\zeta)\bigr)}_{L_p(D)^n} \, \mathrm{d}\zeta \\
    & \le n^2KM \int\limits_{0}^{t} e^{\gamma(t-\zeta)} \norm{ \mathcal{C}^0_{a}(\zeta) u_{p}(\zeta)}_{L_p(D)^n} \, \mathrm{d}\zeta \\
    & \le n^2KM e^{\gamma T}  \int\limits_{0}^{t} \norm{u_{p}(\zeta)}_{L_p(D)^n} \, \mathrm{d}\zeta
\end{align*}
and
\begin{align*}
    \norm{I_2(t)}_{L_p(D)^n} & \le \int\limits_{0}^{t} \norm{U_{a_0}^0(t, \zeta ) \bigl(\mathcal{C}^1_{a}(\zeta ) u_{p}(\zeta  - 1)\bigr)}_{L_p(D)^n}  \, \mathrm{d}\zeta  \\
    & \le n^2KM e^{\gamma T} \int\limits_{0}^{t}  \norm{u_{p}(\zeta - 1)}_{L_p(D)^n} \, \mathrm{d}\zeta
\end{align*}
hold. Therefore we can write that
\begin{align*}
     \norm{u_p(t)}_{L_p(D)^n}& \le M e^{\gamma T} \norm{u_{0}}_{L_p(D)^n \times L_r((- 1, 0), L_p(D)^n)} \\
     &+  n^2KM e^{\gamma T}  \int\limits_{0}^{t}  \norm{u_{p}(\zeta - 1)}_{L_p(D)^n} \, \mathrm{d}\zeta \\
     & + n^2KM e^{\gamma T}  \int\limits_{0}^{t} \norm{u_{p}(\zeta)}_{L_p(D)^n} \, \mathrm{d}\zeta.
\end{align*}
As the function
\begin{align*}
    \biggl[\,[ 0, T ] \ni t \mapsto  & M e^{\gamma T} \norm{u_{0}}_{L_p(D)^n \times L_r((-1,0), L_p(D)^n)}
    \\
    & +  n^2KM e^{\gamma T}  \int\limits_{0}^{t}  \norm{u_{p}(\zeta - 1)}_{L_p(D)^n} \, \mathrm{d}\zeta \in \RR \, \biggr]
\end{align*}
is nondecreasing, an application of the standard Gronwall inequality gives that
\begin{align*}
    & \norm{u_p(t)}_{L_p(D)^n}
    \\
     & \le \Biggl( M e^{\gamma T} \norm{u_{0}}_{L_p(D)^n \times L_r((-1,0), L_p(D)^n)} +  n^2KM e^{\gamma T}  \int\limits_{0}^{t}  \norm{u_{p}(\zeta - 1)}_{L_p(D)^n} \, \mathrm{d}\zeta \Biggr)  \\
     &\times \exp{ \biggl(  n^2KM e^{\gamma T}  \int\limits_{0}^{t}  \, \mathrm{d}\zeta \biggr) }.
\end{align*}
 For $T \in (0, 1]$, applying the Hölder inequality we obtain that
\begin{multline*}
    \int\limits_{0}^{t} \norm{u_{p}(\zeta - 1)}_{L_p(D)^n} \, \mathrm{d}\zeta \le \int\limits_{0}^{1} \norm{u_{p}(\zeta - 1)}_{L_p(D)^n} \, \mathrm{d}\zeta
    \\
    \le \norm{u_{0}^{(2)}}_{L_r((-1,0), L_p(D)^n)} \le \norm{u_{0}}_{L_p(D)^n \times L_r((-1,0), L_p(D)^n)}.
\end{multline*}
Consequently, for such $t$ the main inequality holds with
\begin{equation*}
 Me^{\gamma }(1+n^2K)\exp{\bigl( n^2KMe^{\gamma } \bigr)}.
\end{equation*}

If $T > 1$ we can repeat the above estimation to obtain the desired result.
\end{proof}

\subsection{Regularizing action of the solution operator}
\label{subsect:regularization}

\begin{proposition}
\label{prop:regularization}
    For any $a \in Y$, $1 \le p \le q < \infty$, $1 < r \le \infty$ satisfying
    \begin{equation}
    \label{eq:ass-p-q}
        \frac{N}{2}\biggl(\frac{1}{p}- \frac{1}{q}\biggr) < \frac{1}{r'},
    \end{equation}
    and $u_0 \in L_p(D)^n \times L_r((- 1, 0), L_p(D)^n)$ there holds:

        $u_p(\cdot;0, a, u_0)$ is a continuous function from $(0, T]$ into $L_q(D)^n$, satisfying, for each $\theta \in (0, T)$, the equation
       \begin{equation}
       \label{eq:integralhull-p-q}
       \begin{aligned}
            u_p(t; 0, a, u_0) = {} & U_{a_0}^0(t, \theta) u_p(\theta;0, a, u_0) + \int\limits_{\theta}^{t} U_{a_0}^0(t, \zeta) \bigl( \mathcal{C}^0_{a}(\zeta) u_{p}(\zeta; 0, a, u_0)\bigr) \, \mathrm{d}\zeta
            \\
            & + \int\limits_{\theta}^{t} U_{a_0}^0(t, \zeta) \bigl( \mathcal{C}^1_{a}(\zeta) u_{p}(\zeta - 1; 0, a, u_0)\bigr) \, \mathrm{d}\zeta, \quad t \in [\theta, T],
        \end{aligned}
        \end{equation}
        in the $L_q(D)^n$\nobreakdash-\hspace{0pt}sense.
\end{proposition}
\begin{proof}
    Fix $a \in Y$, $1 \le p \le q < \infty$, $1 < r \le \infty$, and $u_0 \in L_p(D)^n \times L_r((- 1, 0), L_p(D)^n)$, where we assume additionally that \eqref{eq:ass-p-q}is satisfied.

    It follows from Lemmas~\ref{lm:mult-estim-1} and~\ref{lm:I_i} that $I_1, I_2 \in C([0, T], L_q(D)^n)$, where $I_i$, $i = 0, 1, 2$, are defined in~\eqref{eq:I}. As $I_0$ is, by Propositions~\ref{prop:higher_orders-skew_product-p} and~\ref{prop:higher_orders-skew_product-p-q-estimates}(i), continuous from $(0, T]$ to $L_q(D)^n$, the required continuity of $u_p(\cdot; 0, a, u_0)$ follows.

    \eqref{eq:integralhull-p-q} is satisfied in the $L_p(D)^n$\nobreakdash-\hspace{0pt}sense.  Since all the Bochner integrals appearing there are convergent in the $L_q(D)^n$\nobreakdash-\hspace{0pt}norm, the satisfaction of~\eqref{eq:integralhull-p-q} in the $L_q(D)^n$\nobreakdash-\hspace{0pt}sense follows.
\end{proof}

\begin{proposition}
\label{prop:regularization_up_to_infty}
    For any $a \in Y$, $1 \le p  < \infty$, $1 < r \le \infty$ satisfying
    \begin{equation}
        \frac{N}{2p} < \frac{1}{r'},
    \end{equation}
    and $u_0 \in L_p(D)^n \times L_r((- 1, 0), L_p(D)^n)$ there holds: $u_p(t;0, a, u_0) \in L_{\infty}(D)^n$ for all $t \in (0, T]$, and satisfies, for each $\theta \in (0, T)$, the equation
       \begin{equation}
       \label{eq:integralhull-q-infty}
       \begin{aligned}
            u_p(t; 0, a, u_0) = {} & U_{a_0}^0(t, \theta) u_p(\theta;0, a, u_0)
            \\
            & + \textup{G -}\hspace{-1pt}\int\limits_{\theta}^{t} U_{a_0}^0(t, \zeta) \bigl( \mathcal{C}^0_{a}(\zeta) u_{p}(\zeta; 0, a, u_0)\bigr)\, \mathrm{d}\zeta
            \\
            & + \textup{G -}\hspace{-1pt}\int\limits_{\theta}^{t} U_{a_0}^0(t, \zeta) \bigl( \mathcal{C}^1_{a}(\zeta) u_{p}(\zeta - 1; 0, a, u_0)\bigr) \, \mathrm{d}\zeta, \quad t \in [\theta, T],
        \end{aligned}
        \end{equation}
        in the $L_{\infty}(D)^n$\nobreakdash-\hspace{0pt}sense.
\end{proposition}

\begin{proof}
   It is clear from Proposition~\ref{prop:system-estimates} that $ U_{a_0}^0(t, \theta) u_p(\theta;0, a, u_0)\in L_{\infty}(D)^n$ and the mappings
\begin{gather*}
    \Bigl[\, [\theta, t] \ni \zeta \mapsto  U_{a_0}^0(t, \zeta) \bigl( \mathcal{C}^0_{a}(\zeta) u_{p}(\zeta; 0, a, u_0)\bigr) \,\Bigr], \\
    \Bigl[\, [\theta, t] \ni \zeta \mapsto U_{a_0}^0(t, \zeta) \bigl( \mathcal{C}^1_{a}(\zeta) u_{p}(\zeta-1; 0, a, u_0)\bigr) \,\Bigr]
\end{gather*}
    are $L_{\infty}(D)^n$\nobreakdash-\hspace{0pt}valued. Furthermore, we claim they satisfy the assumptions of the Lemma~\ref{lm:L-_infty}. Indeed,
    \begin{gather*}
    \Bigl[\, [\theta, t] \ni \zeta \mapsto  \mathcal{C}^0_{a}(\zeta) u_{p}(\zeta; 0, a, u_0)\in L_{p}(D)^n \,\Bigr]\in L_{r}((\theta,t),L_p(D)^n), \\
    \Bigl[\, [\theta, t] \ni \zeta \mapsto  \mathcal{C}^1_{a}(\zeta) u_{p}(\zeta-1; 0, a, u_0)\in L_{p}(D)^n \,\Bigr] \in L_{r}((\theta,t),L_p(D)^n)
\end{gather*}
   since $u_{p}(\cdot; 0, a, u_0)$ is continuous on $[\theta,T]$ for $1 < p < \infty$ (or bounded and continuous on $(\theta,T]$ when $p = 1$), hence together with Lemma~\ref{lm:mult-estim}  we conclude the claim.  It suffices now to check that for any $v\in L_1(D)^n$ and $t \in [\theta, T]$ equality
\begin{multline*}
      \langle  v, u_p(t; 0, a, u_0) \rangle_{L_1(D)^n,L_{\infty}(D)^n}  =   \langle v, U_{a_0}^0(t, \theta)u_p(\theta;0, a, u_0) \rangle_{L_1(D)^n,L_{\infty}(D)^n}
            \\
            + \Big\langle v, \textup{G -}\hspace{-1pt}\int\limits_{\theta}^{t} U_{a_0}^0(t, \zeta) \bigl( \mathcal{C}^0_{a}(\zeta) u_{p}(\zeta; 0, a, u_0)\bigr)\, \mathrm{d}\zeta \Big\rangle_{L_1(D)^n,L_{\infty}(D)^n}
            \\
            + \Big\langle v, \textup{G -}\hspace{-1pt}\int\limits_{\theta}^{t} U_{a_0}^0(t, \zeta) \bigl( \mathcal{C}^1_{a}(\zeta) u_{p}(\zeta - 1; 0, a, u_0)\bigr) \, \mathrm{d}\zeta \Big\rangle_{L_1(D)^n,L_{\infty}(D)^n}
\end{multline*}
holds, or equivalently, straight from the definition of the Gel$\,\textprime$fand integral,
\begin{multline*}
      \langle  v, u_p(t; 0, a, u_0) \rangle_{L_1(D)^n,L_{\infty}(D)^n}  =   \langle v, U_{a_0}^0(t, \theta)u_p(\theta;0, a, u_0) \rangle_{L_1(D)^n,L_{\infty}(D)^n}
            \\
            +  \int\limits_{\theta}^{t} \big\langle v,  U_{a_0}^0(t, \zeta) \bigl( \mathcal{C}^0_{a}(\zeta) u_{p}(\zeta; 0, a, u_0)\bigr) \big\rangle_{L_1(D)^n,L_{\infty}(D)^n} \, \mathrm{d}\zeta
            \\
            +  \int\limits_{\theta}^{t} \big\langle v,  U_{a_0}^0(t, \zeta) \bigl( \mathcal{C}^1_{a}(\zeta) u_{p}(\zeta-1; 0, a, u_0)\bigr) \big\rangle_{L_1(D)^n,L_{\infty}(D)^n} \, \mathrm{d}\zeta.
\end{multline*}
Let $(v_m)_{m=1}^{\infty}\subset L_{p'}(D)^n$ be a sequence of approximates of $v$ in the $L_1(D)^n$\nobreakdash-\hspace{0pt}norm. By Lemma~\ref{lm:solution-translate} and the Hille theorem (\cite[Thm.~II.2.6]{DiUhl}) we have
\begin{multline*}
      \langle u_p(t; 0, a, u_0), v_m \rangle_{L_{p}(D)^n,L_{p'}(D)^n}  =   \langle U_{a_0}^0(t, \theta)u_p(\theta;0, a, u_0),  v_m \rangle_{L_{p}(D)^n,L_{p'}(D)^n}
            \\
            +  \int\limits_{\theta}^{t} \big\langle  U_{a_0}^0(t, \zeta) \bigl( \mathcal{C}^0_{a}(\zeta) u_{p}(\zeta; 0, a, u_0)\bigr),  v_m \big\rangle_{L_{p}(D)^n,L_{p'}(D)^n}   \, \mathrm{d}\zeta
            \\
            +  \int\limits_{\theta}^{t} \big\langle  U_{a_0}^0(t, \zeta) \bigl( \mathcal{C}^1_{a}(\zeta) u_{p}(\zeta-1; 0, a, u_0)\bigr),  v_m \big\rangle_{L_{p}(D)^n,L_{p'}(D)^n}   \, \mathrm{d}\zeta.
\end{multline*}
We can reinterpret $\langle\cdot,\cdot\rangle$:
\begin{multline*}
      \langle  v_m, u_p(t; 0, a, u_0) \rangle_{L_{1}(D)^n,L_{\infty}(D)^n}  =   \langle v_m, U_{a_0}^0(t, \theta)u_p(\theta;0, a, u_0) \rangle_{L_{1}(D)^n,L_{\infty}(D)^n}
            \\
            +  \int\limits_{\theta}^{t} \big\langle v_m,  U_{a_0}^0(t, \zeta) \bigl( \mathcal{C}^0_{a}(\zeta) u_{p}(\zeta; 0, a, u_0)\bigr) \big\rangle_{L_{1}(D)^n,L_{\infty}(D)^n}   \, \mathrm{d}\zeta
            \\
            +  \int\limits_{\theta}^{t} \big\langle v_m,  U_{a_0}^0(t, \zeta) \bigl( \mathcal{C}^1_{a}(\zeta) u_{p}(\zeta-1; 0, a, u_0)\bigr) \big\rangle_{L_{1}(D)^n,L_{\infty}(D)^n}   \, \mathrm{d}\zeta.
\end{multline*}
Therefore, it remains to notice that
\begin{multline*}
     \langle v_m, U_{a_0}^0(t, \theta)u_p(\theta;0, a, u_0) \rangle_{L_{1}(D)^n,L_{\infty}(D)^n} \\ \to  \langle v, U_{a_0}^0(t, \theta)u_p(\theta;0, a, u_0) \rangle_{L_{1}(D)^n,L_{\infty}(D)^n},
\end{multline*}
\begin{multline*}
    \int\limits_{\theta}^{t} \big\langle v_m,  U_{a_0}^0(t, \zeta) \bigl( \mathcal{C}^0_{a}(\zeta) u_{p}(\zeta; 0, a, u_0)\bigr) \big\rangle_{L_{1}(D)^n,L_{\infty}(D)^n}   \, \mathrm{d}\zeta \\ \to \int\limits_{\theta}^{t} \big\langle v,  U_{a_0}^0(t, \zeta) \bigl( \mathcal{C}^0_{a}(\zeta) u_{p}(\zeta; 0, a, u_0)\bigr) \big\rangle_{L_{1}(D)^n,L_{\infty}(D)^n}   \, \mathrm{d}\zeta,
\end{multline*}
\begin{multline*}
    \int\limits_{\theta}^{t} \big\langle v_m,  U_{a_0}^0(t, \zeta) \bigl( \mathcal{C}^0_{a}(\zeta) u_{p}(\zeta-1; 0, a, u_0)\bigr) \big\rangle_{L_{1}(D)^n,L_{\infty}(D)^n}   \, \mathrm{d}\zeta \\ \to
    \int\limits_{\theta}^{t} \big\langle v,  U_{a_0}^0(t, \zeta) \bigl( \mathcal{C}^1_{a}(\zeta) u_{p}(\zeta-1; 0, a, u_0)\bigr) \big\rangle_{L_{1}(D)^n,L_{\infty}(D)^n}   \, \mathrm{d}\zeta
\end{multline*}
as $m\to \infty$.
\end{proof}

It should be remarked that the equality
\begin{equation}
\label{eq:translation-p-q}
    u_p(t; 0, a, u_0) = u_p(t; \theta, a, R(u_p(\cdot; 0, a, u_0))[\theta]), \quad 0 \le \theta \le t \le T,
\end{equation}
need not hold in the $L_q(D)^n$\nobreakdash-\hspace{0pt}sense for $T \in (0, 1]$ (since $R(u_p(\cdot; 0, a, u_0))[\theta]$ need not belong to $L_q(D)^n \times L_r((-1,0), L_q(D)^n)$).  However, \eqref{eq:translation-p-q} holds in the $L_q(D)^n$\nobreakdash-\hspace{0pt}sense for $T > 1$.

\begin{proposition}
  \label{prop:delay-estimates-1}
  There is a constant $\overline{M} > 0$ such that
    \begin{equation*}
      \lVert u_p(t;0, a, u_0) \rVert_{L_q(D)^n} \le \overline{M} t^{-\tfrac{N}{2}\bigl(\tfrac{1}{p}- \tfrac{1}{q}\bigr)} \lVert u_0 \rVert_{L_p(D)^n \times L_r((-1,0), L_p(D)^n)}, \quad t \in (0, T],
    \end{equation*}
    for all $a \in Y$, provided that $1 \le p \le q < \infty$ are such that
    \begin{equation}
    \label{ineq:auxiliary}
        \frac{N}{2} \left( \frac{1}{p} - \frac{1}{q} \right) < \frac{1}{r'}.
    \end{equation}
\end{proposition}
\begin{proof}
    We consider first the case $T \in (0, 1]$, and write $u_p(\cdot)$ for $u_p(\cdot; 0, a, u_0)$.

    Applying Propositions~\ref{prop:system-estimates} and Lemma~\ref{lm:mult-estim} to~\eqref{eq:integral} we obtain
    \begin{equation}
    \label{eq:p-to-q}
        \begin{aligned}
        \norm{u_p(t)}_{L_q(D)^n} \le {} & M t^{-\tfrac{N}{2} \bigl(\tfrac{1}{p} - \tfrac{1}{q}\bigr)} e^{\gamma t} \norm{u_0^{(1)}}_{L_p(D)^n}
        \\
        & + n^2 K M  \int\limits_{0}^{t} (t-\zeta)^{-\tfrac{N}{2} \bigl(\tfrac{1}{p} - \tfrac{1}{q}\bigr)} e^{\gamma (t - \zeta)} \norm{u_p(\zeta)}_{L_p(D)^n} \, \mathrm{d}\zeta
        \\
        & + n^2 K M \int\limits_{0}^{t} (t-\zeta)^{-\tfrac{N}{2} \bigl(\tfrac{1}{p} - \tfrac{1}{q}\bigr)} e^{\gamma (t -\zeta)} \norm{u_p(\zeta - 1)}_{L_p(D)^n} \, \mathrm{d}\zeta, \quad t \in [0, 1],
        \end{aligned}
    \end{equation}
    By Proposition~\ref{prop:estimate-mild},
    \begin{multline*}
     \int\limits_{0}^{t} (t-\zeta)^{-\tfrac{N}{2} \bigl(\tfrac{1}{p} - \tfrac{1}{q}\bigr)} e^{\gamma (t - \zeta)} \norm{u_p(\zeta)}_{L_p(D)^n} \, \mathrm{d}\zeta
     \\
     \le    \frac{e^{\gamma}  t^{1 -\tfrac{N}{2} \bigl(\tfrac{1}{p} - \tfrac{1}{q}\bigr)}}{1 -\tfrac{N}{2} \bigl(\tfrac{1}{p} - \tfrac{1}{q}\bigr) } Me^{\gamma }(1+n^2K)\exp{\bigl( n^2KMe^{\gamma } \bigr)} \lVert u_0 \rVert_{L_p(D)^n \times L_r((-1,0), L_p(D)^n)} \\
      \le    \frac{  Me^{2\gamma }(1+n^2K)}{1 -\tfrac{N}{2} \bigl(\tfrac{1}{p} - \tfrac{1}{q}\bigr) } \exp{\bigl( n^2KMe^{\gamma} \bigr)} \lVert u_0 \rVert_{L_p(D)^n \times L_r((-1,0), L_p(D)^n)}.
    \end{multline*}
    Regarding the third term, we see that $\zeta \to (t - \zeta)^{-\tfrac{N}{2} \bigl(\tfrac{1}{p} - \tfrac{1}{q}\bigr)}$ belongs to $L_{r'}((0, t))$, whereas $\zeta \to \norm{u_p(\zeta - 1)}_{L_p(D)^n} = \norm{u_0^{(2)}(\zeta - 1)}_{L_p(D)^n}$ belongs to $L_{r}((0, t))$.  By the H\"older inequality,
    \begin{align*}
    & \int\limits_{0}^{t} (t-\zeta)^{-\tfrac{N}{2} \bigl(\tfrac{1}{p} - \tfrac{1}{q}\bigr)} e^{\gamma (t - \zeta)} \norm{u_0^{(2)}(\zeta - 1)}_{L_p(D)^n} \, \mathrm{d}\zeta
    \\
    & \le e^{\gamma} \biggl( \int_{0}^{t} (t-\zeta)^{-\tfrac{N r'}{2} \bigl(\tfrac{1}{p} - \tfrac{1}{q}\bigr)} \, \mathrm{d}\zeta \biggr)^{\!\!1/r'} \, \lVert u_0^{(2)} \rVert_{L_r((-1, 0), L_p(D)^n)}
    \\
    & \le \frac{e^{\gamma}}{\bigl( 1 -\tfrac{N}{2} (\tfrac{1}{p} - \tfrac{1}{q} ) r' \bigr)^{\!\! 1/r'}} \lVert u_0^{(2)} \rVert_{L_r((-1, 0), L_p(D)^n)}.
    \end{align*}
    Putting together the results obtained we arrive at the inequality
    \begin{equation*}
    \lVert u_p(t;0, a, u_0) \rVert_{L_q(D)^n}  \le \overline{M} t^{-\tfrac{N}{2} \bigl(\tfrac{1}{p} - \tfrac{1}{q}\bigr)} \lVert u_0 \rVert_{L_p(D)^n \times L_r((-1,0), L_p(D)^n)}
    \end{equation*}
    for $t \in (0, 1]$, where
    \begin{equation*}
        \overline{M} = M  e^{\gamma} \biggl(  1 +  \frac{  e^{\gamma }(1+n^2K) n^2 K M}{1 -\tfrac{N}{2} \bigl(\tfrac{1}{p} - \tfrac{1}{q}\bigr) } \exp{\bigl( n^2KMe^{\gamma} \bigr)}   +    \frac{n^2 K}{\bigl( 1 -\tfrac{N}{2} (\tfrac{1}{p} - \tfrac{1}{q} ) r' \bigr)^{\!\! 1/r'}}    \biggr).
    \end{equation*}
In view of~\eqref{eq:translation-p-q} the statement follows for $T > 0$ by Proposition~\ref{prop:estimate-mild}.
\end{proof}

In fact, the above proposition can be extended to the case $q=\infty$ without changes in the constant $\overline{M}$.
\begin{proposition} \label{prop:delay-estimates-infty-q}
  There holds
    \begin{equation*}
      \lVert u_p(t;0, a, u_0) \rVert_{L_{\infty}(D)^n} \le \overline{M} t^{-\tfrac{N}{2p}} \lVert u_0 \rVert_{L_p(D)^n \times L_r((-1,0), L_p(D)^n)}, \quad t \in (0, T],
    \end{equation*}
    for all $a \in Y$, provided that $1 \le p  < \infty$ are such that
    \begin{equation}
    \label{ineq:auxiliary-p}
       \frac{N}{2p}<\frac{1}{r'} \ ,
    \end{equation}
    where $\overline{M}$ is the constant as in Proposition~\ref{prop:delay-estimates-1}.
\end{proposition}
\begin{proof}
     For fixed $1\le p< \infty$ and any $p<q<\infty$ by Proposition~\ref{prop:delay-estimates-1} we can write
     \begin{equation*}
      \lVert u_p(t;0, a, u_0) \rVert_{L_q(D)^n} \le \overline{M} t^{-\tfrac{N}{2}\bigl(\tfrac{1}{p}- \tfrac{1}{q}\bigr)} \lVert u_0 \rVert_{L_p(D)^n \times L_r((-1,0), L_p(D)^n)}.
    \end{equation*}
    Moreover, due to the finiteness of the measure of $D$, we have $ \lVert u_p(t;0, a, u_0) \rVert_{L_q(D)^n} \to    \lVert u_p(t;0, a, u_0) \rVert_{L_{\infty}(D)^n}$ as $q\to\infty$.
    Therefore the inequality between limits remains, and we will obtain
    \begin{equation*}
          \lVert u_p(t;0, a, u_0) \rVert_{L_{\infty}(D)^n} \le \overline{M} t^{-\tfrac{N}{2p}} \lVert u_0 \rVert_{L_p(D)^n \times L_r((-1,0), L_p(D)^n)},
    \end{equation*}
    for each $t \in (0, T]$.
\end{proof}

For an arbitrary $q \in (p, \infty)$, we take $p = p_0 < p_1 < \ldots < p_{m_0 - 1} < p_{m_0} = q$ such that
\begin{equation*}
    \frac{N}{2}\biggl(\frac{1}{p_{m}}- \frac{1}{p_{m+1}}\biggr) \le \frac{1}{2r'}, \quad m = 0, \ldots, m_0 - 1.
\end{equation*}
Observe that, since $0 < 1/q < 1/p \le 1$, the number $m_0$ depends only on $r$, and is independent of $1 \le p < q < \infty$.  Indeed, we can take
\begin{equation*}
    m_0 = \lceil Nr' \rceil.
\end{equation*}
In~particular, if $m_0$ is good for some $r_1 > 1$ then it is good for any $r_2 \in [r_1, \infty)$.  Consequently we have the following.
\begin{proposition}
\label{prop:regularization-q}
    For any $1 < r_0 < \infty$ there exist $\Theta = \lceil \tfrac{N r_0}{r_0 - 1} \rceil$ such that, assuming $T > \Theta + 1$, for any $a \in Y$, any $1 \le p < q < \infty$, any $r_0 \le r \le \infty$ and any $u_0 \in L_p(D)^n \times L_r((- 1, 0), L_p(D)^n)$ there holds:

        $u_p(\cdot;0, a, u_0)$ is a continuous function from $[\Theta, T]$ into $L_q(D)^n$, satisfying, for each $\theta \in [\Theta + 1, T)$, the equation
       \begin{equation}
       \label{eq:integralhull-p-q-1}
       \begin{aligned}
            u_p(t; 0, a, u_0) = {} & U_{a_0}^0(t, \theta) u_p(\theta;0, a, u_0) + \int\limits_{\theta}^{t} U_{a_0}^0(t, \zeta) \bigl( \mathcal{C}^0_{a}(\zeta) u_{p}(\zeta; 0, a, u_0) \bigr) \, \mathrm{d}\zeta
            \\
            & + \int\limits_{\theta}^{t} U_{a_0}^0(t, \zeta) \bigl( \mathcal{C}^1_{a}(\zeta) u_{p}(\zeta - 1;0, a, u_0) \bigr) \, \mathrm{d}\zeta, \quad t \in [\theta, T],
        \end{aligned}
        \end{equation}
        in the $L_q(D)^n$\nobreakdash-\hspace{0pt}sense.
\end{proposition}
We can say that the equality
\begin{equation*}
    u_p(t; 0, a, u_0) = u_p(t; \theta, a, R(u_p(\cdot; 0, a, u_0))[\theta]), \quad \Theta + 1 \le \theta \le t \le T
\end{equation*}
holds in the $L_q(D)^n$\nobreakdash-\hspace{0pt}sense.

Proposition~\ref{prop:regularization-q} allows us to suppress the symbol $p$ in the notation for solution, at~least for $t$ sufficiently large.  From now on, we will write $u(\cdot; 0, a, u_0)$.

\smallskip
We proceed to obtaining estimates for the regularizing action from $p$ to $q$ where $1 \le p < q < \infty$ are arbitrary.
\begin{proposition}
  \label{prop:delay-estimates-2}
  For any $1 < r_0 \le \infty$ there exist $\Theta = \lceil \tfrac{N r_0}{r_0 - 1} \rceil$ and $\widetilde{M} > 0$ such that, assuming $T > \Theta + 1$,
  \begin{equation*}
      \lVert u(t; 0, a, u_0) \rVert_{L_q(D)^n} \le \widetilde{M} \lVert u_0 \rVert_{L_p(D)^n \times L_r((-1,0), L_p(D)^n)}
  \end{equation*}
  for all $1 \le p < q < \infty$, $r_0 \le r \le \infty$, $a \in Y$, $u_0 \in L_p(D)^n \times L_r((- 1, 0), L_p(D)^n)$ and $t \in [\Theta + 1, T]$.
\end{proposition}
\begin{proof}
    This is a consequence of Propositions~\ref{prop:delay-estimates-1} and~\ref{prop:regularization-q}.
\end{proof}

\begin{proposition}
 \label{prop:delay-estimates-infty}
 For any $1 < r_0 \le \infty$ there exist $\Theta = \lceil \tfrac{N r_0}{r_0 - 1} \rceil$ and $\widetilde{M} > 0$ such that, assuming $T > \Theta + 1$,
  \begin{equation*}
      \lVert u(t; 0, a, u_0) \rVert_{L_{\infty}(D)^n} \le \widetilde{M} \lVert u_0 \rVert_{L_p(D)^n \times L_r((-1,0), L_p(D)^n)}
  \end{equation*}
  for all $1 \le p < \infty$, $r_0 \le r \le \infty$, $a \in Y$, $u_0 \in L_p(D)^n \times L_r((- 1, 0), L_p(D)^n)$ and $t \in [\Theta + 1, T]$.
\end{proposition}
\begin{proof}
    We let $q \to \infty$ as in the proof of Proposition~\ref{prop:delay-estimates-infty-q}.
\end{proof}

\subsection{Compactness}
\label{subsect:compact}
\begin{lemma}
\label{lm:compact-aux-at_once}
    For any $0 < t \le 1$, $1 \le  p \le q < \infty$, $1 < r \le \infty$ satisfying
    \begin{equation*}
        \frac{N}{2}\biggl(\frac{1}{p}- \frac{1}{q}\biggr) < \frac{1}{r'},
    \end{equation*}
    the set
    \begin{equation*}
        \biggl\{ \int\limits_{0}^{t} U_{a_0}^0(t, {\zeta}) u({\zeta}) \, \dd{\zeta} : a_0 \in Y_0,  \norm{u}_{L_{r}((0, t), L_p(D)^n)} \le 1 \biggr\}
    \end{equation*}
    is precompact in $L_q(D)^n$.
\end{lemma}
\begin{proof}
    Let $t$, $p$, $q$ and $r$ be as in the statement.  Fix, for the moment, $\theta \in (0, t)$.  We show first that the set
    \begin{equation*}
        \biggl\{ \int\limits_{0}^{\theta} U_{a_0}^0(t, {\zeta}) u({\zeta}) \, \dd{\zeta} : a_0 \in Y_0, \norm{u}_{L_{1}((0, t), L_p(D)^n)} \le 1 \biggr\}
    \end{equation*}
    is precompact in $L_q(D)^n$.  Let $C$ stand for the closed convex hull of the closure in $L_q(D)^n$ of the set $\{\, U_{a_0}^0(t, {\zeta}) \tilde{u}: a \in Y_0, \zeta \in [0, \theta], \ \norm{\tilde{u}}_{L_p(D)^n} \le 1 \,\}$.  By Proposition~\ref{prop:higher_order-compactness} and Mazur's theorem (\cite[Thm.~II.2.12]{DiUhl}), $C$ is compact.  Assume that $u \in L_{1}((0, \theta), L_p(D)^n)$ is a simple function: $u = \sum\limits_{k = 1}^r \mathbbm{1}_{G_k} \tilde{u}_k$ where $\tilde{u}_k \in L_p(D)^n$ and pairwise disjoint Borel sets $G_k$ are such that $\bigcup\limits_{k = 1}^r G_k = (0, \theta)$.  Then
    \begin{align*}
        \int\limits_{0}^{\theta} U_{a_0}^0(t, {\zeta}) u({\zeta}) \, \dd{\zeta} & \in \sum_{k = 1}^r \Bigl( \lambda(G_k) \norm{\tilde{u}_k}_{L_p(D)^n} C \Bigr)
        \\
        & = \Bigl( \sum_{k = 1}^r \lambda(G_k) \norm{\tilde{u}_k}_{L_p(D)^n} \Bigr) C = \norm{u}_{L_{1}((0, \theta), L_p(D)^n)} C,
    \end{align*}
   where $\lambda(\cdot)$ stands for the one-dimensional Lebesgue measure.
   For each $u \in L_{1}((0, \theta), L_p(D)^n)$ there exists a sequence $(u_{l})_{l = 1}^{\infty}$ of simple functions such that $\lim\limits_{l \to \infty} \norm{u_{l} - u}_{L_{1}((0, \theta), L_p(D)^n)} = 0$ and $\lim\limits_{l \to \infty} \norm{u_{l}(\zeta) - u(\zeta)}_{L_p(D)^n} = 0$ for $\lambda$\nobreakdash-\hspace{0pt}a.e.\ $\zeta \in (0, \theta)$ (see~\cite[X.2]{AmEsch}).  It follows from Proposition~\ref{prop:norm_continuity-all} that $\lim\limits_{l \to \infty} \norm{U_{a_0}^0(t, {\zeta}) \bigl(u_{l}(\zeta) - u(\zeta)\bigr)}_{L_q(D)^n} = 0$ for $\lambda$\nobreakdash-\hspace{0pt}a.e.\ $\zeta \in (0, \theta)$.  The estimate in Proposition~\ref{prop:system-estimates}(ii), together with our assumptions on $p$, $q$ and $r$, give us that (remember that $0 < \theta < t \le 1$)
   \begin{align*}
       \norm{U_{a_0}^0(t, \zeta) u_{l}(\zeta)}_{L_q(D)^n}
       & \le M (t - \theta)^{- \frac{N}{2} (\frac{1}{p} - \frac{1}{q})} e^{\gamma (t - \theta)} \norm{u_{l}(\zeta)}_{L_p(D)^n}
       \\
       & \le M (t - \theta)^{- \frac{1}{r'}} e^{\gamma t} \norm{u_{l}(\zeta)}_{L_p(D)^n}
   \end{align*}
   for $\lambda$\nobreakdash-\hspace{0pt}a.e.\ $\zeta \in (0, \theta)$.  Since $\lim\limits_{l \to \infty} \norm{u_{l} - u}_{L_{1}((0, \theta), L_p(D)^n)} = 0$, it follows that \newline $\norm{u_{l}}_{L_{1}((0, \theta), L_p(D)^n)} \le 2$ for sufficiently large $l$.  We have thus obtained that the $L_{1}((0, \theta), L_q(D)^n)$\nobreakdash-\hspace{0pt}norms of the functions $[\, \zeta \mapsto U_{a_0}^0(t, \zeta) u_{l}(\zeta) \,]$ are bounded uniformly in $l$, which allows us to apply the Dominated Convergence Theorem for Bochner integrals (see~\cite[Thm.~X.3.12]{AmEsch}) to conclude that
   \begin{equation*}
       \int\limits_{0}^{\theta} U_{a_0}^0(t, {\zeta}) u_{l}({\zeta}) \, \dd{\zeta} \to \int\limits_{0}^{\theta} U_{a_0}^0(t, {\zeta}) u({\zeta}) \, \dd{\zeta} \quad \text{in the } L_q(D)^n \text{ sense},
   \end{equation*}
   consequently
   \begin{align*}
       \int\limits_{0}^{\theta} U_{a_0}^0(t, {\zeta}) u({\zeta}) \, \dd{\zeta} & \in \Bigl(\lim_{l \to \infty} \norm{u_{l}}_{L_{1}((0, \theta), L_p(D)^n)}\Bigr) C
       \\
       & = \norm{u}_{L_{1}((0, \theta), L_p(D)^n)} C.
   \end{align*}
   In~particular, we have that the set
    \begin{equation*}
        \biggl\{ \int\limits_{0}^{\theta} U_{a_0}^0(t, {\zeta}) u({\zeta}) \, \dd{\zeta} : a_0 \in Y_0, \norm{u}_{L_{r}((0, t), L_p(D)^n)} \le 1 \biggr\}
    \end{equation*}
    is precompact in $L_q(D)^n$, for any $\theta \in (0, t)$.  Now we will let $\theta$ approach $t$.  We estimate, applying Proposition~\ref{prop:system-estimates}(ii) and H\"older's inequality,
   \begin{equation}
   \label{eq:compact-aux}
       \begin{aligned}
            & \phantom{\le} \biggl\lVert \int_{\theta}^t U_{a_0}^0(t, {\zeta}) u({\zeta}) \, \dd{\zeta} \biggr\rVert_{L_q(D)^n} \le \int_{\theta}^t \lVert U_{a_0}^0(t, {\zeta}) u({\zeta}) \rVert_{L_q(D)^n} \, \dd{\zeta}
            \\
            & \le M e^{\gamma} \int_{\theta}^t (t - \zeta)^{- \frac{N}{2} (\frac{1}{p} - \frac{1}{q})} \norm{u(\zeta)}_{L_p(D)^n} \, \dd\zeta
            \\
            & \le M e^{\gamma} \biggl(  \int_{\theta}^t (t - \zeta)^{- \frac{N}{2} (\frac{1}{p} - \frac{1}{q}) r'} \, \dd\zeta \biggr)^{\!\! 1/r'} \biggl(  \int_{\theta}^t \norm{u(\zeta)}_{L_p(D)^n}^{r} \, \dd\zeta \biggr)^{\!\! 1/r}
            \\
            & \le \frac{M e^{\gamma}}{\bigl(1 - \frac{N}{2} (\frac{1}{p} - \frac{1}{q}) r'\bigr)^{1/r'}} (t - \theta)^{\frac{1}{r'} - \frac{N}{2} (\frac{1}{p} - \frac{1}{q})}.
       \end{aligned}
   \end{equation}
    Let $(a_{0,m})_{m=1}^{\infty} \subset Y_0$ and  $(u_{m})_{m=1}^{\infty} \subset L_{r}((0, t), L_p(D)^n)$, $\norm{u_m}_{L_{r}((0, t), L_p(D)^n)} \le 1$. For each $\epsilon > 0$ we can thus find $\theta \in (0, t)$ so close to $t$ that
   \begin{equation*}
       \biggl\lVert \int_{\theta}^t U_{a_{0,m}}^0(t, {\zeta}) u_m({\zeta}) \, \dd{\zeta} \biggr\rVert_{L_p(D)^n} < \epsilon
   \end{equation*}
   for all $m \in \NN$.  On the other hand, by passing to a subsequence we can assume that the sequence $\bigl(\int\limits_{0}^{\theta} U_{a_{0,m}}^0(t, {\zeta}) u_m({\zeta}) \, \dd{\zeta}\bigr)_{m = 1}^{\infty}$ is convergent in $L_q(D)^n$.  Consequently, we can find $m_0 \in \NN$ such that
    \begin{equation*}
        \biggl\lVert \int\limits_{0}^{\theta} U_{a_{0,m_1}}^0(t, \zeta) u_{m_1}(\zeta) \, \dd\zeta - \int\limits_{0}^{\theta} U_{a_{0,m_2}}^0(t, \zeta) u_{m_2}(\zeta) \, \dd\zeta \, \biggr\rVert_{L_q(D)^n} < \epsilon
    \end{equation*}
    for any $m_1, m_2 \ge m_0$. Therefore
    \begin{equation*}
        \biggl\lVert \int\limits_{0}^{t} U_{a_{0,m_1}}^0(t, \zeta) u_{m_1}(\zeta) \, \dd\zeta - \int\limits_{0}^{t} U_{a_{0,m_2}}^0(t, \zeta) u_{m_2}(\zeta) \, \dd\zeta \, \biggr\rVert_{L_q(D)^n} < 3 \epsilon
    \end{equation*}
    for any $m_1, m_2 \ge m_0$.  From this it follows that
    \begin{equation*}
        \biggl( \int\limits_{0}^{t} U_{a_{0,m}}^0(t, \zeta)  u_{m}(\zeta) \, \dd\zeta \biggr)_{\!\! m = 1}^{\!\! \infty}
    \end{equation*}
    is a Cauchy sequence in $L_q(D)^n$.
\end{proof}

\begin{lemma}
\label{lm:compact-aux-2}
    Assuming $0 < T \le 1$, for any $1 \le p \le q < \infty$, $1 < r \le \infty$ satisfying
    \begin{equation*}
        \frac{N}{2}\biggl(\frac{1}{p}- \frac{1}{q}\biggr) < \frac{1}{r'},
    \end{equation*}
    and any bounded $F \subset L_{r}((0, T), L_p(D)^n)$ the set
    \begin{equation*}
        \biggl\{ \Bigl[\, [0, T] \ni t \mapsto \int\limits_{0}^{t} U_{a_0}^0(t, {\zeta}) u({\zeta}) \, \mathrm{d}{\zeta} \in L_{q}(D) \,\Bigr]: a_0 \in Y_{0}, \ u \in F \biggr\}
    \end{equation*}
    is precompact in $C([0, T], L_{q}(D)^n)$.
\end{lemma}
\begin{proof}
    By Lemma~\ref{lm:I_i2}, the function
    \begin{equation*}
        \Bigl[\, [0, T] \ni t \mapsto \int_{0}^{t} U^0_{a_0}(t, \zeta) u(\zeta) \, \mathrm{d}\zeta \,\Bigr]
    \end{equation*}
    belongs to $C([0, T], L_q(D)^n)$, for any $u \in F$.

    By the Arzel\`a--Ascoli theorem, it suffices, taking Lemma~\ref{lm:compact-aux-at_once} into account, to show that for any $\epsilon > 0$ there is $\delta > 0$ such that, if $0 \le t_1 \le t_2 \le T$, $t_2 - t_1 < \delta$, then
    \begin{equation*}
        \biggl\lVert \int_{0}^{t_2} U_{a_0}^0(t_2, {\zeta}) u({\zeta}) \, \mathrm{d}{\zeta} - \int_{0}^{t_1} U_{a_0}^0(t_1, {\zeta}) u({\zeta}) \, \mathrm{d}{\zeta} \biggr\rVert_{L_{q}(D)^n} < \epsilon
    \end{equation*}
    for all $a_0 \in Y_{0}$ and all $u \in F$.  In order not to introduce too many constants we assume that $F$ equals the unit ball in $L_{r}((0, T), L_p(D)^n)$.

    Let $\epsilon > 0$ be fixed.  We write
    \begin{multline*}
        \int_{0}^{t_2} U_{a_0}^0(t_2, {\zeta}) u({\zeta}) \, \mathrm{d}{\zeta} - \int_{0}^{t_1} U_{a_0}^0(t_1, {\zeta}) u({\zeta}) \, \mathrm{d}{\zeta}
        \\
        = \int_{0}^{t_1} \bigl( U_{a_0}^0(t_2, {\zeta}) - U_{a_0}^0(t_1, {\zeta})\bigr) u({\zeta}) \, \mathrm{d}{\zeta} + \int_{t_1}^{t_2} U_{a_0}^0(t_2, {\zeta}) u({\zeta}) \, \mathrm{d}{\zeta}.
    \end{multline*}
    Repeating the reasoning as in~\eqref{eq:compact-aux} we see that we can find $\delta_1 > 0$ such that \begin{equation}
    \label{eq:compact-aux-1}
        \biggl\lVert \int_{t_1}^{t_2} U_{a_0}^0(t_2, {\zeta}) u({\zeta}) \, \mathrm{d}{\zeta} \biggr\rVert_{L_{q}(D)^n} < \frac{\epsilon}{3},
    \end{equation}
    provided $t_2 - t_1 < \delta_1$.

    Further, we write, for some $\eta \in (0, t_1)$,
    \begin{multline*}
        \int_{0}^{t_1} \bigl( U_{a_0}^0(t_2, {\zeta}) - U_{a_0}^0(t_1, {\zeta})\bigr) u({\zeta}) \, \mathrm{d}{\zeta}
        \\
        = \int_{0}^{t_1 - \eta} \bigl( U_{a_0}^0(t_2, {\zeta}) - U_{a_0}^0(t_1, {\zeta})\bigr) u({\zeta}) \, \mathrm{d}{\zeta} + \int_{t_1 - \eta}^{t_1} \bigl( U_{a_0}^0(t_2, {\zeta}) - U_{a_0}^0(t_1, {\zeta})\bigr) u({\zeta}) \, \mathrm{d}{\zeta}.
    \end{multline*}
    Again, copying the reasoning as in~\eqref{eq:compact-aux} we find $\eta > 0$ such that
    \begin{equation}
    \label{eq:compact-aux-2}
        \biggl\lVert \int_{t_1 - \eta}^{t_1} \bigl( U_{a_0}^0(t_2, {\zeta}) - U_{a_0}^0(t_1, {\zeta})\bigr) u({\zeta}) \, \mathrm{d}{\zeta} \biggr\rVert_{L_{q}(D)^n} < \frac{\epsilon}{3}.
    \end{equation}
    It follows from Proposition~\ref{prop:continuity-in-norm-top} that the assignment
    \begin{equation*}
        \bigl[\, \{\, (t, \zeta) \in \dot{\Delta} : t - \zeta \ge \eta \,\} \ni (t, \zeta) \mapsto U^0_{a_0}(t, \zeta) \in \mathcal{L}(L_p(D)^n, L_q(D)^n) \,\bigr]
    \end{equation*}
    is uniformly continuous, consequently there exists $\delta \in (0, \delta_1]$ such that if $\zeta + \eta \le t_1$ and $0 \le t_2 - t_1 < \delta$ then
    \begin{equation*}
        \lVert U^0_{a_0}(t_2, \zeta) - U^0_{a_0}(t_1, \zeta) \rVert_{\mathcal{L}(L_p(D)^n, L_q(D)^n)}  < \frac{\epsilon}{3 \eta^{1/r'}}.
    \end{equation*}
   Therefore, by applying H\"older's inequality,
    \begin{equation}
    \label{eq:compact-aux-3}
    \begin{aligned}
        & \phantom{\le} \biggl\lVert \int_{0}^{t_1 - \eta} \bigl( U_{a_0}^0(t_2, {\zeta}) - U_{a_0}^0(t_1, {\zeta})\bigr) u({\zeta}) \, \mathrm{d}{\zeta} \biggr\rVert_{L_{q}(D)^n}
        \\
        & < \frac{\epsilon}{3 \eta^{1/r'}} \int_{0}^{t_1 - \eta} \norm{u(\zeta)}_{L_p(D)^n} \, \dd\zeta \\
        & \le \frac{\epsilon}{3 \eta^{1/r'}} \biggl( \int_{0}^{t_1 - \eta} \norm{u(\zeta)}_{L_p(D)^n}^r \, \dd\zeta \biggr)^{\!\! 1/r} \eta^{1/r'}
        \le \frac{\epsilon}{3 } \norm{u}_{L_r((0,T), L_p(D)^n)}
        \\
        & \le \frac{\epsilon}{3}.
    \end{aligned}
    \end{equation}
    The estimates~\eqref{eq:compact-aux-1}, \eqref{eq:compact-aux-2} and~\eqref{eq:compact-aux-3} do not depend on the choice of $a_0 \in Y_0$, so gathering them gives the required property.
\end{proof}

\begin{proposition}
\label{prop:delay-compactness-X_1-at_once}
    Assuming $T \le 1$, for any $0 < t_0 \le T$, $1 \le p \le q < \infty$, $1 < r \le \infty$ satisfying
    \begin{equation*}
        \frac{N}{2}\biggl(\frac{1}{p}- \frac{1}{q}\biggr) < \frac{1}{r'},
    \end{equation*}
    and any bounded $E \subset L_{r}((0, T), L_p(D)^n)$ the set
    \begin{equation*}
        \Bigl\{\, \bigl[\, [t_0, T] \ni t \mapsto u(t; 0, a, u_0) \,\bigr] : a \in Y, u_0 \in E \,\Bigr\}
    \end{equation*}
  is precompact in $C([t_0, T], L_q(D)^n)$.
\end{proposition}
\begin{proof}
    We will use the notation $I_i(t; a, u_0)$, $i = 0, 1, 2$, as in~\eqref{eq:I}, taking account of the parameter $a$ and the initial value $u_0$.  The precompactness of the set
    \begin{equation*}
        \Bigl\{\, \bigl[\, [t_0, T] \ni t \mapsto I_0(t; a, u_0) \,\bigr] : a \in Y, u_0 \in E \,\Bigr\}
    \end{equation*}
    in $C([t_0, T], L_q(D)^n)$ is a consequence of Proposition~\ref{prop:higher_order-compactness}.

    To prove the precompactness of
    \begin{equation*}
        \Bigl\{\, \bigl[\, [0, T] \ni t \mapsto I_i(t; a, u_0) \,\bigr] : a \in Y, u_0 \in E \,\Bigr\}, \quad i = 1, 2,
    \end{equation*}
    in $C([0, T], L_p(D)^n)$ we apply Lemma~\ref{lm:compact-aux-2}.  In case of $I_1$ we take $F$ to be equal to $\{\, \mathcal{C}^{0}_{a}(\cdot) \, u(\cdot; 0, a, u_0) : a \in Y, u_0 \in E \,\}$, which is bounded in $L_{\infty}((0, T), L_p(D)^n)$ in view of the definition of solution, Proposition~\ref{prop:estimate-mild} and Lemma~\ref{lm:mult-estim-1}.  In case of $I_2$ for $F$ we take $\{\, \mathcal{C}^{1}_{a}(\cdot) \, u_0^{(2)} : a \in Y, u_0 \in E \,\}$,
    which is bounded in $L_r((0, T), L_q(D)^n)$ in view of Lemma~\ref{lm:mult-estim-1}.

    We conclude the proof by noting that the sum-set of three precompact sets is precompact.
\end{proof}

\subsection{Continuous dependence}
Assume moreover~\ref{AS-ae-converge}.

\label{subsubsect:continuous-dependence}
\begin{proposition}
\label{prop:continuity-wrt-parameters-1}
    Assuming $T \le 1$, for any $0 < t_0 \le T$, $1 \le p \le q < \infty$, $1 < r \le \infty$ satisfying
    \begin{equation}
    \label{eq:ass-p-q-2}
        \frac{N}{2}\biggl(\frac{1}{p}- \frac{1}{q}\biggr) < \frac{1}{r'},
    \end{equation}
    and any $u_0 \in L_p(D)^n \times L_r((-1,0), L_p(D)^n)$ the mapping
    \begin{equation*}
        \bigl[\, Y \ni a \mapsto u(\cdot; 0, a, u_0){\restriction}_{[t_0, T]}  \in C([t_0, T], L_{q}(D)^n) \,\bigr]
    \end{equation*}
    is continuous.

    If $1 < p = q$ then the mapping
    \begin{equation*}
        \bigl[\, Y \ni a \mapsto u(\cdot; 0, a, u_0)  \in C([0, T], L_{p}(D)^n) \,\bigr]
    \end{equation*}
    is continuous.
\end{proposition}
\begin{proof}
     We assume first that $1 \le p < q < \infty$.  Fix $p$, $q$ and $u_0$.  Let $(a_{m})_{m=1}^{\infty} \subset Y$ converge to $a$, and put $u_m(\cdot)$ for $u(\cdot; 0, a_m, u_0)$ and $u(\cdot)$ for $u(\cdot; 0, a, u_0)$.  We can take $a_{0,m}$ to converge to $a_0$ in $Y_0$.

    It follows from Proposition~\ref{prop:delay-compactness-X_1-at_once} via diagonal process that we can assume, by taking a subsequence if~necessary, that $u_{m}{\restriction}_{(0, T]}$ converge to some continuous $\tilde{u} \colon (0, T] \to L_{q}(D)^n$ and the convergence is uniform on compact subsets of $(0, T]$.

    In view of assumption~\eqref{eq:ass-p-q-2}, Proposition~\ref{prop:delay-estimates-1} implies, via the Dominated Convergence Theorem for Bochner integrals (see~\cite[Thm.~X.3.12]{AmEsch}), that $u_{m}$ converge to $\tilde{u}$ in the $L_{r'}((0, T), L_q(D)^n)$\nobreakdash-\hspace{0pt}norm.

    We proceed to showing that, for each $t \in (0, T]$,
    \begin{align}
    \label{cont-1}
        U^0_{a_{0,m}}(t, 0) u_{0}^{(1)} & \to U^0_{a_{0}}(t, 0) u_{0}^{(1)}
        \\
    \label{cont-2}
        \int\limits_{0}^{t} U^0_{a_{0,m}}(t, \zeta) \bigl( \mathcal{C}^0_{a_{m}}(\zeta) u_{m}(\zeta) \bigr) \, \mathrm{d}\zeta & \to \int\limits_{0}^{t} U^0_{a_{0}}(t, \zeta) \bigl( \mathcal{C}^0_{a}(\zeta) \tilde{u}(\zeta) \bigr) \, \mathrm{d}\zeta
        \\
    \label{cont-3}
        \int\limits_{0}^{t} U^0_{a_{0,m}}(t, \zeta) \bigl( \mathcal{C}^1_{a_{m}}(\zeta) u_{0}^{(2)}(\zeta - 1) \bigr) \, \mathrm{d}\zeta & \to \int\limits_{0}^{t} U^0_{a_{0}}(t, \zeta) \bigl( \mathcal{C}^1_{a}(\zeta) u_{0}^{(2)}(\zeta - 1) \bigr) \, \mathrm{d}\zeta
    \end{align}
    in the $L_{q}(D)^n$\nobreakdash-\hspace{0pt}norm as $m \to \infty$.  Until revoking, let $t \in (0, T]$ be fixed.

    \eqref{cont-1} is a consequence of Proposition~\ref{prop:norm_continuity-all}.

    To prove~\eqref{cont-2}, write
    \begin{equation}
    \label{cont-4}
        \begin{aligned}
            & \int\limits_{0}^{t} U^0_{a_{0,m}}(t, \zeta) \bigl( \mathcal{C}^0_{a_{m}}(\zeta) u_{m}(\zeta) \bigr) \, {\mathrm{d}}\zeta - \int\limits_{0}^{t} U^0_{a_{0}}(t, \zeta) \bigl( \mathcal{C}^0_{a}(\zeta) \tilde{u}(\zeta) \bigr) \, \mathrm{d}\zeta
            \\
            = & \int\limits_{0}^{t} \bigl( U^0_{a_{0,m}}(t, \zeta) - U^0_{a_{0}}(t, \zeta) \bigr) \bigl( \mathcal{C}^0_{a_{m}}(\zeta) u_{m}(\zeta) \bigr) \, \mathrm{d}\zeta
            \\
            & + \int\limits_{0}^{t} U^0_{a_{0}}(t, \zeta) \bigl( \mathcal{C}^0_{a_{m}}(\zeta) (u_{m}(\zeta) - \tilde{u}(\zeta)) \bigr) \, \mathrm{d}\zeta
            \\
            & + \int\limits_{0}^{t} U^0_{a_{0}}(t, \zeta) \bigl( (\mathcal{C}^0_{a_{m}}(\zeta) - \mathcal{C}^0_{a}(\zeta)) \tilde{u}(\zeta) \bigr) \, \mathrm{d}\zeta.
        \end{aligned}
    \end{equation}
    Write $J^{(i)}_{m}(t)$, $i = 1, 2, 3$, for the $i$\nobreakdash-\hspace{0pt}th term on the right-hand side of~\eqref{cont-4}.

    As, by Proposition~\ref{prop:system-estimates}(i) and Lemma~\ref{lm:mult-estim}, $U^0_{a_{0}}(t, \zeta) \in \mathcal{L}(L_p(D)^n, L_q(D)^n)$ and $\mathcal{C}^0_{a_{m}}(\zeta) \in \mathcal{L}(L_p(D)^n)$ with norms bounded uniformly in $m$, we have that for Lebesgue\nobreakdash-\hspace{0pt}a.e. $\zeta \in (0, t)$ the integrands in $J_m^{(2)}(t)$ converge to zero in $L_q(D)^n$.  Further, by Proposition~\ref{prop:system-estimates}(ii), Lemma~\ref{lm:mult-estim} and Proposition~\ref{prop:delay-estimates-1},
    \begin{equation*}
        \norm{U^0_{a_{0}}(t, \zeta) \bigl( \mathcal{C}^0_{a_{m}}(\zeta) (u_{m}(\zeta) - \tilde{u}(\zeta)) \bigr)}_{L_q(D)^n} \le C (t - \zeta)^{- \frac{N}{2}(\frac{1}{p} - \frac{1}{q})} \zeta^{- \frac{N}{2}(\frac{1}{p} - \frac{1}{q})},
    \end{equation*}
    where $C > 0$ is independent of $m$.  In view of our assumptions on $p$ and $q$, the Dominated Convergence Theorem for Bochner integrals provides the convergence of $J^{(2)}_m(t)$ to $0$ in the $L_q(D)^n$\nobreakdash-\hspace{0pt}norm as $m \to \infty$.

    Regarding the convergence of $J^{(1)}_{m}(t)$ and $J^{(3)}_{m}(t)$ to $0$, we need to use a slightly subtler approach.  It follows from Lemma~\ref{lm:compact-aux-2} that the sets
    \begin{equation*}
        \{\, J^{(1)}_{m}(t) : m \in \NN \,\} \text{ and } \{\, J^{(3)}_{m}(t) : m \in \NN \,\}
    \end{equation*}
    are precompact in $L_{q}(D)^n$.  We get the required norm convergence if we prove that
    \begin{align}
    \label{cont-3-1}
        \langle J^{(1)}_{m}(t), v \rangle & \to 0 \text{ as } m \to \infty,
    \intertext{and}
    \label{cont-3-3}
        \langle J^{(3)}_{m}(t), v \rangle & \to 0 \text{ as } m \to \infty,
    \end{align}
    for any $v \in L_{{q'}}(D)^n$, where $\langle \cdot, \cdot \rangle$ stands for the duality pairing between $L_{q}(D)^n$ and $L_{q'}(D)^n$.  By the Hille theorem (\cite[Thm.~II.2.6]{DiUhl}) and Proposition~\ref{prop:dual-Lp-q},
    \begin{multline*}
        \langle J^{(1)}_{m}(t), v \rangle = \int_{0}^{t} \bigl\langle \bigl(U_{a_{0,m}}^{0}(t, \zeta) - U_{a_0}^{0}(t, \zeta)\bigr) \bigl( \mathcal{C}^0_{a_{m}}(\zeta) u_{m}(\zeta) \bigr), v \bigr\rangle \, \mathrm{d}\zeta
        \\
        = \int_{0}^{t} \bigl\langle  \mathcal{C}^0_{a_{m}}(\zeta) u_{m}(\zeta), \bigl(U_{a_{0,m}}^{0,*}(\zeta, t) - U_{a_{0}}^{0,*}(\zeta, t)\bigr) v \bigr\rangle \, \mathrm{d}\zeta
    \end{multline*}
    and
    \begin{multline*}
        \langle J^{(3)}_{m}(t), v \rangle =  \int_{0}^{t} \bigl\langle U_{a_0}^{0}(t, \zeta) \bigl( (\mathcal{C}^0_{a_{m}}(\zeta) - \mathcal{C}^0_{a}(\zeta)) \tilde{u}(\zeta) \bigr), v \bigr\rangle \, \mathrm{d}\zeta
        \\
        = \int_{0}^{t} \bigl\langle  \bigl(\mathcal{C}^0_{a_{m}}(\zeta) - \mathcal{C}^0_{a}(\zeta)\bigr) \tilde{u}(\zeta), U_{a_{0}}^{0,*}(\zeta, t) v \bigr\rangle \, \mathrm{d}\zeta.
    \end{multline*}
    The integrands $\bigl\langle  \mathcal{C}^0_{a_{m}}(\zeta) u_{m}(\zeta), \bigl(U_{a_{0,m}}^{0,*}(\zeta, t) - U_{a_{0}}^{0,*}(\zeta, t)\bigr) v \bigr\rangle$ have, by Lemma~\ref{lm:mult-estim}, Proposition~\ref{prop:delay-estimates-1} and Proposition~{\ref{prop:higher_orders-skew_product-p-q-adjoint}}, their $L_1((0, t))$\nobreakdash-\hspace{0pt}norms bounded uniformly in $m$, and the adjoint equation analog of Proposition~\ref{prop:higher_orders-skew_product-continuity} gives that $U_{a_{0,m}}^{0,*}(\zeta, t) v \to U_{a_{0}}^{0,*}(\zeta, t) v$ in $L_{q'}(D)^n$ for all $\zeta \in (0, t)$, consequently, by the Dominated Convergence Theorem for Lebesgue integrals on $[0, t]$, $\langle J^{(1)}_{m}(t), v \rangle$ converges to $0$ as $m \to \infty$.

    We write
    \begin{align*}
        & \int_{0}^{t} \bigl\langle  \bigl(\mathcal{C}^0_{a_{m}}(\zeta) - \mathcal{C}^0_{a}(\zeta)\bigr) \tilde{u}(\zeta), U_{a_{0}}^{0,*}(\zeta, t) v \bigr\rangle \, \mathrm{d}\zeta
        \\
        & = \int_{0}^{t} \biggl( \sum_{k, l = 1}^n \int_{D} (c_{0,m}^{kl}(\zeta, x) - c_0^{kl}(\zeta, x)) \, \tilde{u}^k(\zeta)[x] \, (U_{a_{0}}^{l,*}(\zeta, t) v^l)[x] \, \mathrm{d}x \biggr)  \mathrm{d}\zeta
    \end{align*}
    Fix, for a moment, $1 \le k, l \le n$.

    We have that $\tilde{u}^k  \in L_{r'}((0, t), L_q(D)) \subset L_1((0, t), L_q(D))$.  By Proposition~\ref{prop:higher_orders-skew_product-p-adjoint}(ii), the function $\bigl[ \, [0, t] \ni \zeta \mapsto U_{a_{0}}^{l,*}(\zeta, t) v^l \in L_{q'}(D) \, \bigr]$ belongs to $C([0, t], L_{q'}(D))$.  Therefore, $\tilde{u}^k(\zeta) (U_{a_{0}}^{l,*}(\zeta, t) v^l) \in L_1(D)$ for Lebesgue\nobreakdash-\hspace{0pt}a.e.\ $\zeta \in (0, t)$, with the function  $\bigl[ \, [0, t] \ni \zeta \mapsto \norm{\tilde{u}^k(\zeta) (U_{a_{0}}^{l,*}(\zeta, t) v^l)}_{L_1(D)} \, \bigr]$ in $L_1((0, t))$.  Further, since $L_q(D)$ and $L_{q'}(D)$ embed continuously in $L_1(D)$, we have that both $\tilde{u}^k$ and $\bigl[ \, [0, t] \ni \zeta \mapsto U_{a_{0}}^{l,*}(\zeta, t) v^l \in L_{1}(D) \, \bigr]$ are in $L_1((0, t), L_1(D))$.  By Lemma~\ref{lm:Dunford-Schwartz}(a), both
    \begin{equation*}
         \bigl[ \, (0, t) \times D \ni (\zeta, x) \mapsto \tilde{u}^k(\zeta)[x] \in \RR \, \bigr]
    \end{equation*}
    and
    \begin{equation*}
         \bigl[ \, (0, t) \times D \ni (\zeta, x) \mapsto (U_{a_0}^{l,*}(\zeta, t) v^l)[x] \in \RR \, \bigr]
    \end{equation*}
    belong to $L_1((0, t) \times D)$, hence are $(\mathfrak{L}((0, t) \times D), \mathfrak{B}(\RR))$\nobreakdash-\hspace{0pt}measurable. Consequently their product
     \begin{equation}
     \label{eq:aux1}
        \bigl[ \, (0, t) \times D \ni (\zeta, x) \mapsto \tilde{u}^k(\zeta)[x] \, (U_{a_0}^{l,*}(\zeta, t) v^l)[x] \in \RR \, \bigr]
    \end{equation}
    is $(\mathfrak{L}((0, t) \times D), \mathfrak{B}(\RR))$\nobreakdash-\hspace{0pt}measurable.  As we have shown earlier that its $\zeta$\nobreakdash-\hspace{0pt}sections are Lebesgue\nobreakdash-\hspace{0pt}a.e.\ in $L_1(D)$, Lemma~\ref{lm:Dunford-Schwartz}(b)  gives that $\tilde{u}^k(\cdot) (U_{a_0}^{l,*}(\cdot, t) v^l) \in L_1((0, t), L_1(D))$, so, again by Lemma~\ref{lm:Dunford-Schwartz}(a), \eqref{eq:aux1} belongs to $L_1((0, t) \times D)$.
    Putting it to be zero on $(t, 1) \times D$ we have that the extended function belongs to $L_1((0, 1) \times D)$.  Since $c_{0,m}^{kl}$ converge, as $m \to \infty$, to $c_{0}^{kl}$ in the weak-$*$ topology of $(L_1((0, t) \times D))^{*}$, the convergence of $\langle J_{m}^{(3)}(t), v \rangle$ to zero follows.

    For proving~\eqref{cont-3}, we write
    \begin{equation}
    \label{cont-5}
        \begin{aligned}
            & \int\limits_{0}^{t} U_{a_{0,m}}^0(t, \zeta) \bigl( \mathcal{C}^1_{a_{m}}(\zeta) u_{0}^{(2)}(\zeta - 1) \bigr) \, \mathrm{d}\zeta - \int\limits_{0}^{t} U_{a_{0}}^0(t, \zeta) \bigl(\mathcal{C}^1_{a}(\zeta) u_{0}^{(2)}(\zeta - 1) \bigr) \, \mathrm{d}\zeta
            \\
            = & \int\limits_{0}^{t} (U_{a_{0,m}}^0(t, \zeta) - U_{a_{0}}^0(t, \zeta)) \bigl( C^1_{a_{m}}(\zeta) u_{0}^{(2)}(\zeta - 1) \bigr) \, \mathrm{d}\zeta
            \\
            & + \int\limits_{0}^{t} U_{a_{0}}^0(t, \zeta) \bigl((\mathcal{C}^1_{a_{m}}(\zeta) - \mathcal{C}^1_{a}(\zeta)) u_{0}^{(2)}(\zeta - 1) \bigr)\, \mathrm{d}\zeta
        \end{aligned}
    \end{equation}
    and proceed as in the proof for $J_{m}^{(1)}(t)$ (in the case of the first summand on the right-hand side) and in the proof for $J_{m}^{(3)}(t)$ (in the case of the second summand).

    We have thus proved that for any accumulation point $\tilde{u}$ of $(u_{m})_{m=1}^{\infty}$ there holds
    \begin{equation}
    \label{eq:auxiliary}
        \begin{aligned}
            \tilde{u}(t) = {} & U^0_{a_{0}}(t, 0) u_{0}^{(1)} + \int\limits_{0}^{t} U^0_{a_{0}}(t, \zeta) \bigl( \mathcal{C}^0_{a}(\zeta) \tilde{u}(\zeta) \bigr) \, \mathrm{d}\zeta
        \\
        &  + \int\limits_{0}^{t} U^0_{a_{0}}(t, \zeta) \bigl( \mathcal{C}^1_{a}(\zeta) u_{0}^{(2)}(\zeta - 1) \bigr) \, \mathrm{d}\zeta, \quad t \in (0, T]
        \end{aligned}
     \end{equation}
    in the $L_q(D)^n$\nobreakdash-\hspace{0pt}sense, hence in the $L_p(D)^n$\nobreakdash-\hspace{0pt}sense.  We extend $\tilde{u}$ to $[0, T]$ by putting $\tilde{u}(0) = u_{0}^{(1)}$.  We shall prove that $\tilde{u}$ so extended is an $L_p$\nobreakdash-\hspace{0pt}solution.  By Proposition~\ref{prop:continuiuty-to-C(ii)}, the first summand on the right\nobreakdash-\hspace{0pt}hand side of \eqref{eq:auxiliary} belongs to \newline $C([0, T], L_p(D)^n)$, and, by Lemma~\ref{lm:I_i2}, the second and third summands do.   This proves that $\tilde{u}$ is indeed a mild $L_p$\nobreakdash-\hspace{0pt}solution of \eqref{main-eq}$_{a}$\textup{+}\eqref{main-bc}$_{a}$ satisfying the initial condition $\tilde{u}(0) = u_{0}^{(1)}$.  By uniqueness, not only some subsequence but the whole sequence $u(\cdot; 0, a_m, u_0)$ must converge to $\tilde{u}$ in the $C([0, T], L_p(D)^n)$\nobreakdash-\hspace{0pt}norm.  Going back to $L_q(D)^n$, we have, by Proposition~\ref{prop:continuiuty-to-C(i)}, that the first summand in~\eqref{eq:auxiliary} belongs to $C([t_0, T], L_q(D)^n)$ for any $t_0$, whereas the second and third summands belong to $C([0, T], L_q(D)^n)$, by Lemma~\ref{lm:I_i2}.

    The above reasoning holds true in the case $1 < p = q$.  When $1 = p = q$, we take $\tilde{q} > 1$ but so close to $1$ that \eqref{eq:ass-p-q-2} is satisfied.  Then we apply the previous reasoning to prove that, for each $t_0 \in (0, T]$ the mapping
    \begin{equation*}
        \bigl[\, Y \ni a \mapsto u(\cdot; 0, a, u_0){\restriction}_{[t_0, T]}  \in C([t_0, T], L_{q'}(D)^n) \,\bigr]
    \end{equation*}
    is continuous, hence the mapping
    \begin{equation*}
        \bigl[\, Y \ni a \mapsto u(\cdot; 0, a, u_0){\restriction}_{[t_0, T]}  \in C([t_0, T], L_{1}(D)^n) \,\bigr]
    \end{equation*}
    is continuous, too.
\end{proof}

\begin{lemma}
\label{lm:continuity-wrt-parameters-2}
    Assuming $T \le 1$, for any $0 < t_0 \le T$, $1 \le p \le q < \infty$, $1 < r \le \infty$ satisfying
    \begin{equation*}
        \frac{N}{2}\biggl(\frac{1}{p}- \frac{1}{q}\biggr) < \frac{1}{r'},
    \end{equation*}
    the mapping
    \begin{multline*}
        \Bigl[\, Y \times \bigl(L_p(D)^n \times L_r((-1,0), L_p(D)^n)\bigr) \ni (a, u_0)
        \\
        \mapsto u(\cdot; 0, a, u_0){\restriction}_{[t_0, T]} \in C([t_0, T], L_{q}(D)^n) \,\Bigr]
    \end{multline*}
    is continuous.
\end{lemma}
\begin{proof}
    Fix $p$, $q$ and $t_0$.  Let $(a_{m})_{m=1}^{\infty}$ converge to $a$ in $Y$ and let $(u_{0,m})_{m=1}^{\infty}$ converge to $u_0$ in $C([- 1, 0], L_p(D)^n)$.  Put $u_m(\cdot)$ for $u(\cdot; 0, a_m, u_{0,m})$ and $u(\cdot)$ for $u(\cdot; 0, a, u_0)$. For each $m \in \NN$ we have
    \begin{align*}
        & \norm{(u_m - u) {\restriction}_{[t_0, T]} }_{C([t_0, T], L_{q}(D)^n)}
        \\
        & \le \norm{(u_m - u(\cdot; 0, a_m, u_0)) {\restriction}_{[t_0, T]} }_{C([t_0, T], L_{q}(D)^n)}
        \\
        & \phantom{\le} + \norm{(u(\cdot; 0, a_m, u_0) - u) {\restriction}_{[t_0, T]}}_{C([t_0, T], L_{q}(D)^n)}.
    \end{align*}
    As $m \to \infty$, the first term on the right-hand side converges to zero by Proposition~\ref{prop:delay-estimates-1} and the second term converges to zero by Proposition~\ref{prop:continuity-wrt-parameters-1}.
\end{proof}
In view of Lemma~\ref{lm:solution-translate} we have the following.
\begin{proposition}
\label{prop:continuity-wrt-parameters}
    For any $1 \le p \le q < \infty$, $1 < r \le \infty$ with
     \begin{equation*}
        \frac{N}{2}\biggl(\frac{1}{p} - \frac{1}{q}\biggr) < \frac{1}{r'},
    \end{equation*}
    and $0 < t_0 \le T$ the mapping
    \begin{multline*}
        \Bigl[\, Y \times \bigl( L_p(D)^n \times L_r((-1,0), L_p(D)^n) \bigr) \ni (a, u_0)
        \\
        \mapsto u(\cdot; 0, a, u_0){\restriction}_{[t_0, T]} \in C([t_0, T], L_{q}(D)^n) \,\Bigr]
    \end{multline*}
    is continuous.
\end{proposition}

As stated above Proposition~\ref{prop:regularization-q}, for any $1 < p < q < \infty$ we can take $p = p_0 < p_1 < \ldots < p_{m_0 - 1} < p_{m_0} = q$ such that
\begin{equation*}
    \frac{N}{2} \biggl( \frac{1}{p_m} - \frac{1}{p_{m-1}} \biggr) < \frac{1}{2 r'}, \quad m = 0, \dots, m_0 - 1,
\end{equation*}
so applying Proposition~\ref{prop:continuity-wrt-parameters} finitely many times we obtain the following.
\begin{theorem}
\label{thm:continuity-wrt-parameters-regularizing}
    Assuming $T > 2$, there exists $\widehat{T} \in (2, T)$ such that for any $1 < p \le q < \infty$ the mapping
    \begin{multline*}
        \Bigl[\, Y \times \bigl(L_p(D)^n \times L_r((-1,0), L_p(D)^n)\bigr) \ni (a, u_0)
        \\
        \mapsto u(\cdot; 0, a, u_0){\restriction}_{[\widehat{T}, T]} \in C([\widehat{T}, T], L_q(D)^n) \,\Bigr]
    \end{multline*}
    is continuous.
\end{theorem}

\subsection*{Acknowledgments}
The authors would like to thank the anonymous referee for their careful reading the manuscript and for valuable comments and suggestions.

M. K. was supported by the National Science Centre, Poland (NCN) under the grant Sonata Bis with a number NCN 2020/38/E/ST1/00153.  J. M. was partially supported by the Faculty of Pure and Applied Mathematics, Wrocław University of Science and Technology.

\subsection*{Data availability}
No data was used for the research described in the article.

\end{document}